\documentclass[12pt, a4paper]{article}
\usepackage[margin=1in]{geometry}
\usepackage{amsfonts}
\usepackage{setspace}
\usepackage[utf8]{inputenc}
\usepackage{enumitem}
\usepackage{biblatex, ulem}
\usepackage{graphicx,subcaption} 

\usepackage{soul}
\usepackage{amsmath, amsthm, amssymb, mathtools, mathabx, amscd, graphicx, color}
\usepackage{accents}
\usepackage{bbm}
\usepackage{parskip}
\usepackage{textcomp}
\usepackage{comment}
\usepackage[colorlinks=true,urlcolor=blue,citecolor=blue,linkcolor=blue]{hyperref}

\usepackage{fancyhdr}
\newtheorem{theorem}{Theorem}[section]
\newtheorem{lemma}[theorem]{Lemma}

\newtheorem{corollary}[theorem]{Corollary}
\theoremstyle{definition}
\newtheorem{definition}[theorem]{Definition}

\newtheorem{claim}[theorem]{Claim}
\theoremstyle{remark}
\newtheorem{remark}[theorem]{Remark}
\numberwithin{equation}{section}
\newcommand{\ubar}[1]{\underline{#1}}
\newcommand{\absmod}[1]{\left|#1\right|}

\title{CLT for Linear Spectral Statistics of sample covariances of superimposition of Separable Covariance Components}

\author{Javed Hazarika and Debashis Paul}

\usepackage{setspace}
\begin{document}
\maketitle

\begin{abstract}
Suppose $\{Z_r\}_{r=1}^K$ are $p\times n$ random matrices with entries having zero mean and unit variance and $p/n \to c \in (0,\infty)$. We construct the sample covariance matrix $S_n = \frac{1}{n}XX^*$ where $X =\sum_{r=1}^K X_r$, where $X_r = A_r^\frac{1}{2}Z_rB_r^\frac{1}{2}$, $K \in \mathbb{N}$ and $A_r,B_r$ are sequences of positive semi-definite matrices of dimensions $p\times p$ and $n\times n$, respectively that are are simultaneously digaonalizable within themselves. As an extension of our previous work related to the L.S.D. of such classes of matrices, we show that the linear spectral statistics arising from the spectra of $S_n$ have Gaussian limits under assumptions of fourth moment bounds on the entries of $\{Z_r\}_{r=1}^K$ and convergence of the joint spectral distribution of $\{A_r\}_{r=1}^K$ and $\{B_r\}_{r=1}^K$ to probability distributions over $\mathbb{R}_+^K$ with compact support.
\end{abstract}

\section{Introduction}

Linear spectral statistics (LSS) of large random matrices play a central role in modern multivariate statistics, high-dimensional inference and random matrix theory. Given a $p\times p$ random matrix with eigenvalues $\lambda_1,\ldots,\lambda_p$, an LSS is a statistic of the form
$$L_p(f)=\sum_{i=1}^p f(\lambda_i)$$
where $f$ is a sufficiently regular test function. The asymptotic behavior of such statistics provides considerably finer information than the empirical spectral distribution alone. While the latter describes the first-order spectral behavior, a central limit theorem (CLT) for LSS characterizes the fluctuations around the deterministic limiting spectral distribution and forms the basis for a number of second-order statistical procedures.

There is now an extensive body of literature focusing on the limiting behavior of the spectra of sample covariance matrices corresponding to data matrices with a separable covariance structure under high dimensional asymptotic regimes. The earliest comprehensive work, in particular establishing the limiting spectral distribution (LSD) of sample covariances, was by \cite{Lixin06}, which has been followed by \cite{PaulSilverstein2009},  \cite{SepCovar}, who covered various aspects of the limiting spectrum. More recently, 
CLTs for linear spectral statistics of sample covariances
have been established under a separable covariance model \cite{BaiLiPan2019}. Similar work has also been done in \cite{XieJohnstone} for random effects model with applications in genetics. This was preceded by work done by \cite{BaiSilv04} and \cite{NajimYao} which established this line of work for large dimensional sample covariance matrices. Moreover, \cite{WangPanJiang} have studies the asymptotic spectral behavior of high-dimensional random
Gram matrices with sparsity and a variance profile, motivated by applications in wireless communications.

More recently, models with a more general variance-profile or separable structure have attracted considerable attention. In these models, the population structure cannot necessarily be described by a single covariance matrix on each side. Instead, the interaction between the two dimensions may involve several population components. A convenient formulation as done in \cite{XieJohnstone} is the data model where the $j-$th observation of the data matrix X has a variance given by
$$T_{j}=\sum_{r=1}^K l_{jr}^2\Sigma_r,$$
where $K$ is a fixed integer, $l_{jr}$'s are loadings and the $\Sigma_r$'s represent the $K$ levels of variation. This framework allows for substantial heterogeneity in the population structure.

In this manuscript, we extend the result to sample covariances of data matrices arising out of a mixture several signals each having a separable covariance structure. This will build upon the LSD type results of such classes of matrices done in \cite{JHDP26}. The matrix class we are working with is expressed as follows:
$$S_n = \frac{1}{n}XX^*;  \text{ where } X = \sum_{r=1}^KX_r, \quad X_r = A_r^\frac{1}{2}Z_rB_r^\frac{1}{2}.$$
Here, $Z_r \in \mathbb{C}^{p\times n}$ has independent entries with zero mean and unit variance and the row and column scaling matrices $\{A_r \in \mathbb{C}^{p \times  p}\}, \{B_r \in \mathbb{C}^{n \times n}\}$ are p.s.d. matrices satisfying some mild conditions.

Such models arise naturally when observations exhibit dependence across both rows and columns. The limiting spectral distribution of separable covariance matrices has been extensively studied, and CLTs for their LSS have been established under increasingly general assumptions. In particular, the work of \cite{BaiSilv04}, \cite{BaiLiPan2019} and subsequent developments provide a powerful framework based on resolvent methods, martingale decompositions, and deterministic equivalents.

Primary goal of this paper is to establish a CLT for linear spectral statistics in this generalized setting. The main result identifies both the limiting mean and covariance of the centered and scaled LSS and shows that the fluctuations converge to a Gaussian process. The resulting formulas are expressed in terms of the deterministic quantities associated with the generalized self-consistent equations governing the limiting resolvent. In particular, the covariance structure is characterized through a finite-dimensional system involving the interaction of the (K) population components.

We point out that a major difficulty in the extension from the classical separable model (\cite{BaiLiPan2019}) to our current model lies in the fact that the deterministic equivalent is no longer described by a single scalar self-consistent equation. Instead, one obtains a finite dimensional system of coupled equations. When two spectral parameters $z_1,z_2\in\mathbb{C}\backslash\mathbb{R}$ are considered simultaneously, the corresponding covariance calculations lead naturally to matrix-valued quantities such as
$$U_{2,2}(z_1,z_2),\qquad W_2(z_1,z_2),$$
and to a finite-dimensional stability operator of the form
$$I-c\,U_{2,2}(z_1,z_2)W_2(z_1,z_2).$$
Controlling this operator \textit{uniformly} over the spectral domain is one of the principal technical ingredients of the proof. The existence and uniqueness of the underlying deterministic solution, as well as the stability of the associated Dyson map, provide the mechanism through which these matrix-valued quantities can be controlled.

Our proof follows the general resolvent and martingale approach developed for CLTs of LSS, but several steps require substantial modification in the generalized model. First, the deterministic equivalent is obtained from a coupled system of (K)-dimensional equations rather than a scalar equation. We establish the required stability properties of this system and use them to obtain uniform bounds for the associated resolvents and their derivatives. Second, the covariance calculation produces products of matrix-valued deterministic equivalents evaluated at two spectral parameters. A key step is therefore to control the corresponding stability matrix and to establish the invertibility of the operator appearing in the covariance equations. Third, the bias term requires a careful finite-dimensional linearization of the deterministic equations. This leads to a system whose solution determines the limiting mean of the LSS.

The resulting theory also clarifies the relationship between the different levels of generality considered in the literature. When $K=1$, the coupled system reduces to the familiar scalar deterministic-equivalent equations for separable covariance matrices, and the corresponding CLT recovers the classical structure. Thus, the present framework is a finite-dimensional generalization of the standard separable covariance theory.

The rest of the paper is organized as follows. We start with the setup and results related to the LSD (from \cite{JHDP26}) in Section \ref{sec:ModelPreliminaries}. Section \ref{sec:MainResult} contains Theorem \ref{mainTheorem_CLT} which is the main result of this manuscript. We provide a detailed proof sketch in Section \ref{sec:proofSketch}. Assuming Gaussian innovations, Section \ref{sec:GaussianCase} details the process of arriving at the in-probability limit of the covariance function of the underlying random process. In Section \ref{sec:BiasTerm}, we derive the limit of the bias component. Sections \ref{sec:Truncation_Centralization_rescaling},  \ref{sec:NonGaussian} and \ref{sec:Tightness} cover respectively, the proof of truncation, centralization and scaling of innovations, the extension to the more general non-Gaussian case and the proof of tightness.

\section{Model and preliminaries}\label{sec:ModelPreliminaries}
\textbf{Notation:} For $n \in \mathbb{N}$, we use the symbol $[n]$ to denote the set $\{1, \ldots, n\}$. $A_{\cdot j}$ refers to the $j^{th}$ column of the matrix $A$. $\mathbb{C}_+ = \{u + \mathbbm{i}v: u \in \mathbb{R}, v > 0\}$ denotes the upper half of the complex plane excluding the real axis. $\|A\|$ denotes the operator norm of the matrix $A$. $\lambda_{\max}(A)$ ($\lambda_{\min}(A)$) denotes the largest (smallest) eigenvalue of the square matrix $A$. $\mathbbm{1}_n$ denotes the vector $(1, \ldots,1) \in \mathbb{C}^n$.

Let $K \in \mathbb{N}$ be fixed. For each $r \in [K]$, suppose $\{Z_{n}^{(r)}\}_{n=1}^{\infty}$ are sequences of random matrices, each having dimension $p \times n$ such that $p/n \rightarrow c \in (0, \infty)$. The entries have zero mean, unit variance and uniformly bounded moments of some order. Let $\{A_{n}^{(r)} \in \mathbb{C}^{p \times p}, B_n^{(r)} \in \mathbb{C}^{n \times n}; r \in [K]\}_{n=1}^\infty $ be sequences of random positive semi-definite matrices. We are interested in the limiting behavior (as $p, n \rightarrow \infty$) of the ESDs of matrices of the type 
\begin{align}\label{defining_Sn}
S_n &= \frac{1}{n}X_nX_n^*,  \text{ where } X_n = \sum_{r=1}^KX_{n}^{(r)},\\ \notag
\text{ and } 
X_{n}^{(r)} &= U_{n}^{(r)}Z_{n}^{(r)}V_{n}^{(r)}\,, U_n^{(r)} = (A_n^{(r)})^\frac{1}{2}\,, V_n^{(r)} = (B_n^{(r)})^\frac{1}{2}.
\end{align}
We define the following central objects of our analysis. Henceforth, the $Z_n^{(r)}$ matrices will be referred to as \textit{innovation} matrices and $A_n^{(r)},B_n^{(r)}$ as \textit{scaling} matrices.
\begin{definition}\label{defining_JESD}
Let $\{\lambda_{j}^{(1)}, \ldots, \lambda_{j}^{(K)}\}_{j=1}^p$ and $\{\theta_{j}^{(1)}, \ldots, \theta_{j}^{(K)}\}_{j=1}^p$ represent the eigenvalues of $\mathbb{A}_n := (A_{n}^{(1)}, \ldots, A_{n}^{(K)})$ and $\mathbb{B}_n := (B_{n}^{(1)}, \ldots, B_{n}^{(K)})$ respectively. Then we define the joint empirical spectral distribution (J.E.S.D.) of the $K$ matrices as follows.
\begin{align*}
    H_n(s_1,\ldots, s_K) &= F^{\mathbb{A}_n}(s_1,\ldots, s_K) := \frac{1}{p}\sum_{j=1}^p \prod_{r=1}^K\mathbbm{1}_{\{0 \leq \lambda_{j}^{(r)} \leq s_r\}},\\
    G_n(t_1,\ldots, t_K) &= F^{\mathbb{B}_n}(t_1,\ldots, t_K) := \frac{1}{n}\sum_{j=1}^n \prod_{r=1}^K\mathbbm{1}_{\{0 \leq \theta_{j}^{(r)} \leq t_r\}}.
\end{align*}
\end{definition}
\begin{definition}
The dual of the Sample Covariance Matrix is defined as
\begin{align}\label{defining_dual}
    \ubar{S}_n := \frac{1}{n}X_n^*X_n
\end{align}
For $j \in [n]$, the leave-one-out version of the sample covariance matrix is given by
\begin{align}\label{defining_Snj}
    S_{nj} : = \frac{1}{n}\sum_{i \neq j}X_{\cdot i}X_{\cdot i}^*, \quad j \in [n]
\end{align}
\end{definition}

\begin{definition}
For $z \in \mathbb{C}_+:= \{u+\mathbbm{i}v: v>0, u \in \mathbb{R}\}$, we define the following resolvent-type objects. The resolvent of $S_n$ is given by
\begin{align}\label{defining_Qz}
    Q(z) := (S_n - zI_p)^{-1},
\end{align}
that of $S_{nj}$ and $\ubar{S}_n$ are given by
\begin{align}\label{defining_Qzj}
    Q_{j}(z) := (S_{nj} - zI_p)^{-1}, \quad \ubar{Q}(z) := (\ubar{S}_n - zI_n)^{-1}.
\end{align}
\end{definition}
\begin{definition}
We now introduce the associated Stieltjes transforms for these resolvents. The Stieltjes transform for $F^{S_n}$ is
\begin{align}\label{defining_s_n}
s_n(z) := \frac{1}{p}\operatorname{tr}(Q(z)),
\end{align}
and that for $F^{\ubar{S}_n}$ is given by
\begin{align}\label{defining_s_n_bar}
    \ubar{s}_n(z) := \frac{1}{p}\operatorname{tr}(\ubar{Q}(z))
\end{align}
The following quantities will play a vital role in the matrix Dyson system.
\begin{align}\label{defining_hn}
\textbf{h}_{n}(z) := (h_{1n}(z), \ldots, h_{Kn}(z))^\top, \text{ where } h_{rn}(z) := \frac{1}{p}\operatorname{tr}\{A_{n}^{(r)}Q(z)\}, \quad r \in  [K].
\end{align}
\begin{align}\label{defining_gn}
\textbf{g}_{n}(z) := (g_{1n}(z), \ldots, g_{Kn}(z))^\top, \text{ where } g_{rn}(z) := \frac{1}{n}\operatorname{tr}\{B_{n}^{(r)}\ubar{Q}(z)\}, \quad r \in [K].
\end{align}
\end{definition}
\begin{remark}\label{h_rn_is_a_ST}
     Note that, $h_{rn}(\cdot)$ is the Stieltjes Transform of some measure. To see this, let $S_n = P_n\Lambda_nP_n^*$ be a spectral decomposition of $S_n$ and let $D_n^{(r)} = P_n^*A_n^{(r)}P_n$. Denote the diagonal elements of $D_n^{(r)}$ by $d_{jj}^{(r)}$ and denote the (real) eigenvalues of $S_n$ by $\lambda_j$ for $j \in [p]$. Then, observe that
     \begin{align}
         h_{rn}(z) =& \frac{1}{p}\operatorname{tr}(A_n^{(r)}Q(z)) \\ \notag
         =& \frac{1}{p}\operatorname{tr}(A_n^{(r)}P_n(\Lambda_n - zI_p)^{-1}P_n^*)
         = \frac{1}{p}\operatorname{tr}(D_n^{(r)}(\Lambda_n - zI_p)^{-1})
         = \frac{1}{p}\sum_{j=1}^p \frac{d_{jj}^{(r)}}{\lambda_j - z}.
     \end{align}
The RHS is the Stieltjes Transform of the measure that allocates a mass of $d_{jj}^{(r)}/p$ to the point $\lambda_j$. Similar results also hold for $g_{rn}(\cdot)$.
\end{remark}

\subsection{LSD under pairwise commutativity of scaling matrices}
We state the LSD result of (\ref{defining_Sn}) established in \cite{JHDP26}. For a probability measure $\mu$ on $\mathbb{R}_+^K$, $z \in \mathbb{C}_+$ and $\textbf{p} \in \mathbb{C}_+^K$, we define the following complex valued vector function:
\begin{align}\label{defining_Operator_O}
    \textbf{O}(z, \textbf{p}, \mu) := \int \frac{\boldsymbol{\lambda}\mu(d\boldsymbol{\lambda})}{-z(1 + \boldsymbol{\lambda}^\top \textbf{p})}.
\end{align}
Here, $\boldsymbol{\lambda}$ represents an arbitrary point in $\mathbb{R}_+^K$. $O_r(z, \textbf{p}, \mu)$ is the $r^{th}$ coordinate of $\textbf{O}(z, \textbf{p}, \mu)$.
\begin{theorem}\label{mainTheorem_LSD}
    Suppose the following conditions hold.
    \begin{enumerate}
        \item[\textbf{T1}] $c_n := p/n \rightarrow c \in (0, \infty)$
        \item[\textbf{T2}] Entries of $Z_n^{(1)},\ldots, Z_n^{(K)}$ are independent and have zero mean, unit variance and have a uniform bound on moments of order $2+\eta_0$ for some $\eta_0 > 0$.
         \item[\textbf{T3}] (Commutativity) $A_{n}^{(r)}A_{n}^{(s)} = A_{n}^{(s)}A_{n}^{(r)}$ and $B_{n}^{(r)}B_{n}^{(s)} = B_{n}^{(s)}B_{n}^{(r)}$ for $r,s \in [K]$
        \item[\textbf{T4}] $\{A_{n}^{(r)}, B_{n}^{(r)}: r \in [K]\}_{n=1}^\infty$ are sequences of p.s.d. matrices such that $F^{\textbf{A}_{n}} \xrightarrow{d} H$ a.s. and $F^{\textbf{B}_{n}} \xrightarrow{d} G$ a.s.  where $H,G$ are non-random probability distributions on $\mathbb{R}_+^K$ such that $H \neq \delta_\textbf{0}$, $G\neq \delta_\textbf{0}$ and they are not supported entirely on any of the $K$ axes. 
        \item[\textbf{T5}]         
        There exists a constant $C_0>0$ such that $$\underset{n \rightarrow \infty}{\limsup} ~\underset{r \in [K]}{\max} \max\bigg\{\frac{1}{p}\operatorname{tr}(A_{n}^{(r)}), \frac{1}{n}\operatorname{tr}(B_{n}^{(r)})\bigg\} < C_0.$$
    \end{enumerate}
Then, $F^{S_n} \overset{d}{\longrightarrow} F^{c,H,G}$ a.s. where $s_F$, the Stieltjes transform of $F$ at $z \in \mathbb{C}_+$ satisfies the following:
    \begin{align}\label{s_main_eqn}
        s_F(z) = \displaystyle \int_{\mathbb{R}_+^K} \frac{dH(\boldsymbol{\lambda})}{-z(1 + \boldsymbol{\lambda}^\top \textbf{O}(z, c\textbf{h}(z), G))}.
    \end{align}
    where, $\textbf{h}(z) = (h_1(z),\ldots, h_K(z))^\top  \in \mathbb{C}_+^K$ is the unique solution within the class of Stieltjes transforms of measures and satisfies
    \begin{align}\label{h_main_equation}
         \textbf{h}(z) &= \displaystyle \int_{\mathbb{R}_+^K} \frac{\boldsymbol{\lambda} dH(\boldsymbol{\lambda})}{-z(1 + \boldsymbol{\lambda}^\top \textbf{O}(z, c\textbf{h}(z), G))}
    \end{align}    
\end{theorem}

\begin{corollary}\label{corollary_dual_results}
    Under \textbf{T1-T5} of Theorem \ref{mainTheorem_LSD}, we have an alternate characterization for the Stieltjes Transform of the LSD $F$ of $S_n$.
\begin{align}\label{s_alt_char_2}
    s_F(z) &= \frac{1}{c}\int \frac{dG(\boldsymbol{\theta})}{-z(1 +c\boldsymbol{\theta}^\top \textbf{O}(z, \textbf{g}(z), H))} + \bigg(\frac{1}{c}-1\bigg)\frac{1}{z}
\end{align}
where $\textbf{g}(z) = (g_1(z), \ldots, g_K(z))^\top  \in \mathbb{C}_+^K$ is the unique solution within the class of Stieltjes transforms of measures and satisfies
\begin{align}\label{g_main_equation}
    \textbf{g}(z) &= \int \frac{\boldsymbol{\theta} dG(\boldsymbol{\theta})}{-z(1 + c\boldsymbol{\theta}^\top \textbf{O}(z, \textbf{g}(z), H))}
\end{align}
Moreover \textbf{h} (from \ref{h_main_equation}) and \textbf{g} satisfy the below.
\begin{align}\label{hg_recursive}
    \textbf{g}(z) = \textbf{O}(z, c\textbf{h}(z), G); \quad \textbf{h}(z) = \textbf{O}(z, \textbf{g}(z), H)
\end{align}
Yet another characterization of the Stieltjes Transform of the L.S.D. is given by 
\begin{align}\label{s_alt_char_1}
     s(z) =\frac{1}{z}\bigg(\frac{1}{c}-1\bigg) - \frac{1}{cz}\int \frac{dG(\boldsymbol{\theta})}{1 + c\boldsymbol{\theta}^\top \textbf{h}(z)} = -z^{-1} - \textbf{h}^\top (z)\textbf{g}(z)
\end{align}
\end{corollary}

\subsection{Properties of Finite Sample Solutions}
Let $(s_n^0(z), \textbf{h}_n^0(z), \textbf{g}_n^0(z))$ be the finite sample (unique) solutions of the Matrix Dyson system corresponding to parameters $c_n, H_n$ and  $G_n$. From \cite{JHDP26}, it is known that for $r \in [K]$, they satisfy the following:
\begin{align}\label{finite_sample_properties}
    s_n^0(z) &= \frac{1}{p}\sum_{j=1}^p \frac{1}{1 + \boldsymbol{\lambda}_j^\top \textbf{g}_n^0(z)}\\ \notag
    -zh_{rn}^0(z) &= \frac{1}{p}\sum_{j=1}^p \frac{\lambda_{rj}}{1 + \boldsymbol{\lambda}_j^\top \textbf{g}_n^0(z)}\\ \notag    
    -zg_{rn}^0(z) &= \frac{1}{n}\sum_{j=1}^n \frac{\theta_{rj}}{1 + c_n\boldsymbol{\theta}_j^\top \textbf{h}_n^0(z)}.
\end{align}

Let $\ubar{s}_n(z) = \frac{1}{n}\operatorname{tr}\ubar{Q}(z)$ represent the Stieltjes transform of the dual of the sample covariance. Then,
\begin{align}\label{primal_dual_relationship}
    \operatorname{tr}Q(z) = \operatorname{tr}\ubar{Q}(z) + (p-n)\frac{1}{-z}
    \iff  s_n(z) = \frac{1}{c_n}\ubar{s}_n(z) + \bigg(1 - \frac{1}{c_n}\bigg)\frac{1}{-z}.
\end{align}

\color{black}
\subsection{Main Result}\label{sec:MainResult}
Let $\lambda_j^{S_n};j\in [p]$ denote the eigenvalues of $S_n$. For bounded continuous real-valued functions $f$ on $\mathbb{R}$, the linear spectral statistics of $S_n$ are of the following form:
\begin{align}\label{defining_LSS}
    L_n(f) = \frac{1}{p}\sum_{j=1}^pf(\lambda_j^{S_n}).
\end{align}
We will show that under the assumption that $\underset{i,j,r,n}{\sup}\mathbb{E}|z_{ij}^{(r)}|^4 < \infty$ and analyticity of $f$, the quantity 
$$\absmod{\int f(x)dF^{S_n}(x)dx - \int f(x)dF^{c_n,H_n,G_n}(x)dx},$$
converges to $0$ at a rate of $1/n$. Define the \textit{centered} and scaled empirical process:
\begin{align}\label{defining_Enx}
    \xi_n(x) := p(F^{S_n}(x) - F^{c_n, H_n, G_n}(x)).
\end{align}
Our main result involves functions which are analytic on an open set containing a compact interval associated with the support of the LSD. We describe this interval below. For $r \in [K]$, let $\theta_{rn} \leq \ldots \leq \theta_{r1}$ bet the eigenvalues of $B_r$. Then, define
\begin{align}\label{defining_interval}
    e_l &:= K^2 \min_{1\leq r \leq K}\liminf_{n \rightarrow \infty}\theta_{rn}\lambda_{\min}(A_r)I_{(0,1)}(c)(1-\sqrt{c})^2,\\ \notag
    e_r &:= K \sum_{r=1}^K\limsup_{n \rightarrow \infty}{\theta_{r1}\lambda_{\max}(A_r)(1+\sqrt{c})^2}.
\end{align}
These are essentially a.s. bounds for the eigenvalues of $S_n$ and they follow from Weyl's and Cauchy-Schwarz inequalities. Clearly, if the spectral norms of $A_r,B_r$ are uniformly bounded, $[e_l, e_r]$ is a compact interval of $\mathbb{R}_+$. The statement of our main result is extremely notation-heavy and we dedicate Subsection \ref{sec:Glossary} entirely to introduce these sequentially.
\begin{theorem}\label{mainTheorem_CLT}
Suppose the following conditions hold:
\begin{enumerate}
    \item Let $z_{ij}^{(r)} = [Z_n^{(r)}]_{ij}$, for $i \in [p], j \in [n]$ represent the innovation entries that satisfy the following Lindeberg-type condition. For any $\delta>0$, as $n \rightarrow \infty$, 
\begin{align}\label{L4}
    \frac{1}{Knp\delta^4}\sum_{i,j,r}\mathbb{E}|z_{ij}^{(r)}|^{4}\mathbbm{1}_{\{|z_{ij}^{(r)}|>\delta\sqrt{n}\}} \longrightarrow 0.
\end{align}
    \item For some $N \in \mathbb{N}$, let $(f_1, \ldots, f_N)$ be analytic functions on $\mathbb{R}$ analytic on an open interval containing $[e_l, e_r]$. Let $$I_j := \int f_j(x)d\xi_n(x); j \in [N].$$
    \item $\{A_r\}_{r=1}^K$ and $\{B_r\}_{r=1}^K$ are non-random matrices and for each $n \in \mathbb{N}$ their spectral norms satisfy
\begin{align}\label{bounded_spectral_norm}
    \sup_{n\in \mathbb{N}}\, \max_{r \in [K]} \, \{\|A_r\|, \|B\|_r\} \leq \tau < \infty.
\end{align}
\end{enumerate}
Then, the following holds.
\begin{enumerate}
    \item[1] \textbf{Real Case:} If $Z_r, A_r, B_r$ are real and $\mathbb{E}|z_{ij}^{(r)}|^4 = 3$, 
then,
$(I_1, \ldots, I_N)$ converges weakly to a Gaussian vector ($\mathcal{G}_1, \ldots, \mathcal{G}_N$) with mean given by
\begin{align}\label{mean_Xf}
    \mathbb{E}\mathcal{G}_i = -\frac{1}{2\pi\mathbbm{i}}\oint_{\mathcal{C}} f_i(z) m_2(z)dz,
\end{align}
    and covariance given by
\begin{align}\label{cov_Xf_Xg}
    \operatorname{Cov}(\mathcal{G}_i, \mathcal{G}_j) = -\frac{1}{2\pi^2}\oint_{\mathcal{C}_1}\oint_{\mathcal{C}_2} f_i(z_1)f_j(z_2)\operatorname{ker}(z_1, z_2) dz_1dz_2,
\end{align}
where $m_2(z), \operatorname{ker}(z_1, z_2)$ are defined in  (\ref{defining_population_mean_function}) and (\ref{defining_population_covariance_function}). The contours in (\ref{mean_Xf}) and (\ref{cov_Xf_Xg}) are assumed to be closed, non-overlapping and positively oriented such that they enclose the support of $F^{c, H, G}$.
\item[2] \textbf{Complex Case:} If $\{Z_r = W_r + \mathbbm{i}Y_r\}_{r=1}^K$ are complex matrices with $\mathbb{E}w_{ij}^{(r)} = 0 = \mathbbm{E}y_{ij}^{(r)}$, $\mathbb{E}(w_{ij}^{(r)})^2 = 1/2 = \mathbb{E}(y_{ij}^{(r)})^2$, and $\mathbb{E}(w_{ij}^{(r)})^4 = 3/4 = \mathbb{E}(y_{ij}^{(r)})^4$, then the result holds except that the mean function in (\ref{mean_Xf}) is zero and the covariance function is half of that given in (\ref{cov_Xf_Xg}).
\end{enumerate}
\end{theorem}

\begin{remark}
    The expressions (\ref{mean_Xf}) and (\ref{cov_Xf_Xg}) reduce to the ones of \cite{BaiLiPan2019} when $K = 1$.
\end{remark}

The proof relies on establishing the LSS as the smoothed average of a complex random process (that converges weakly to a Gaussian Process) on a contour surrounding the support of the LSD. By \textit{Cauchy’s formula}, we have 
\begin{align}\label{cauchys_formula}
    \int f(x)dG(x) = -\frac{1}{2\pi i}\int_\gamma f(z)m_G(z)dz,
\end{align}
where $G$ is a cumulative distribution function, $m_G(\cdot)$ is the Stieltjes transform of $G$ and $f$ is analytic on an open set containing the support of $G$. The complex integral on the RHS is over any positively oriented
contour $\gamma$ enclosing the support of $G$ and on which $f$ is analytic. As a result, we express the central object of Theorem \ref{mainTheorem_CLT} as follows:
\begin{align*}
    \int_\mathcal{C} f(\lambda)d\xi_n(\lambda) &= p\bigg(\int_\mathcal{C} f(\lambda)dF^{S_n}(\lambda) - \int_\mathcal{C} f(\lambda)dF^{c_n, H_n, G_n}(\lambda)\bigg)\\
    &= p\int_\mathcal{C} f(z)(s_n(z) - s_n^0(z))dz.
\end{align*}
Here $s_n^0(\cdot)$ represents the finite sample solution of (\ref{s_main_eqn}) corresponding to $F^{c_n,H_n,G_n}$. Hence, the proof of Theorem \ref{mainTheorem_CLT} relies on establishing limiting results for
\begin{align}\label{decomposing_Mn}
    M_n(z) &= p(s_n(z) - s_n^0(z))\\ \notag
           &= p(s_n(z) - \mathbb{E}s_n(z)) + p(\mathbb{E}s_n(z) - s_n^0(z))\\ \notag
           &:= M_{n1}(z) + M_{n2}(z).
\end{align}
So we have expressed an integral over the (positive) real line to a complex contour integral involving the familiar Stieltjes transform. $M_{n1}$ is the random component which will be shown to converge to a Gaussian process. The bias term that remains is the limit of the deterministic part $M_{n2}$. We will show that it suffices to consider a truncated contour to avoid points near the real axis where stability issues creep in.

\subsection{Construction of the Contour}
By the assumption of Theorem \ref{mainTheorem_CLT}, we have $$\sup_n\underset{r \in [K]}\max\{\|A_r\|,\|B_r\|\} \leq \tau.$$ 
We will now construct a rectangular contour using three parameters $x_l$, $x_r$ and $v_0>0$. Let $\mathcal{C}_u = \{x + \mathbbm{i}v_0: x \in [x_l,x_r]\}$. Then, we define the following contour: 
\begin{align}\label{defining_the_contour}
  \mathcal{C} = \{x_l + \mathbbm{i}v: v \in [0, v_0]\} \cup \mathcal{C}_u \cup \{x_r + \mathbbm{i}v: v \in [0, v_0]\}.  
\end{align}
Here, $x_l$ and $x_r$ denote the left and right ends respectively, on the real axis. Similarly, $\pm v_0$ denote the top and bottom points respectively, on the imaginary axis. To ensure that $\mathcal{C}$ encloses the support of $F^{c,H,G}$ completely, we impose some additional constraints on $v_0, x_l, x_r$.
\begin{enumerate}
    \item If the right end point of the interval (\ref{defining_interval}), i.e. $e_r$ is zero, we choose $x_r$ to be any positive number. Otherwise choose
        $$x_r \in (e_r,\infty).$$
    \item If the left end point of (\ref{defining_interval}), i.e. $e_l$ is zero, we choose $x_l$ to be any negative number. Otherwise choose
    \begin{align}
        x_l \in (0, e_l).
    \end{align}
\end{enumerate}
To avoid dealing with small values of $\Im(z)$, we truncate $M_n(z)$ on a contour $\mathcal{C}$ of the complex plane. We define the subsets $\mathcal{C}_n$ of $\mathcal{C}$ on which $\widehat{M}_n(\cdot)$. Choose sequence $\{\epsilon_n\}$ decreasing to $0$ satisfying for some $\alpha \in (0,1)$
\begin{align}\label{choice_of_epsilon_n}
  \epsilon_n \geq n^{-\alpha}  
\end{align}
Let $$\mathcal{C}_l =\{x_l+\mathbbm{i}v: v \in [\epsilon_n/n, v_0]\}, \quad \mathcal{C}_r =\{x_r+\mathbbm{i}v: v \in [\epsilon_n/n, v_0]\}.$$
Then, the truncated contour is 
\begin{align}\label{script_Cn}
 \mathcal{C}_n = \mathcal{C}_l \cup \mathcal{C}_u \cup \mathcal{C}_r.   
\end{align}
For $z = x  + \mathbbm{i}v$, the process $\widehat{M}_n(\cdot)$ is now defined as
\begin{align}\label{defining_truncated_contour}
    \widehat{M}_n(\cdot) = \begin{cases}
        M_n(z), \hspace*{\fill} &\text{ for } z \in \mathcal{C}_n,\\
        M_n(x_l + \mathbbm{i}n^{-1}\epsilon_n), &\text{ for } x = x_l, v \in [0, n^{-1}\epsilon_n],\\
        M_n(x_r + \mathbbm{i}n^{-1}\epsilon_n), &\text{ for } x = x_r, v \in [0, n^{-1}\epsilon_n].
    \end{cases}
\end{align}
We now state a wrapper lemma that will serve as a stepping stone to prove Theorem \ref{mainTheorem_CLT}.
\begin{lemma}\label{lemma:GaussianProcess}
    Recall the definitions of $U(\cdot, \cdot)$ and $W(\cdot, \cdot)$ from (\ref{defining_population_U_matrices}) and (\ref{defining_population_W_matrices}). Under the conditions of Theorem \ref{mainTheorem_CLT}, $\widehat{M}_n(z)$ converges weakly to a two-dimensional Gaussian process $M(\cdot)$ for $z \in \mathcal{C}$ under the \textbf{real case} with
    \begin{align}
        \mathbb{E}M(z) = m_2(z),
    \end{align}
    and for $z_1,z_2 \in \mathcal{C} \cup \bar{\mathcal{C}}$ with $\bar{\mathcal{C}} = \{\bar{z}: z \in \mathcal{C}\}$,
    \begin{align}\label{covariance_function_real}
        \operatorname{Cov}(M(z_1), M(z_2))
        &= 2\operatorname{ker}(z_1, z_2).
    \end{align}
    Under the \textbf{complex case}, $\mathbb{E}M(z) = 0$ and the covariance function is half of (\ref{covariance_function_real}). $m_2(\cdot)$ and $\operatorname{ker}(\cdot, \cdot)$ are defined in (\ref{defining_population_mean_function}) and (\ref{defining_population_covariance_function}).
\end{lemma}

\textbf{Choice of} $[\eta_l, \eta_r]$:

From \cite{BaiSilv04} and \cite{BaiYinKri}, for any integer $l$, $\eta_1 > (1+\sqrt{c})^2, \eta_2 < (1-\sqrt{c})^2$, we have
\begin{align*}
    &\mathbb{P}\bigg(\lambda_{\max}\bigg(\frac{1}{n}ZZ^*\bigg)\geq \eta_1\bigg) = o(n^{-l}),\quad
\mathbb{P}\bigg(\lambda_{\min}\bigg(\frac{1}{n}ZZ^*\bigg)\leq \eta_2\bigg) = o(n^{-l}).
\end{align*}
Using this, we will construct $[\eta_l, \eta_r]$ which will contain the eigenvalues of $S_n$ with very high probability. Let 
\begin{align}\label{eta_l_r}
    \eta_l \in (x_l, e_lI_{(0,1)}(c)), \quad \eta_r \in (e_r, x_r)
\end{align}
In particular, we can choose $x_l,x_r,\eta_l, \eta_r$ around the points $e_l, e_r$ such that 
\begin{align}\label{contour_gap_technical_condition}
    \min\{|x_r - \eta_r|,  |\eta_l - x_l| \} > \max\{4K\tau^2, 4Kc\tau^2\}.
\end{align}
We therefore have
\begin{align}\label{extreme_eigenvalues}
    \mathbb{P}(\lambda_{\max}(S_n) \geq \eta_r) = o(n^{-l}), \quad \mathbb{P}(\lambda_{\min}(S_n) \leq \eta_l) = o(n^{-l}).
\end{align}

\subsection{Proof of Theorem \ref{mainTheorem_CLT}}
Assuming Lemma \ref{lemma:GaussianProcess}, which will be proved sequentially throughout Sections \ref{sec:GaussianCase} and \ref{sec:NonGaussian}, we are in a position to prove Theorem \ref{mainTheorem_CLT}. We establish the following extreme eigenvalue bounds using \cite{BaiYin} and \cite{BaiYinKri}. For any $r \in [K]$, the following holds:
$$\lambda_{\max}\bigg(\frac{1}{n}Z_rZ_r^*\bigg) \rightarrow (1 + \sqrt{c})^2 \quad \text{a.s. and}  \quad \lambda_{\min}\bigg(\frac{1}{n}Z_rZ_r^*\bigg) \rightarrow (1 - \sqrt{c})^2 \quad \text{a.s.}$$
We have,
\begin{align*}
    \lambda_{\max}(S_n) &=\frac{1}{n}\|X\|^2= \bigg\|\sum_{r=1}^K A_r^\frac{1}{2} \bigg(\frac{1}{\sqrt{n}}Z_r\bigg) B_r^\frac{1}{2}\bigg\|^2\\
        &\leq K \sum_{r=1}^K \|A_r\| \times \|B_r\| \times \|\frac{1}{n}Z_rZ_r^*\|,
\end{align*}
so that, $\limsup \lambda_{\max}(S_n) \leq K^2 \tau^2(1+\sqrt{c})^2$ a.s.
Similarly, we also have
\begin{align*}
   \liminf \lambda_{\min}(S_n) \geq K^2 \min_{1\leq r \leq K}(\liminf \lambda_{p}(A_r)\theta_{rn})I_{(0,1)}(c)(1-\sqrt{c})^2, \quad a.s.
    \end{align*}

Combining the definitions of $x_l,x_r,$ we find with probability $1$, the following holds:
$$ \underset{n \rightarrow \infty}{\liminf} \, \, \min \,\, \{x_r - \lambda_{\max}(S_n), \lambda_{\min}(S_n) - x_l\} > 0.$$
Define the \textit{eigen-gap} as follows:
$$e_{gap} := \max\bigl\{\max(\lambda_{\max}(S_n), e_r) - x_r, \min(\lambda_{\min}(S_n), e_l) - x_l\bigr\},$$
where $e_l$ and $e_r$ are defined in (\ref{defining_interval}). Since $F^{S_n} \rightarrow F^{c,H,G}$ with probability $1$, the support of $F^{c_n,H_n,G_n}$ is contained in the interval (\ref{defining_interval}) with probability $1$. Thus, by (\ref{cauchys_formula}), for $f \in \{f_1, \ldots, f_N\}$ and large $n$, with probability $1$,
    $$\oint f(x)d\xi_n(x) = -\frac{1}{2\pi\mathbbm{i}}\oint f(z) M_n(z)dz$$
where the contour integral is over $\mathcal{C} \cup \bar{\mathcal{C}}$. Then for $v \in [0, \epsilon_n/n]$, we have
$$\|Q(x_r + \mathbbm{i}v)\| \leq \frac{1}{e_{gap}}, \quad a.s.$$
As a result, for some absolute constant $C>0$, we have
\begin{align*}
    |M_n(x_r + \mathbbm{i}v) - M_n(x_r + \mathbbm{i}\epsilon_n/n)| 
    &\leq \frac{Cn}{e_{gap}},
\end{align*}
and similarly,
\begin{align*}
    |M_n(x_l + \mathbbm{i}v) - M_n(x_l + \mathbbm{i}\epsilon_n/n)| \leq \frac{Cn}{e_{gap}}.
\end{align*}
Since $\sup_{z \in \mathcal{C}}|f(z)| < \infty$, the following holds with probability $1$:
\begin{align}\label{contour_truncation_result}
    \absmod{\oint f(z)(M_n(z) - \widehat{M}_n(z))dz}\leq \frac{4C\epsilon_n}{e_{gap}} \longrightarrow 0.
\end{align}
Since the mapping from $C(\mathcal{C}, \mathbb{R}^2)$ to $\mathbb{C}^N$ given by
$$\mathbbm{i}\widehat{M}_n(\cdot) \longmapsto -\frac{1}{2\pi}\bigg( \oint f_1(z)\widehat{M}_n(z)dz, \ldots, \oint f_N(z)\widehat{M}_n(z)dz\bigg)$$
is continuous, the proof of Theorem \ref{mainTheorem_CLT} is complete using Lemma \ref{lemma:GaussianProcess}. The rest of the paper is devoted to proving the wrapper lemma.

\subsection{Applications}
We propose an application that will test (assuming the true model has separable structure) the hypothesis that the row and column variances are spherical in nature. Note that under the null hypothesis, we can assume $K=1$.

Consider the ratio $\mu_{2n}/\mu_{1n}^2$, from the observed data, where $\mu_{rn} := \mathbb{E}\int \lambda^rdF^{S_n}(\lambda)$. We then have an approximation to the ratio in terms of analogous ratios of spectral moments of the row and column covariance matrices as expressed below.
\begin{align}\label{eq:ratio_moment}
    \frac{\mu_{2n}}{\mu_{1n}^2}
    = \frac{\alpha_2}{\alpha_1^2} + c_n \frac{\beta_2}{\beta_1^2} + \frac{1}{p}\frac{c_n\alpha_2\beta_2}{\alpha_1^2\beta_1^2} \approx \frac{\alpha_2}{\alpha_1^2} + c \frac{\beta_2}{\beta_1^2},
\end{align}
where $\alpha_j = \int \lambda^j dF^{A}(\lambda)$ and $\beta_j = \int \theta^j dF^{B}(\theta)$, with $F^{A}$ and $F^{B}$ being the limiting spectral distributions of $A$ and $B$, respectively. A proof of (\ref{eq:ratio_moment}) is provided in (\ref{proof:eq:ratio_moment}). 

Under the null hypothesis of sphericity,  
$$H_0: A = a_1I_p, B = a_2I_n, \quad \text{for unknown constants } a_1, a_2 >0,$$
we expect the estimated value to be close to $1 + c$, which is also the lower bound on the right hand side of (\ref{eq:ratio_moment}). Therefore, large exceedance of this value should indicate departure from sphericity. The CLT for linear spectral statistics derived here can then be used to study the asymptotic  behavior of the corresponding test statistic under certain classes of departure from the null, namely, when the data matrices are sums of independent separable components, with simultaneously diagonalizable scaling matrices.

As another potential application, notice that the approximation (\ref{eq:ratio_moment}) also serves as an estimating equation for the quantity $(\alpha_2/\alpha_1^2) + c (\beta_2/\beta_1^2)$. The CLT derived here can be used to obtain the asymptotic distribution of its method-of-moments estimator under the assumption of separability.

\section{Proof Sketch}\label{sec:proofSketch}
The following serves as an outline of the proof of Lemma \ref{lemma:GaussianProcess}. Items (1), (2) and (4) follow standard techniques and hence have been moved to the Appendix.
\begin{enumerate}
    \item In Section \ref{sec:Truncation_Centralization_rescaling}, we show that under Condition (1) of Theorem \ref{mainTheorem_CLT}, we can assume that the innovation entries $\{z_{ij}^{(r)}\}_{i,j,r}$ are truncated at $\delta_n\sqrt{n}$, centered and scaled to unit variance. Here, $\{\delta_n\}$ is a sequence that goes to 0 such that $n\delta_n^4\to\infty$. Let $LSS(Z), LSS(\widehat{Z})$ represent the LSS's arising from the original innovations as in Theorem \ref{mainTheorem_CLT} and their truncated versions respectively. Using similar techniques as \cite{BaiSilv04}, we show that
    $$|LSS(Z) -  LSS(\widehat{Z})| \overset{p}{\longrightarrow} 0.$$
    \item In Section \ref{sec:NonGaussian}, we show that 
    $$d_L(LSS(W), LSS(\widehat{Z})) \longrightarrow 0,$$
    where LSS(W) represents the LSS arising from assuming all the innovations to be i.i.d. Gaussian and $d_L(\cdot, \cdot)$ is the Levy metric between distribution functions. The key is to show that the characteristic function of the LSS arising out of the truncated innovations ($\widehat{Z}$) and the one arising from the Gaussian innovations ($\widehat{Z}$) are arbitrarily close to each other for large $n$. This is achieved by constructing entry wise linear combinations of $\widehat{Z}$ and $W$ using variance preserving trigonometric functions as the weights.   

    Therefore, without loss of generality, we can assume the entries to be i.i.d. Gaussian. In particular, this implies that we can assume $A_r,B_r$ to be diagonal matrices. Letting $LSS(\widehat{W})$ represent the LSS arising by truncating the Gaussian innovations at $\delta_n\sqrt{n}$, by Step 1, we have 
    $$|LSS(W) - LSS(\widehat{W})| \overset{p}{\longrightarrow} 0.$$
    As a result, we assume the innovations to be i.i.d. truncated Gaussian in nature.
    
    \item Section \ref{sec:essential_quantities_bounds} introduces essential quantities such as leave-one-out (LOO) deflators and some uniform bounds and concentration results. As a foundational result, we establish that the norm of resolvent-type objects have uniformly bounded moments even as we approach the real line along the vertical arms of the contour.
    This itself leads to a string of powerful results stated as various lemmas throughout the section.

    \item The tightness of the process $M_{n1}(z)$ over $z \in \mathcal{C}_n$ is established in Section \ref{sec:Tightness} using Lemma \ref{lemma:tightness_condition}.
    \item We show that the finite dimensional distributions of $M_{n1}(z)$ have Gaussian limits. We start by expressing $M_{n1}(z)$ as a telescopic sum of martingale differences which are functions of LOO deflators and centered quadratic forms and this satisfies (ii) of Theorem \ref{theorem:martingale_CLT}.
    \item After replacing the random LOO deflators with their in-probability deterministic limits, we are left with an expression of the form
    $$\mathbb{E}_j\epsilon_j(z_1)\mathbb{E}_j\epsilon_j(z_2).$$
    On using Lemma \ref{lemma:covariance_quadratic_forms}, we recover a scaled tracial term composed of elements of the form $$\operatorname{tr}A_r\mathbb{E}_jQ_j(z_1)A_s\mathbb{E}_jQ_j(z_2), \quad r,s, \in [K].$$
    \item The rest of Section \ref{sec:GaussianCase} proceeds to sequentially replace the LOO resolvents/ deflators with their LTO counterparts to derive a system of self-consistent equation involving this random tracial expression and some deterministic expressions and some random error of small magnitude. This allows us to \textit{solve} for this quantity of interest and apply limiting results. The solution of this system of equation involves the invertibility of an operator which we explicitly show in Section \ref{sec:invertibility_statbility_operator} 
    \item The limit of $M_{n2}(z)$ is recovered in Section \ref{sec:BiasTerm} by performing decompositions of $\mathbb{E}Q(z) - R(z)$ and $R(z) - R^0(z)$. The basic idea is to express the quantity of interest in terms of a set of tracial expressions that solve a self-consistent system of equations upto some small order terms and follow steps similar to those that led to the expression for the covariance function. The bias term vanishes in the limit when the innovations are complex. However, when the innovations are real, we have a residual limiting bias term. This is completely an artefact of the covariance of quadratic forms taking different forms under the two scenarios (Lemma \ref{lemma:covariance_quadratic_forms}).
\end{enumerate}

\section{Essential Quantities and Results}\label{sec:essential_quantities_bounds}
By Items (1) and (2) of the proof sketch, it suffices to prove Lemma \ref{lemma:GaussianProcess} when the innovations are i.i.d. standard Gaussian, i.e. $\{z_{ij}^{(r)}\}, i \in [p], j \in [n], r \in [K]$ are i.i.d. standard real or complex Gaussian random variables. Since $\{B_r\}_{r=1}^K$ are pairwise commuting  p.s.d. matrices, there exist unitary matrices $U_n$ such that
$$B_r = U_n\operatorname{diag}(\theta_{r1},\ldots, \theta_{rn})U_n^* = U_n\Theta_rU_n^*.$$
Since $Z_rU_n \overset{d}{=}Z_r$ we note that
\begin{align*}
    S_n &= \frac{1}{n}\sum_{r,s=1}^K A_r^\frac{1}{2}Z_rB_r^\frac{1}{2}B_s^\frac{1}{2}Z_s^*A_s^\frac{1}{2} = \frac{1}{n}\sum_{r,s=1}^n A_r^\frac{1}{2}Z_rQ_n\Theta_r^\frac{1}{2}\Theta_s^\frac{1}{2}Q_n^*Z_s^*A_s^\frac{1}{2}\\
    &\overset{d}{=} \frac{1}{n}\sum_{r,s=1}^n A_r^\frac{1}{2}Z_r\Theta_r^\frac{1}{2}\Theta_s^\frac{1}{2}Z_s^*A_s^\frac{1}{2}.
\end{align*}
As a consequence, we will henceforth assume the $B_r$ matrices to be diagonal. Again using item (1) of the proof sketch, we will further assume that the innovations $z_{ij}^{(r)}$ are truncated at $\delta_n\sqrt{n}$, centralized and renormalized. Moreover, under the complex case, we assu me that
$\mathbb{E}|z_{ij}^{(r)}|^2 = o(n^{-1}), \mathbb{E}|z_{ij}^{(r)}|^2 = 2 + o(1),$ and $\mathbb{E}|z_{ij}^{(r)}|^4 = 3 + o(1)$ under the real case.

For $j \in [n]$, define $\mathcal{F}_j = \sigma(\{Z_{\cdot i}^{(r)}; i \in [j], r \in [K])$ and let $\mathbb{E}_j$ denote expectation operator with respect to the sigma field $\mathcal{F}_j$. Also by $\mathbb{E}_{-j}$, we will denote the conditional expectation with respect to the sigma field generated by all the columns \textit{except} the $j-$th. We further define the difference operator
$$\Delta_j:= \mathbb{E}_j - \mathbb{E}_{j-1}.$$

We will now sequentially define a list of essential quantities that we shall use in our analysis. Let $B_r = \operatorname{diag}(\{\theta_{rj}\}_{j=1}^n)$ for $r \in [K]$ and $\boldsymbol{\theta}_j = (\theta_{1j}, \ldots, \theta_{Kj})^\top$ for $j \in [n]$. Let 
$$x_{r,j} := n^{-\frac{1}{2}}A_r^\frac{1}{2}Z_{\cdot j}^{(r)}, x_j := \sum_{r=1}^K \sqrt{\theta_{rj}}x_{r,j}, \quad r \in [K], j \in [n].$$ Then, $\sqrt{n}x_j$ denotes the $j^{th}$ column of $X$ and we can express $S_n$ as
$$S_n = \sum_{j=1}^n x_jx_j^*.$$
The \textit{weighted} average of the $A_r$ matrices is defined as follows:
\begin{align}\label{weighted_average_A}
    \textbf{A}_j = \sum_{r=1}^K \theta_{rj}A_r, \quad j \in [n].
\end{align}
With this notation, $x_j$ is a random vector with zero mean and variance matrix $n^{-1}\textbf{A}_j$.

We introduce the basic resolvent and the leave-one-out (LOO) resolvent objects which are the critical elements of our work.
 \begin{align}\label{defining_Qz_Qjz}
     Q(z) = (S_n - zI_p)^{-1}, \quad S_{nj} = \sum_{i\neq j}x_ix_i^*, \quad Q_j(z) = (S_{nj} - zI_p)^{-1}.
 \end{align}
The following is a basic property of resolvents and will be used in many places.
\begin{align}\label{Leave_one_out_property}
    Q(z) - Q_j(z) = -\beta_j(z)Q_j(z)x_jx_j^*Q_j(z).
\end{align}
Observe that $\mathbb{E} x_j^*Q_jx_j  = n^{-1}\sum_{j=1}^n \theta_{rj}\operatorname{tr}(A_rQ_j).$ This tracial quantity will play a key role in our analysis and we denote the matrix under the trace operator as follows:
\begin{align}\label{defining_script_Aj}
    \mathcal{A}_j(z) = \frac{1}{n}\textbf{A}_jQ_j(z) = \frac{1}{n}\sum_{r=1}^K \theta_{rj}A_rQ_j(z), \quad j \in [n].
\end{align}
Analogous to $\mathcal{A}_j$ which depends on the LOO-resolvent, we also define below its counterpart which depends on the \textit{full}-resolvent
\begin{align}\label{script_Bj}
    \mathcal{B}_j(z) = \frac{1}{n}\textbf{A}_jQ(z)=\frac{1}{n}\sum_{r=1}^K \theta_{rj}A_rQ(z), \quad j \in [n].
\end{align}

Finally, we introduce the essential centered quadratic forms and the various versions of the LOO-deflators as follows:
\begin{align}\label{essential_quantities}
    &\epsilon_j(z) := x_j^*Q_jx_j - \operatorname{tr}(\mathcal{A}_j), \,\,\,\gamma_{j}(z) := x_j^*Q_j^2(z)x_j - \alpha_j(z), \, \alpha_{j}(z) := n^{-1}\sum_{r=1}^K\theta_{rj}\operatorname{tr}(A_rQ_j^2)\\ \notag
    &\beta_j(z) := (1 + x_j^*Q_jx_j)^{-1}, \quad \widetilde{\beta}_j(z) := (1 + \operatorname{tr}(\mathcal{A}_j))^{-1}, \quad b_j(z) := (1 + \mathbb{E}\operatorname{tr}(\mathcal{A}_j))^{-1}\\ \notag
    &\psi_j(z) := (1 + \mathbb{E}\operatorname{tr}(\mathcal{B}_j))^{-1} = (1 + c_n\boldsymbol{\theta}_{j}^\top \mathbb{E}\textbf{h}_n(z))^{-1}, \quad b_j^0(z) = (1 + c_n\boldsymbol{\theta}_{j}^\top \textbf{h}_n^0(z))^{-1}.
\end{align}
where, $\textbf{h}_n(z)$ and $\textbf{h}_n^0(z)$ were defined in (\ref{defining_hn}) and (\ref{finite_sample_properties}) respectively. Here, 
$\beta_j$ and $\widetilde{\beta}_j$ are the random LOO deflators which will appear in the martingale style decomposition (\ref{decomposing_Mn1}) of $M_{n1}(z)$. $b_j$,$\psi_j$ and $b_j^0$ are their deterministic counterparts and we will use these to gradually establish limiting quantities. 

A very important fact we establish in Lemma \ref{lemma:derivative_quadratic_forms} is that
\begin{align}\label{derivative_quadratic_forms}
    \frac{d}{dz}\big(\epsilon_j(z)\big) = \gamma_j(z).
\end{align}
We start with some basic bounds for a few important quantities. The proof can be found in Section \ref{sec:proofLemma4.1}.
\begin{lemma}\label{lemma:Bounds_Important_Quantities}
    For any $z \in \mathbb{C}_+$, the following hold:
    $$\max\{|\beta_j(z)|, |\widetilde{\beta}_j(z)|, |b_j(z)|, |\psi_j(z)|\} \leq |z|/v, \quad \max\{\|R(z)\|, \|R_j(z)\| \}\leq 1/v.$$
\end{lemma}

These results help only for $z$ with large imaginary component (say $\Im(z) \geq v_0 > 0$). A particular issue which arises is on $z \in \mathcal{C}_v := \mathcal{C}_l \cup \mathcal{C}_r$, is the instability of resolvent type matrices as we approach the real axis. Since our truncated contour becomes arbitrarily close to the real line, we will need results sharper than those established in Lemma \ref{lemma:Bounds_Important_Quantities}. 

Depending on the locations of the random eigenvalues of $S_n/S_{nj}$ with respect to the interval $[\eta_l, \eta_r]$ defined in (\ref{eta_l_r}), we define the \textit{good} events to be 
\begin{align}\label{defining_good_event}
    \mathcal{G} &:= \{\eta_l \leq \lambda_{\min}(S_n) \leq \lambda_{\max}(S_n) \leq \eta_r\},\\ \notag
    \mathcal{G}_j &:= \{\eta_l \leq \lambda_{\min}(S_{nj}) \leq \lambda_{\max}(S_{nj}) \leq \eta_r\}, \quad j \in [n]\\ \notag
    \mathcal{G}_0 &:= \bigcap_{j=1}^n  \mathcal{G}_j,
\end{align}
and naturally the \textit{bad} events to be their complements. From (\ref{extreme_eigenvalues}), for any $l \in \mathbb{N}$, we therefore can establish the following:
\begin{align}\label{bad_event_probability}
    \mathbb{P}(\mathcal{G}^c) = o(n^{-l}), \quad \mathbb{P}( \mathcal{G}_0^c) = o(n^{-l}).
\end{align}
Another quantity which we shall frequently encounter is a bound on the norms of resolvent type objects on events of the \textit{good} nature.
\begin{align}{\label{resolvent_norm_bound_good_event}}
    M_{good} := \max\bigg \{\frac{1}{|x_r - \eta_r|}, \frac{1}{|\eta_l - x_l|}\bigg\}.
\end{align}
Note that by our choice of $[\eta_l, \eta_r]$, we have $M_{good} < 1/4Kc\tau^2$. We see that 
$$\sup_{n:z \in \mathcal{C}_n}\|Q(z)\| \leq M_{good}, \quad \sup_{n; z \in \mathcal{C}_n} \max_{j \in [n]}\|Q_j(z)\| \leq M_{good},$$ on $\mathcal{G}$ and $\mathcal{G}_0$ respectively. On the bad events, we will use (\ref{bad_event_probability}) in conjunction with suitable bounds as per the quantity we are dealing with. 

We now state a few uniform across $\mathcal{C}_n$ bounds that will be used throughout the paper.

\begin{lemma}\label{lemma:bounds_tightness_section}
The following are some uniform bounds over $z \in \mathcal{C}_n$ for the norms of resolvent-type objects, LOO deflators and centered quadratic forms defined in (\ref{essential_quantities}).
    \begin{enumerate}
        \item $\underset{n,z \in \mathcal{C}_n} {\sup}\mathbb{E}\|Q_*(z)\|^k < \infty$, where $Q_*$ can be $Q, Q_j$ or $Q_{ij}$, $i,j \in [n], i \neq j$.
        As a result for non-random $M$ with $\|M\|<\infty$, we have
        $$\sup_{z\in\mathcal{C}_n}\mathbb{E}|\frac{1}{n}\operatorname{tr}MQ_*(z)| < \infty.$$
        \item $\underset{n, j \in [n],z \in \mathcal{C}_n}{\sup}|\psi_j(z)| < \infty, \quad \underset{n,j \in [n], z \in \mathcal{C}_n}{\sup}|b_j(z)| < \infty$.
        \item $\underset{n, j \in [n],z \in \mathcal{C}_n}{\sup}\mathbb{E}|\beta_j(z)|^k \leq C_k, \quad \underset{n, j \in [n],z \in \mathcal{C}_n}{\sup} \mathbb{E}|\widetilde{\beta}_j(z)|^k \leq C_k$.
        \item $\underset{n, j \in [n],z \in \mathcal{C}_n}{\sup}\mathbb{E}|\epsilon_j(z)|^k = O(\delta_n^{2k-4}n^{-1}), \quad \underset{n, j \in [n],z \in \mathcal{C}_n}{\sup}\mathbb{E}|\gamma_j(z)|^k = O(\delta_n^{2k-4}n^{-1}).$
        
        If in addition, the innovation entries are i.i.d. Gaussian in nature, we have 
        $$\underset{n, j \in [n],z \in \mathcal{C}_n}{\sup}\mathbb{E}|\epsilon_j(z)|^k = O(n^{-\frac{k}{2}}), \quad \underset{n, j \in [n],z \in \mathcal{C}_n}{\sup}\mathbb{E}|\gamma_j(z)|^k = O(n^{-\frac{k}{2}}).$$
    \end{enumerate}
\end{lemma}
\begin{lemma}\label{lemma:bound_centered_trace_resolvent_squared}
Let $M$ be a $p \times p$ non-random matrix. For any $j \in [n]$, we have $$\sup_{z\in\mathcal{C}_n}\mathbb{E}\absmod{\operatorname{tr}(MQ(z)) - \mathbb{E}\operatorname{tr}MQ(z)}^2 \leq C\|M\|^2.$$
\end{lemma}
\begin{lemma}\label{lemma:hn_Ehn_hn0_uniform_results} Concentration of empirical and finite sample solutions:
     \begin{enumerate}
        \item $\sup_{z\in \mathcal{C}_n}\|\textbf{h}_n(z) - \mathbb{E}\textbf{h}_n(z)\| \xrightarrow{a.s.} 0.$
        \item $\sup_{z\in\mathcal{C}}\|\textbf{h}_n^0(z) - \textbf{h}(z)\| \to 0.$
\end{enumerate}
\end{lemma}
\begin{lemma}\label{lemma:exp_betaj_psij_bj_difference}
Concentration of LOO deflators:
\begin{enumerate}
    \item $\underset{n, j \in [n],z \in \mathcal{C}_n}{\sup} |b_j(z) - \psi_j(z)| = O(n^{-1}).$
    \item $\underset{n, j \in [n],z \in \mathcal{C}_n}{\sup}\mathbb{E}|\beta_j(z) - \psi_j(z)| = O(n^{-1})$.
\end{enumerate}
\end{lemma}
\begin{lemma}\label{lemma:bounds_bias_section}
Recall a few definitions from (\ref{essential_quantities1}). An additional set of bounds for $z \in \mathcal{C}_n$.
    \begin{enumerate}
        \item $\underset{n;z \in \mathcal{C}_n}{\sup}\|\mathbb{E}\textbf{h}_n(z) - \textbf{h}_n^0(z)\| \longrightarrow 0.$
        \item $\underset{n;j \in [n], z \in \mathcal{C}_n}{\sup}|b_j^0(z)| < \infty, \quad \underset{n;j \in [n], z \in \mathcal{C}_n}{\sup}|\psi_j(z) - b_j^0(z)| = o(1)$. 
    \end{enumerate}
\end{lemma}

\section{Finite Dimensional Convergence of \texorpdfstring{$\{M_{n1}\}$}{}}\label{sec:GaussianCase}
The aim of this section is to find the limiting distribution of $M_{n1}(z)$ or equivalently to show that for arbitrary $k \in \mathbb{N}$ and $a_1, \ldots, a_k, z_1, \ldots, z_k \in \mathbb{C}$, there exists a deterministic sequence $\{\sigma_n\}$ such that the quantity defined below
\begin{align}\label{finite_dimensional_sum}
    \sigma_n^{-1}\sum_{j=1}^k a_j M_{n1}(z_j), \quad \Im(z) \neq 0,
\end{align}
converges weakly to a standard Gaussian random variable. 
\begin{remark}\label{covariance_horizontal_only}
We point out that for the covariance function, it suffices to consider only $z = u + \mathbbm{i}v_0 \in \mathcal{C}_u$. This follows from the fact (see \ref{proof_vanishing_contribution_vertical_arms} for proof ) that
\begin{align}\label{vanishing_contribution_vertical_arms}
    \underset{v_0 \downarrow 0}{\lim}\, \underset{n \rightarrow \infty}{\limsup}\,\mathbb{E}\absmod{\int_{\mathcal{C}_l \cup \mathcal{C}_r} f(z)M_{n1}(z)dz}^2 \longrightarrow 0,
\end{align}
\end{remark}

Using the property (\ref{Leave_one_out_property}), for fixed $z \in \mathbb{C}_+$, we have
\begin{align}\label{decomposing_Mn1}
    M_{n1}(z) &=
    \operatorname{tr}(Q(z)) - \mathbb{E}\operatorname{tr}(Q_j(z))\\ \notag
    &=\sum_{j=1}^n(\mathbb{E}_j-\mathbb{E}_{j-1})\operatorname{tr}(Q(z) - \operatorname{Q}_j(z))\\ \notag
    &= -\sum_{j=1}^n \Delta_j\beta_j(z)x_j^*Q_j^2(z)x_j\\ \notag
    &= -\sum_{j=1}^n \Delta_j \beta_j(z) (\gamma_j(z) + \alpha_j(z)) := I_1 + I_2.
\end{align}

For simplicity, we will avoid the argument ($z \in \mathbb{C}$) in quantities such as $\beta_j, \widetilde{\beta}_j, \epsilon_j, \gamma_j,\alpha_j$ whenever the context is clear. Using the identity,
\begin{align}\label{beta_betaTilde_relationship}
    \beta_j - \widetilde{\beta}_j 
    &= -\beta_j\widetilde{\beta}_j\epsilon_j
    = -\widetilde{\beta}_j^2\epsilon_j + \beta_j\widetilde{\beta}_j ^2\epsilon_j^2,
\end{align}
we simplify $I_1$ of (\ref{decomposing_Mn1}) to get
\begin{align}\label{decomposing_I1}
    I_1 =  -\sum_{j=1}^n \mathbb{E}_j \widetilde{\beta}_j \gamma_j -\sum_{j=1}^n \Delta_j \beta_j\widetilde{\beta}_j \epsilon_j\gamma_j.
\end{align}
Here we used the fact that 
\begin{align}\label{tower_property_usage1}
    \mathbb{E}_{j-1}(\widetilde{\beta}_j\gamma_j) = \mathbb{E}_{j-1}(\mathbb{E}_{-j}(\widetilde{\beta}_j\gamma_j)) = \mathbb{E}_{j-1}(\widetilde{\beta}_j\mathbb{E}_{-j}\gamma_j)= \mathbb{E}_{j-1}(\widetilde{\beta}_j\cdot 0) = 0.
\end{align}
Note that the random quantities $\{\Delta_j\beta_j\widetilde{\beta}_j \epsilon_j \gamma_j\}_{j=1}^n$ are pairwise uncorrelated. This can be explained as follows. Let $T_j = \beta_j\widetilde{\beta}_j \epsilon_j \gamma_j$. Let $j \neq k$. Without loss of generality, assume $j < k$. Then, we have
\begin{align*}
    \operatorname{Cov}(\Delta_jT_j, \Delta_kT_k)&= 
    \mathbb{E} (\Delta_jT_j)(\Delta_kT_k) = \mathbb{E}\bigg(\mathbb{E}_j(\Delta_jT_j)(\Delta_kT_k)\bigg) = 0,
\end{align*}
where the last equality follows by properties of martingale difference sequences. Using the above fact, Lemma \ref{lemma:bounds_tightness_section} and Cauchy-Schwarz, we bound the second term of $I_1$ as follows:
\begin{align*}
    \mathbb{E}\absmod{\sum_{j=1}^n \Delta_j\beta_j\widetilde{\beta}_j \epsilon_j \gamma_j}^2
    =&\, \mathbb{E} \sum_{j=1}^n \absmod{\Delta_j\beta_j\widetilde{\beta}_j \epsilon_j \gamma_j}^2\\
    \leq&\, 4 \sum_{j=1}^n \mathbb{E}\absmod{\beta_j\widetilde{\beta}_j \epsilon_j \gamma_j}^2\\
    \leq&\, C\sum_{j=1}^n \bigg(\mathbb{E} |\epsilon_j|^8 \mathbb{E} |\gamma_j|^8\bigg)^\frac{1}{4} = O(n^{-1}).
\end{align*}
Therefore the second term in the expansion of $I_1$ in (\ref{decomposing_I1}) converges to 0 in probability and thus
\begin{align}\label{final_expression_I1}
    I_1 = -\sum_{j=1}^n \mathbb{E}_j \widetilde{\beta}_j \gamma_j  + o_P(1).
\end{align}
Now we focus on $I_2$ of (\ref{decomposing_Mn1}). Using property (\ref{beta_betaTilde_relationship}) we get
\begin{align*}
    -I_2 =&\sum_{j=1}^n \Delta_j\beta_j\alpha_j\\
    =& \sum_{j=1}^n\Delta_j \widetilde{\beta}_j\alpha_j - \sum_{j=1}^n\Delta_j \widetilde{\beta}_j^2\epsilon_j\alpha_j + \sum_{j=1}^n\Delta_j \widetilde{\beta}_j^2\beta_j\epsilon_j^2\alpha_j.
\end{align*}
Note that $\widetilde{\beta}_j$ and $\alpha_j$ do not depend on $x_j$. Therefore, we have
\begin{align*}
    \sum_{j=1}^n\Delta_j \widetilde{\beta}_j\alpha_j=0.
\end{align*}
Using arguments similar as those in (\ref{tower_property_usage1}), we have
\begin{align}\label{tower_property_usage2}
     \sum_{j=1}^n\mathbb{E}_{j-1}\widetilde{\beta}_j^2\epsilon_j\alpha_j=0.
\end{align}
Repeating the arguments that led to (\ref{final_expression_I1}), we also have
\begin{align*}
    \sum_{j=1}^n\Delta_j \widetilde{\beta}_j^2\beta_j\epsilon_j^2\alpha_j = o_P(1).
\end{align*}
Combining all the terms gives us
\begin{align}\label{final_expression_I2}
    I_2 = -\sum_{j=1}^n\mathbb{E}_j\widetilde{\beta}_j^2\epsilon_j\alpha_j + o_P(1).
\end{align}
Finally from (\ref{decomposing_Mn1}), (\ref{final_expression_I1}), (\ref{final_expression_I2}) and (\ref{derivative_quadratic_forms}), we have
\begin{align}\label{final_expression_Mn1}
    M_{n1}(z) =& -\sum_{j=1}^n \mathbb{E}_j\bigg(\widetilde{\beta}_j \gamma_j + \widetilde{\beta}_j^2\epsilon_j \alpha_j\bigg) + o_P(1)\\\notag
    :=&\, \sum_{j=1}^n h_j(z) +o_P(1)\\\notag
    =& -\sum_{j=1}^n\mathbb{E}_j\bigg(\frac{d}{dz}  \widetilde{\beta}_j \epsilon_j\bigg)  + o_P(1).
\end{align}
From (\ref{tower_property_usage1}) and (\ref{tower_property_usage2}), it follows that $\mathbb{E}_{j-1}h_j(z) = 0$ and $\{h_j(z)\}_{j=1}^n$ is therefore a martingale difference sequence. 
We first proceed to verify condition (ii) of Theorem \ref{theorem:martingale_CLT}. For $t \in \mathbb{N}$, $\textbf{a} = (a_1,\ldots,a_t)$, $\textbf{z} =(z_1,\ldots,z_t)$, we define the quantity 
$$J_n(\textbf{a}, \textbf{z}) := \sum_{k=1}^t a_k \sum_{j=1}^n h_j(z_k) = \sum_{j=1}^n \sum_{k=1}^t a_k h_j(z_k).$$ 
From Lemma \ref{lemma:bounds_tightness_section}, for some constant $C>0$, we have
\begin{align}\label{Mn1_4th_moment_bound}
    \sum_{j=1}^n \mathbb{E}\absmod{\sum_{k=1}^ta_k h_j(z_k)}^4 \leq C\sum_{j=1}^n \sum_{k=1}^t a_j^4\bigg(C_1\bigl(\mathbb{E} |\gamma_j(z_k)|^8\bigr)^\frac{1}{2} + C_2\bigl(\mathbb{E}|\epsilon_j(z_k)|^8\bigr)^\frac{1}{2}\bigg) = O(n^{-1}).
\end{align}
This implies that for arbitrary $\epsilon > 0$, we have
\begin{align*}
    \sum_{j=1}^n\mathbb{E}\bigg(|J_n(\textbf{a}, \textbf{z})|^2 I\bigg(|J_n(\textbf{a}, \textbf{z})| > \epsilon\bigg)\bigg) \longrightarrow 0.
\end{align*}
\subsection{Deriving the Covariance function}
We substitute the martingale representation (\ref{final_expression_Mn1}) into the finite dimensional sum in (\ref{finite_dimensional_sum}). Thus verifying condition (i) of Theorem \ref{theorem:martingale_CLT} reduces to establishing the existence of a deterministic sequence $\{\sigma_n\}$ such that 
$$\sigma_n^{-2}\sum_{j=1}^n \mathbb{E}_{j-1} \bigg|\sum_{i=1}^t a_ih_j(z_i)\bigg|^2 \xrightarrow{\,\,\, p\,\,} 1.$$
Equivalently, for fixed but arbitrary $z_1, z_2$, we will establish that there exists deterministic $\{\sigma_n\}$ such that 
\begin{align}\label{condition_ii}
    \sigma_n^{-2}\Phi_n(z_1, z_2) = \sigma_n^{-2}\sum_{j=1}^n\mathbb{E}_{j-1}[h_j(z_1)h_j(z_2)] \xrightarrow{\,\,\, p\,\,} 1.
\end{align}
We derive a simplified version of the above expression that utilizes the derivatives. From (\ref{final_expression_Mn1}) and utilizing the boundedness and integrability results established in Lemma \ref{lemma:bounds_tightness_section}, we can apply Dominated Convergence Theorem to get 
\begin{align}\label{ip_limit_simplification1}
    \Phi_n(z_1, z_2) = \frac{\partial^2}{\partial z_2 \partial z_1} \sum_{j=1}^n \mathbb{E}_{j-1} \bigg(\mathbb{E}_j(\widetilde{\beta}_j(z_1)\epsilon_j(z_1))\mathbb{E}_j(\widetilde{\beta}_j(z_2)\epsilon_j(z_2))\bigg).
\end{align}
By Lemma \ref{lemma:transfer_of_derivative}, our object of interest thus becomes
\begin{align}\label{ip_limit_simplification2}
   f_n(z_1, z_2) := \sum_{j=1}^n \mathbb{E}_{j-1} \bigg(\mathbb{E}_j(\widetilde{\beta}_j(z_1)\epsilon_j(z_1))\mathbb{E}_j(\widetilde{\beta}_j(z_2)\epsilon_j(z_2))\bigg).
\end{align}
Similar to $f_n(z_1, z_2)$, we construct
\begin{align}\label{ip_limit_simplification3}
    g_n(z_1, z_2) : = \sum_{j=1}^n b_j(z_1)b_j(z_2)\mathbb{E}_{j-1} \bigg(\mathbb{E}_j(\epsilon_j(z_1))\mathbb{E}_j(\epsilon_j(z_2))\bigg).
\end{align}
From Lemma \ref{lemma:bounds_tightness_section}, we obtain the following
$$\sup_{z\in\mathcal{C}_n}\mathbb{E}|\widetilde{\beta}_j(z) - b_j(z)|^2 \leq C(\mathbb{E}|\widetilde{\beta}_j(z)|^4)^\frac{1}{2}(\mathbb{E}|\epsilon_j(z)|^4)^\frac{1}{2}= O(n^{-1}).$$
Using the above, it can be shown that $$|f_n(z_1, z_2) - g_n(z_1, z_2)| \xrightarrow{p} 0,$$
and thus (\ref{ip_limit_simplification3}) becomes our expression of interest. From Lemma \ref{lemma:exp_betaj_psij_bj_difference}, it is clear that $b_j$ approaches $b_j^0$ uniformly. We therefore, elaborate the expression  $\mathbb{E}_j(\epsilon_j(z_1))\mathbb{E}_j(\epsilon_j(z_2))$. 
\begin{claim}\label{claim:standardized_vector}
    Let $x$ be a random vector with $\mathbb{E}(x) = \textbf{0}$ and $\operatorname{Var}(x) = \textbf{A}$.  $\textbf{A}^\frac{1}{2}$ represents the p.s.d. (not necessarily p.d.) square root of $\textbf{A}.$ Then, there exists a random vector $y$ with mean $\textbf{0}$ and variance $I_p$ such that $$x = \textbf{A}^\frac{1}{2}y.$$
\end{claim} 
Recall that $x_j$ has mean \textbf{0} and variance $n^{-1}\textbf{A}_j$ for $j \in [n]$. For each $x_j$, construct $y_j$ as per the above claim.
Then, it is clear that 
$$x_j^*Q_jx_j = \frac{1}{n}y_j^*\textbf{A}_j^\frac{1}{2}Q_j\textbf{A}_j^\frac{1}{2}y_j = y_j^*\mathcal{M}_jy_j,$$ where
$\mathcal{M}_j(z)$ is a symmetric analog of $\mathcal{A}_j(z)$ (see \ref{defining_script_Aj}) defined as:
\begin{align}\label{defining_Mj_Nj}
    \mathcal{M}_j(z) = \frac{1}{n}\textbf{A}_j^\frac{1}{2}Q_j(z)\textbf{A}_j^\frac{1}{2}.
\end{align}

Since $\operatorname{tr}(\mathcal{M}_j(z)) = \operatorname{tr}(\mathcal{A}_j(z))$, by Lemma \ref{lemma:bounds_tightness_section}, $y_j^*\mathcal{M}_j(z)y_j$ and $y_j^*\mathcal{A}_j(z)y_j$ concentrate around the same quantity which follows from 
\begin{align*}
    &\mathbb{E}\absmod{y_j^*\mathcal{A}_j(z)y_j - y_j^*\mathcal{M}_j(z)y_j}^2 \\
    \leq& \,\, 2\mathbb{E}|y_j^*\mathcal{A}_j(z)y_j - \operatorname{tr}\mathcal{A}_j(z)|^2 + 2\mathbb{E}|y_j^*\mathcal{M}_j(z)y_j - \operatorname{tr}\mathcal{M}_j(z)|^2 = O(n^{-1}).
\end{align*}
Therefore, $x_j^*Q_j(z)x_j = y_j^*\mathcal{M}_j(z)y_j + o_p(n^{-\omega})$ for any $\omega \in [0, 1/2)$. Define 
\begin{align}
    \tilde{\epsilon}_j(z) := y_j^*\mathcal{A}_j(z)y_j - \operatorname{tr}\mathcal{A}_j(z).
\end{align}
Then, $\epsilon_j(z) = \widetilde{\epsilon}_j(z) + o_P(1)$ and by uniform integrability of $\epsilon_j(\cdot)$ and $\widetilde{\epsilon}(\cdot)$, we have
$$\mathbb{E}_{j-1}[\mathbb{E}_j(\epsilon_j(z_1))\mathbb{E}_j(\epsilon_j(z_2))] - \mathbb{E}_{j-1}[\mathbb{E}_j(\widetilde{\epsilon}_j(z_1))\mathbb{E}_j(\widetilde{\epsilon}_j(z_2))] \overset{p}{\longrightarrow} 0.$$
Noting that $Q_j(z_i)$ is free of $x_j$, applying Lemma \ref{lemma:covariance_quadratic_forms} gives us
\begin{align}
    &\mathbb{E}_{j-1}[\mathbb{E}_j(\tilde{\epsilon}_j(z_1))\mathbb{E}_j(\tilde{\epsilon}_j(z_2))]\\ \notag
    =& \mathbb{E}_{j-1}[(y_j^*(\mathbb{E}_{j-1}\mathcal{A}_j(z_1))y_j - \operatorname{tr}\mathcal{A}_j(z_1))(y_j^*(\mathbb{E}_{j-1}\mathcal{A}_j(z_2))y_j - \operatorname{tr}\mathcal{A}_j(z_2))]\\ \notag
    =& \operatorname{tr}\mathbb{E}_{j-1}(\mathcal{A}_j(z_1)\mathcal{A}_j(z_2)) + o(n^{-1})\\ \notag
    =& \frac{1}{n^2}\sum_{r,s=1}^K \theta_{rj}\theta_{sj}\operatorname{tr}(A_r\mathbb{E}_jQ_j(z_1)A_s\mathbb{E}_jQ_j(z_2)) + o(n^{-1}).
\end{align}
The above expansion is valid when the innovations are complex. For the real case, it follows from Lemma \ref{lemma:covariance_quadratic_forms} that
$$\mathbb{E}_{j-1}[\mathbb{E}_j(\epsilon_j(z_1))\mathbb{E}_j(\epsilon_j(z_2))]
    =\frac{2}{n^2}\sum_{r,s=1}^K \theta_{rj}\theta_{sj}\operatorname{tr}(A_r\mathbb{E}_jQ_j(z_1)A_s\mathbb{E}_jQ_j(z_2)) + o\bigg(\frac{1}{n}\bigg)+ o_p(n^{-2\omega}).$$
Therefore, our goal now is to prove the existence of $\{\sigma_n\}$ such that
\begin{align}\label{revised_goal}
    \frac{1}{n^2\sigma_n^2} \sum_{j=1}^n b_j(z_1)b_j(z_2)\bigg(\sum_{r,s=1}^K \theta_{rj}\theta_{sj}\operatorname{tr}(A_r\mathbb{E}_jQ_j(z_1)A_s\mathbb{E}_jQ_j(z_2))\bigg) \longrightarrow 1.
\end{align}
So all that needs to be done is to develop a simplified expression for $$\operatorname{tr}(A_r\mathbb{E}_jQ_j(z_1)A_r\mathbb{E}_jQ_j(z_2))$$ in the limit. To this end, we introduce the following matrix-valued function:
\begin{align}\label{defining_Vjrs}
    V_j(z_1, z_2)[r,s]:=\operatorname{tr} (A_r\mathbb{E}_jQ_j(z_1)A_s\mathbb{E}_jQ_j(z_2)), \quad r,s \in [K].
\end{align}
Henceforth, we will use the $[r,s]$ notation to indicate the $rs-$th element of matrices. For simplicity, we will ignore the arguments $z_1,z_2$ from matrix-valued functions such as $V_j$ and other similar functions to be defined later. So the $rs-$th element of $V_j$ shall be denoted by $V_j[r,s]$. 

Now denote $(\Theta_j)_{(r,s)} = \theta_{rj}\theta_{sj}, r,s \in [K]$ and $V_j(z_1,z_2)$ as the $K\times K$ matrix $\{V_j[r,s]\}_{r,s=1}^K$. By the symmetry of $\Theta_j$, we can express (\ref{revised_goal}) as follows:
\begin{align}\label{revised_goal1}
    \frac{1}{\sigma_n^2} \operatorname{tr}\bigg(\frac{1}{n^2}\sum_{j=1}^nb_j(z_1) b_j(z_2)\Theta_jV_j(z_1, z_2)\bigg) \overset{p}{\longrightarrow} 1.
\end{align}
By Lemma \ref{lemma:transfer_of_derivative}, the second mixed partial derivative of the limit in probability will yield the limit of 
$$\operatorname{Cov}(M_{n1}(z_1), M_{n1}(z_2)).$$

\subsection{Convergence of \texorpdfstring{$V_j[r,s]$}{}}

Similar to LOO resolvents and deflators, we now introduce their leave-two-out (LTO) versions. Recalling that $x_i$ denotes the $i^{th}$ column of $n^{-\frac{1}{2}}X$, we introduce these quantities for $i \neq j, j \in [n]$.
\begin{align}\label{essential_quantities1}
    &Q_{ij}(z) := (S_n - x_ix_i^* - x_jx_j^* - zI_p)^{-1},\\\notag
    &\beta_{ij}(z) := \big(1 + x_{i}^*Q_{ij}(z)x_{i}\big)^{-1}, \quad
    b_{ij}(z) := \big(1 + n^{-1}\mathbb{E}\operatorname{tr}\textbf{A}_iQ_{ij}(z)\big)^{-1}\\\notag
    & R_j(z) := \bigg(\frac{1}{n}\sum_{i \neq j}\psi_i(z)\textbf{A}_i - zI_p\bigg)^{-1}, \quad  R(z) := \bigg(\frac{1}{n}\sum_{i =1}^n\psi_i(z)\textbf{A}_i - zI_p\bigg)^{-1}.
\end{align}
The main idea here is to sequentially substitute the random components in $V_j[r,s]$, i.e. $Q_j(z_1)$ and $Q_j(z_2)$ with their deterministic counterparts $R_j(z_1)$ and $R_j(z_2)$ and compute the non-vanishing terms. Noting that $x_i^*Q_j(z) = \beta_{ij}x_i^*Q_{ij}(z)$, we use the resolvent identity (\ref{Resolvent_Identity}) to split $Q_j(z)$ into deterministic and random fluctuation terms as follows:
\begin{align}\label{second_order_expansion_Qj}
    Q_j(z) - R_j(z) &= R_j(z)(R_j^{-1}(z) - Q_j^{-1}(z)) Q_j(z)\\ \notag
    &= R_j(z)\bigg(\frac{1}{n}\sum_{i\neq j}\psi_i(z)\textbf{A}_i - \sum_{i\neq j}x_ix_i^*\bigg)Q_j(z)\\ \notag
    &= -\sum_{i \neq j} (\beta_{ij}(z) - \psi_i(z))R_j(z)x_ix_i^*Q_{ij}(z)\\ \notag
    &- \sum_{i \neq j}\psi_i(z)R_j(z)(x_ix_i^* - n^{-1}\textbf{A}_i)Q_{ij}(z)\\ \notag
    &- \frac{1}{n}\sum_{i \neq j}\psi_i(z) R_j(z)\textbf{A}_i(Q_{ij}(z) - Q_j(z))\\ \notag
    &:= D_j(z) + C_j(z) + E_j(z).
\end{align}
Replacing $Q_j(z_i)$ with $R_j(z_i)$ in the expression for $V_j[r, s]$ will lead to random fluctuations involving the terms $C_j(z)$, $D_j(z)$ and $E_j(z)$ since $R_j(z)$ is a deterministic object. We first establish a few bounds which will show that only the terms involving $C_j(z)$ lead to non-vanishing contribution in the limit.
\begin{lemma}\label{lemma:tracial_bounds}
    Under Assumptions, for any $z \in \mathbb{C}_0$, the following hold:
    \begin{enumerate}
        \item For any (possibly random) matrix $M$ with a deterministic bound on its spectral norm,
        $$\mathbb{E}|\operatorname{tr}D_j(z)M| = O(n^\frac{1}{2}), \quad |\operatorname{tr}E_j(z)M| = O(1).$$

        \item For any deterministic matrix $M$ with bounded spectral norm, 
        $$\mathbb{E}|\operatorname{tr}C_j(z)M| = O(n^\frac{1}{2}).$$
    \end{enumerate}
\end{lemma}
The proof of this lemma is in Section \ref{sec:proofLemma4.2}. Having proved $\|R_j(z)\|\leq v_0^{-1}$ in Lemma \ref{lemma:Bounds_Important_Quantities}, we use (\ref{second_order_expansion_Qj}) in the expression for $V_j[r,s]$ from (\ref{defining_Vjrs}). For $l =1,2$, we will use $\beta_j^{(l)}$, $\widetilde{\beta}_j^{(l)}$, $\psi_j^{(l)}$, $C_j^{(l)}$, $Q^{(l)}$, $Q_j^{(l)}$, $Q_{ij}^{(l)}$ and $R_j^{(l)}$ to denote $\beta_j(z_l)$,$\widetilde{\beta}_j(z_l)$, $\psi_j(z_l)$, $C_j(z_l)$, $Q(z_l)$, $Q_j(z_l)$, $Q_{ij}(z_l)$ and $R(z_l)$ respectively. Fix $r_0,s_0 \in [K]$. Implementing the bounds established in Lemma \ref{lemma:tracial_bounds}, we get the following:
\begin{align*}
    V_j[r_0, s_0] =& \operatorname{tr}(A_{r_0}\mathbb{E}_jQ_j^{(1)}A_{s_0}\mathbb{E}_jQ_j^{(2)})\\
    =& \operatorname{tr}(A_{r_0}R_j^{(1)}A_{s_0}\mathbb{E}_jQ_j^{(2)}) +  \operatorname{tr}(A_{r_0}\mathbb{E}_jC_j^{(1)}A_{s_0}\mathbb{E}_jQ_j^{(2)}) + a(z_1, z_2),
\end{align*}
where the term $a(z_1, z_2)$\footnote{Throughout this section, we will use $a(z_1, z_2)$ as a wrapper for expressions which satisfying $\mathbb{E}|a(z_1, z_2)|\leq O(n^\frac{1}{2})$. Similarly, $A(z_1, z_2)$ will be used to represent a $K \times K$ complex matrix which satisfies the previous condition entry-wise. Each occurrence of $a(\cdot, \cdot)$ and $A(\cdot, \cdot)$ may represent different objects.} accumulates all the random fluctuations arising out of $D_j(z_1)$ and $ E_j(z_1)$ and therefore satisfies $\mathbb{E}|a(z_1, z_2)| \leq O(n^\frac{1}{2})$ by Lemma \ref{lemma:tracial_bounds}. Using the expression for $Q_j(z_2)$ above, we have
\begin{align}\label{Vjrs_simplification}
    V_j[r_0,s_0] = \operatorname{tr}(A_{r_0}R_j^{(1)}A_{s_0}R_j^{(2)}) +  \operatorname{tr}(A_{r_0}\mathbb{E}_jC_j^{(1)}A_{s_0}\mathbb{E}_jQ_j^{(2)}) + a(z_1, z_2),
\end{align}
where $a(z_1, z_2)$ collects more random fluctuation terms, this time from $D_j(z_2)$ and $E_j(z_2)$. Now we focus on the term $\operatorname{tr}(A_{r_0}\mathbb{E}_jC_j^{(1)}A_{s_0}\mathbb{E}_jQ_j^{(2)})$. We will use $C_j(z_l)$ and $C_j^{(l)}$ for $l=1,2$ interchangeably as space permits. Observe that for $i > j$, the following holds:
\begin{align}
\mathbb{E}_j[\psi_i(z)R_j(z)(x_ix_i^* - n^{-1}\textbf{A}_i)Q_{ij}(z)] = \psi_i(z)R_j(z)\mathbb{E}_j(x_ix_i^* - n^{-1}\textbf{A}_i)Q_{ij}(z) = 0.
\end{align}
Therefore, we have
$$\mathbb{E}_jC_j(z_1) = -\sum_{i < j} \mathbb{E}_j \psi_i(z_1)R_j(z_1)(x_ix_i^* - n^{-1}\textbf{A}_i)Q_{ij}(z_1),$$
which leads to the following expansion:
\begin{align}\label{tracial_expansion_covariance}
    &\operatorname{tr}(A_{r_0}\mathbb{E}_jC_j^{(1)}A_{s_0}\mathbb{E}_jQ_j^{(2)})\\ \notag
    =& -\sum_{i<j} \psi_i^{(1)}\operatorname{tr}R_j^{(1)}(x_ix_i^* - n^{-1}\textbf{A}_i)\mathbb{E}_jQ_{ij}^{(1)}A_{s_0}\mathbb{E}_jQ_j^{(2)}A_{r_0}\\ \notag
    =&-\sum_{i<j} \psi_i^{(1)}\operatorname{tr}R_j^{(1)}(x_ix_{i}^* - n^{-1}\textbf{A}_i)\mathbb{E}_jQ_{ij}^{(1)}A_{s_0}\mathbb{E}_jQ_{ij}^{(2)}A_{r_0} \\ \notag
    -&\sum_{i<j} \psi_i^{(1)}\operatorname{tr}R_j^{(1)}(x_ix_i^* - n^{-1}\textbf{A}_i)\mathbb{E}_jQ_{ij}^{(1)}A_{s_0}\mathbb{E}_j(Q_j^{(2)} - Q_{ij}^{(2)})A_{r_0}\\ \notag
    :=& T_1 + T_2,
\end{align}
where,
\begin{align*}
    T_1 = -\sum_{i<j} \psi_i^{(1)}\bigg(x_i^*\mathbb{E}_jQ_{ij}^{(1)}A_{s_0}\mathbb{E}_jQ_{ij}^{(2)}A_{r_0}R_j^{(1)}x_i - \frac{1}{n}\operatorname{tr}\textbf{A}_i\mathbb{E}_jQ_{ij}^{(1)}A_{s_0}\mathbb{E}_jQ_{ij}^{(2)}A_{r_0}R_j^{(1)}\bigg),
\end{align*}
and
\begin{align*}    
    T_2=& -\sum_{i<j} \psi_i^{(1)}x_i^*\mathbb{E}_jQ_{ij}^{(1)}A_{s_0}\mathbb{E}_j(Q_j^{(2)}-Q_{ij}^{(2)})A_{r_0}R_j^{(1)}x_i\\
    &+ \frac{1}{n}\sum_{i<j} \psi_i^{(1)}\operatorname{tr}\textbf{A}_i\mathbb{E}_jQ_{ij}^{(1)}A_{s_0}\mathbb{E}_j(Q_j^{(2)}-Q_{ij}^{(2)})A_{r_0}R_j^{(1)}\\
    :=&T_3 + T_4.
\end{align*}

We have $\mathbb{E}|T_1| \leq O(n^\frac{1}{2})$ and $\mathbb{E}|T_4| \leq O(1)$. All that remains is the term $T_3$. Now we create another term $T_5$ such that $\mathbb{E}_jT_5 = T_3$ as follows:
\begin{align}
    T_5:=& -\sum_{i<j} \psi_i^{(1)}x_i^*\mathbb{E}_jQ_{ij}^{(1)}A_{s_0}(Q_j^{(2)}-Q_{ij}^{(2)})A_{r_0}R_j^{(1)}x_i\\ \notag
=& \sum_{i<j} \psi_i^{(1)}\beta_{ij}^{(2)}\underbrace{x_i^*\mathbb{E}_jQ_{ij}^{(1)}A_{s_0}Q_{ij}^{(2)}x_i}_{qf_1}\underbrace{x_i^*Q_{ij}^{(2)}A_{r_0}R_j^{(1)}x_i}_{qf_2},
\end{align}
where we used the relationship between $Q_{ij}$ and $Q_j$ evaluated at $z_2$, i.e.
$$Q_j^{(2)} - Q_{ij}^{(2)} =-\beta_{ij}^{(2)}Q_{ij}^{(2)}x_ix_i^*Q_{ij}^{(2)}.$$
Note that the inner summands in $T_5$ are products of two quadratic forms (i.e. $qf_1$ and $qf_2$). We will now expand the product of quadratic forms in $T_5$ using the equality $$XY = (X-\mathbb{E}X)(Y-\mathbb{E}Y) + (X-\mathbb{E}X)\mathbb{E}Y + X(Y - \mathbb{E}Y) +\mathbb{E}X\mathbb{E}Y$$ to get the following expression:
\begin{align*}
    T_5  &= \sum_{i<j} \psi_i^{(1)}\beta_{ij}^{(2)} \bigg(x_i^*\mathbb{E}_jQ_{ij}^{(1)}A_{s_0}Q_{ij}^{(2)}x_i\bigg)\bigg(x_i^*Q_{ij}^{(2)}A_{r_0}R_j^{(1)}x_i\bigg)\\
    &= \sum_{i<j} \psi_i^{(1)}\beta_{ij}^{(2)} \bigg(x_i^*\mathbb{E}_jQ_{ij}^{(1)}A_{s_0}Q_{ij}^{(2)}x_i - n^{-1}\operatorname{tr}\textbf{A}_i\mathbb{E}_jQ_{ij}^{(1)}A_{s_0}Q_{ij}^{(2)}\bigg)\\
    &\qquad \qquad \qquad \qquad\times \bigg(x_i^*Q_{ij}^{(2)}A_{r_0}R_j^{(1)}x_i - n^{-1}\operatorname{tr}\textbf{A}_iQ_{ij}^{(2)}A_{r_0}R_j^{(1)}\bigg)\\
    & + \sum_{i<j} \psi_i^{(1)}\beta_{ij}^{(2)} \bigg(x_i^*\mathbb{E}_jQ_{ij}^{(1)}A_{s_0}Q_{ij}^{(2)}x_i - n^{-1}\operatorname{tr}\textbf{A}_i\mathbb{E}_jQ_{ij}^{(1)}A_{s_0}Q_{ij}^{(2)}\bigg)\bigg(\frac{1}{n} \operatorname{tr}\textbf{A}_iQ_{ij}^{(2)}A_{r_0}R_j^{(1)}\bigg)\\
    & + \sum_{i<j} \psi_i^{(1)}\beta_{ij}^{(2)} \bigg(\frac{1}{n}\operatorname{tr}A_r\mathbb{E}_jQ_{ij}^{(1)}A_{s_0}Q_{ij}^{(2)}\bigg)
   \bigg(x_i^*Q_{ij}^{(2)}A_{r_0}R_j^{(1)}x_i - n^{-1}\operatorname{tr}A_rQ_{ij}^{(2)}A_{r_0}R_j^{(1)}\bigg)\\
   & + \frac{1}{n^2}\sum_{i<j}\psi_i^{(1)}\beta_{ij}^{(2)}\bigg(\operatorname{tr}\textbf{A}_i\mathbb{E}_jQ_{ij}^{(1)}A_{s_0}Q_{ij}^{(2)}\bigg)  \bigg(\operatorname{tr}\textbf{A}_iQ_{ij}^{(2)}A_{r_0}R_j^{(1)}\bigg)\\
    &:= \alpha_1(z_1,z_2)+\alpha_2(z_1,z_2)+\alpha_3(z_1,z_2)+\alpha_4(z_1,z_2).
\end{align*}
Applying the concentration results of Lemma \ref{lemma:bounds_tightness_section} gives us
$$\mathbb{E}|\alpha_1(z_1,z_2)+\alpha_2(z_1,z_2)+\alpha_3(z_1,z_2)| = O(n^\frac{1}{2}).$$
Therefore, (\ref{tracial_expansion_covariance}) simplifies to
\begin{align}\label{object_interest}
    \operatorname{tr}(A_{r_0}\mathbb{E}_jC_j^{(1)}A_{s_0}\mathbb{E}_jQ_j^{(2)}) = \mathbb{E}_j\alpha_4(z_1,z_2) + a(z_1, z_2).
\end{align}
In Lemma \ref{lemma:betaij_minus_psii}, we established that 
$$\mathbb{E}|\beta_{ij} - \psi_i| = O(n^{-\frac{1}{2}}).$$
We substitute $Q_{ij}$ with $Q_j$ to obtain
\begin{align*}
    \alpha_4(z_1, z_2) = \frac{1}{n^2}\sum_{i<j}\psi_i^{(1)}\beta_{ij}^{(2)}\bigg(\operatorname{tr}\textbf{A}_i\mathbb{E}_jQ_j^{(1)}A_{s_0}Q_{ij}^{(2)}\bigg)  \bigg(\operatorname{tr}\textbf{A}_iQ_j^{(2)}A_{r_0}R_j^{(1)}\bigg) + a(z_1, z_2).
\end{align*}
Substituting $Q_j(z_2)$ with $R_j(z_2) + C_j(z_2) + D_j(z_2) + E_j(z_2)$ and using the tracial bounds of $C_j,D_j,E_j$ (Lemma \ref{lemma:tracial_bounds}) and the non-randomness of $R_j$, we get 
\begin{align*}
    \alpha_4(z_1, z_2) :=&\frac{1}{n^2}\sum_{i<j}\psi_i^{(1)}\psi_{i}^{(2)}\bigg(\operatorname{tr}\textbf{A}_i\mathbb{E}_jQ_j^{(1)}A_{s_0}Q_{ij}^{(2)}\bigg)  \bigg(\operatorname{tr}\textbf{A}_iR_j^{(2)}A_{r_0}R_j^{(1)}\bigg) + a(z_1, z_2)\\
    =& \frac{1}{n^2}\sum_{r,s=1}^K\sum_{i < j}\psi_i^{(1)}\psi_i^{(2)} \theta_{ri}\theta_{si}\bigg(\operatorname{tr}A_r\mathbb{E}_jQ_j^{(1)}A_{s_0}Q_{ij}^{(2)}\bigg)  \bigg(\operatorname{tr}A_sR_j^{(2)}A_{r_0}R_j^{(1)}\bigg) + a(z_1, z_2),
\end{align*}
where we used $\textbf{A}_i = \sum_{r=1}^K\theta_{ri}A_r$ and the linearity of the trace operator. Substituting the above into (\ref{object_interest}), and then invoking (\ref{Vjrs_simplification}),   we have
\begin{align}\label{equation_Vrs}
    V_j[r_0,s_0] &= \operatorname{tr}(A_{r_0}R_j^{(1)}A_{s_0}R_j^{(2)}) + \mathbb{E}_j\alpha_4(z_1, z_2) + a(z_1,z_2)\\ \notag
    &=\operatorname{tr}(A_{r_0}R_j^{(1)}A_{s_0}R_j^{(2)}) + a(z_1, z_2) \\ \notag&+ \sum_{r,s=1}^K \bigg(\frac{1}{n}\sum_{i<j}\theta_{ri}\theta_{si}\psi_i^{(1)}\psi_i^{(2)}\bigg)\bigg(\operatorname{tr}A_r\mathbb{E}_jQ_j^{(1)}A_{s_0}\mathbb{E}_jQ_j^{(2)}\bigg)\bigg(\frac{1}{n}\operatorname{tr}A_sR_j^{(2)}A_{r_0}R_j^{(1)}\bigg).
\end{align}
We define the following symmetric matrix:
\begin{align}\label{defining_Ujrs}
    U_j(z_1, z_2)[r,s] = \frac{1}{p}\operatorname{tr}A_rR_j^{(1)}A_sR_j^{(2)}, \quad r,s \in [K].
\end{align}
Then, (\ref{equation_Vrs}) can be written as
\begin{align}\label{equation_Vrs1}
    V_j[r_0,s_0] &= pU_j[r_0,s_0]+c_n\sum_{r,s=1}^K \bigg(\frac{1}{n}\sum_{i<j}\theta_{ri}\theta_{si}\psi_i^{(1)}\psi_i^{(2)}\bigg)V_j[r,s_0]U_j[r_0,r] + a(z_1,z_2).
\end{align}
We now focus on the quantity below 
$$\frac{1}{n}\sum_{i<j}\theta_{ri}\theta_{si} \psi_i^{(1)}\psi_i^{(2)},$$
by developing a limit as a function of $z_1, z_2$ to bypass the dependence on $j$ altogether. By Lemma \ref{lemma:bounds_bias_section}, we have
$$\sup_{j\in[n],z\in\mathcal{C}_n}\absmod{\psi_j(z) - b_j^0(z)} = o(1).$$
Using this we define,
\begin{align}\label{defining_Wjr_Wr}
    W_j^0(z_1, z_2)[r,s] :=  \frac{1}{n}\sum_{i<j}a_i^{(r,s)}, \quad \text{ and }  \quad W^0(z_1, z_2)[r,s] := \frac{1}{n}\sum_{i=1}^na_i^{(r,s)},
\end{align}
where,
\begin{align*}
    a_i^{(r,s)} &:= \theta_{ri}\theta_{si}b_j^0(z_1)b_j^0(z_2), \quad i \in [n].
\end{align*}
Due to \textbf{T3} of Theorem \ref{mainTheorem_LSD}, we can write
$$W^0[r,s] = \int \frac{\theta_r\theta_sdG_n(\boldsymbol{\theta})}{(1 + c_n\boldsymbol{\theta}^\top \textbf{h}_n^0(z_1))(1 + c_n \boldsymbol{\theta}^\top\textbf{h}_n^0(z_2))}.$$
Now we will show that 
\begin{align}\label{Wjr_simplification}
\absmod{\bigg(\frac{1}{n}\sum_{i<j}a_i^{(r,s)}\bigg)U_j[s,r_0] - \frac{j-1}{n}W^0[r,s]U_j[s,r_0]} = o(1).
\end{align}
This will allow us to plug-in the second term i.e. $n^{-1}(j-1)W^0[r,s]U_j[s,r_0]$ into (\ref{equation_Vrs}). Following \cite{BaiLiPan2019}, we fix $\epsilon \in (0, 1/100)$. We will prove the claim separately based on the value of $j \in [n]$.

\textbf{Case 1:} Suppose $j \leq n^{1-\epsilon}$. Then, we observe that $|a_i^{(r,s)} - W^0[r,s]|$ is uniformly bounded and $|U_j[s,r_0]| = O(1)$. By triangle inequality, we have
\begin{align*}
&\absmod{\bigg(\frac{1}{n}\sum_{i<j}a_i^{(r,s)}\bigg)U_j[s,r_0] - \frac{j-1}{n}W^0[r,s]U_j[s,r_0]}\\
\leq&\, 
    \frac{1}{n}\sum_{i<j}\absmod{(a_i^{(r,s)} - 
    W^0[r,s])U_j[s,r_0]}\\
    \leq& \,\,\frac{1}{n}\sum_{i<j}\absmod{a_i^{(r,s)} - W^0[r,s]} \leq \frac{C(z_1,z_2)n^{1-\epsilon}}{n} = O\bigg(\frac{1}{n^\epsilon}\bigg).
\end{align*}
\textbf{Case 2:} Suppose $j > n^{1-\epsilon}$. Then, we have $\frac{1}{j-1}\sum_{i<j}a_i^{(r,s)} = W^0[r,s] + o(1)$. Using these gives us
\begin{align*}
\absmod{\bigg(\frac{1}{n}\sum_{i<j}a_i^{(r,s)}\bigg)U_j[s,r_0] - \frac{j-1}{n}W^0[r,s]U_j[s,r_0]} 
\leq& \absmod{\frac{1}{n}\sum_{i<j} a_i^{(r,s)} - \frac{j-1}{n}W^0[r,s]}\\
=& \frac{j-1}{n}\absmod{\frac{1}{j-1}\sum_{i<j} a_i^{(r,s)} - W^0[r,s]}= o(1).
\end{align*}

Using (\ref{Wjr_simplification}) in (\ref{equation_Vrs1}), we therefore have
\begin{align}\label{equation_Vrs2}
    &V_j[r_0,s_0] = pU_j[r_0,s_0] + c_n\sum_{r=1}^K \frac{j-1}{n}W^0[r,s]U_j[s,r_0]V_j[r,s_0] + a(z_1,z_2)\\\notag
    \iff& V_j = pU_j + c_n\frac{j-1}{n}U_jW^0V_j + a(z_1,z_2)\\\notag
    \iff& \bigg(I - c_n\frac{j-1}{n}U_jW^0\bigg)V_j = pU_j + A(z_1, z_2).
\end{align}
We now seek to replace $U_j$ in the above expressions with another term that is independent of $j$. Recall the definition of $R_j(z)$ from (\ref{essential_quantities}) which is a component of $U_j$. We now define
\begin{align}\label{defining_U0_rs}
    U^0(z_1, z_2)[r,s] := \frac{1}{p}\operatorname{tr}A_rR^0(z_1)A_sR^0(z_2), \quad r,s \in [K],
\end{align}
where
$$R^0(z) := \bigg(\frac{1}{n}\sum_{i =1}^nb_i^0(z)\textbf{A}_i - zI_p\bigg)^{-1},$$
Note that $U^0$ is a symmetric matrix just like $U_j$ but it is free of $j$. From Lemma \ref{lemma:Uj_U0_U_element_wise_difference}, we have 
$$\|U_j - U^0\| = o(1).$$
In light of this, we now have
\begin{align}\label{expression_Vj_1}
    \frac{1}{p}\bigg(I - c_n\frac{j-1}{n}U^0W^0\bigg)V_j = U^0 + \frac{1}{p}A(z_1, z_2).
\end{align}

By Lemma \ref{lemma:spectral_radius_le_1}, we therefore conclude that
\begin{align}\label{expression_Vj}
    \frac{1}{n}V_j = c_n\bigg(I - c_n\frac{j-1}{n}U^0W^0\bigg)^{-1}U^0 + \frac{1}{n}A(z_1,z_2).
\end{align}
Note that, by Lemma \ref{lemma:norm_control}, we have $$\mathbb{E} \bigl\|\frac{1}{p}A(z_1, z_2)\bigr\| = O(n^{-\frac{1}{2}}).$$
\footnote{$A = o_{\|.\|}(1)$ indicates that $\|A\| = o(1)$. The notation $O_{\|.\|}$ is also similarly defined.} 
Using all the simplifications done in this Section, (\ref{revised_goal1}) reduces to 
\begin{align}
    \bigl|\operatorname{Cov}(M_{n1}(z_1), M_{n2}(z_2)) - \operatorname{tr}\bigg(\frac{1}{n^2}\sum_{j=1}^nd_j\Theta_jV_j\bigg)\bigr| \overset{p}{\longrightarrow} 0,
    \end{align}
where, 
\begin{align}
    d_{j} =b_j^0(z_1)b_j^0(z_2), \quad j \in [n].
\end{align}
We now investigate the limit of the quantity
$$\operatorname{tr}\bigg(\frac{1}{n^2}\sum_{j=1}^n d_j\Theta_jV_j\bigg).$$
Let $M_j = \frac{1}{j}\sum_{i=1}^j d_i\Theta_i$ denote the partial averages. Using Abel's Lemma (\ref{Lemma:Abel}), we have
\begin{align*}
    \frac{1}{n^2}\sum_{j=1}^n d_j\Theta_jV_j &= \sum_{j=1}^n\bigg(\frac{1}{n}d_j\Theta_j\bigg) \bigg(\frac{1}{n}V_j\bigg)\\
    &=\bigg(\frac{1}{n}\sum_{j=1}^nd_j\Theta_j\bigg)\frac{1}{n}V_n - \sum_{j=1}^{n-1}\bigg(\frac{1}{n}\sum_{i=1}^j d_i\Theta_i\bigg)\frac{V_{j+1}-V_j}{n}\\
    &=\frac{1}{n}M_nV_n - \sum_{j=1}^{n-1}\frac{j}{n}\frac{M_j(V_{j+1}-V_j)}{n}.
\end{align*}
Note that $M_n = W^0 + o_{\|.\|}(1)$, $n^{-1}V_n = c_n(I-c_nU^0W^0)^{-1}U^0 + o_{\|.\|}(n^{-\frac{1}{2}})$ and by the resolvent identity (\ref{Resolvent_Identity}), we have
$$\frac{V_{j+1} - V_j}{n} = \frac{c_n^2}{n}T_{j+1}U^0W^0T_jU^0 + o_{\|.\|}(1),$$
where $$T_{j+1} := \bigg(I - c_n\frac{j}{n}U^0W^0\bigg)^{-1}.$$
This leads to the following simplification:
\begin{align}\label{trace_expansion}
   \operatorname{tr}\bigg(\frac{1}{n^2}\sum_{j=1}^n d_j\Theta_jV_j\bigg) =\,\operatorname{tr}((I - c_nU^0W^0)^{-1}c_nU^0W^0) - \frac{c_n^2}{n}\operatorname{tr}\bigg(\sum_{j=1}^{n-1}\frac{j}{n} M_jT_{j+1}U^0W^0T_jU^0\bigg) + o(1).
\end{align}
We split the summation in the second term of (\ref{trace_expansion}) as shown below. For $\epsilon \in (0, 1/100)$, we have 
\begin{align}\label{summation_part_1}
    \frac{1}{n}\sum_{j \leq n^{1-\epsilon}}\frac{j}{n}\absmod{\operatorname{tr}\bigg((M_j - W^0)T_{j+1}U^0W^0T_jU^0\bigg)} \leq \frac{C}{n^\epsilon} = o(1).
\end{align}
This follows since each term is uniformly bounded as explained below.
\begin{enumerate}
    \item $U^0,W^0$ are $K\times K$ matrices satisfying $\|W^0\|, \|U^0\| < \infty$.
    \item For each $j \in [n]$, $\|T_j\|$ is bounded by Remark \ref{radius_uniformly_away_from_1_on_contour}. As a result, by the resolvent identity (\ref{Resolvent_Identity}), we also have 
\begin{align*}
    \|T_{j+1} - T_j \| \leq \frac{c_n}{n}\|T_{j+1}\| \|U^0W^0\| \|T_j\| \to 0.
\end{align*}
    \item $\{\|M_j - W^0\|\}_{j=1}^n$ is a bounded sequence.
\end{enumerate}
For the other part, using $|\operatorname{tr}(XY)| \leq \|X\|\,\|Y\|_*$, where $\|\cdot\|_*$ denotes the nuclear norm, we have
\begin{align}\label{summation_part_2}
    \frac{1}{n}\sum_{j > n^{1-\epsilon}}\frac{j}{n}\absmod{\operatorname{tr}\bigg((M_j - W^0)T_{j+1}U^0W^0T_j\bigg)}
    \leq& \, \frac{C}{n}\sum_{j \geq n^{1-\epsilon}} \|M_j-W^0\| = o(1).
\end{align}
Combining (\ref{summation_part_1}) and (\ref{summation_part_2}), we simplify (\ref{trace_expansion}) using Lemma \ref{lemma:tracial_integration} as follows:
\begin{align}\label{trace_expansion_1}
    &\operatorname{tr}\bigg(\frac{1}{n^2} \sum_{j=1}^n d_j\Theta_jV_j\bigg) \\\notag
    =&\,\operatorname{tr}((I - c_nU^0W^0)^{-1}c_nU^0W^0) - \operatorname{tr}\bigg(\frac{c_n^2}{n}\sum_{j=1}^{n-1}\frac{j}{n} W^0T_{j+1}U^0W^0T_jU^0\bigg) + o(1)\\ \notag
    =&\, \operatorname{tr}\bigg((I-c_nU^0W^0)^{-1}c_nU^0W^0\bigg) - \int_0^1 t \operatorname{tr}\bigg((I - c_ntU^0W^0)^{-1}c_nU^0W^0\bigg)^2dt + o(1)\\\notag
    =& -\operatorname{tr}\log(I - c_nU^0W^0) + o(1).
\end{align}

\begin{figure}[!htbp]
\centering
 \begin{subfigure}{0.8\textwidth}
 \centering
\includegraphics[width=\textwidth,height=6cm]{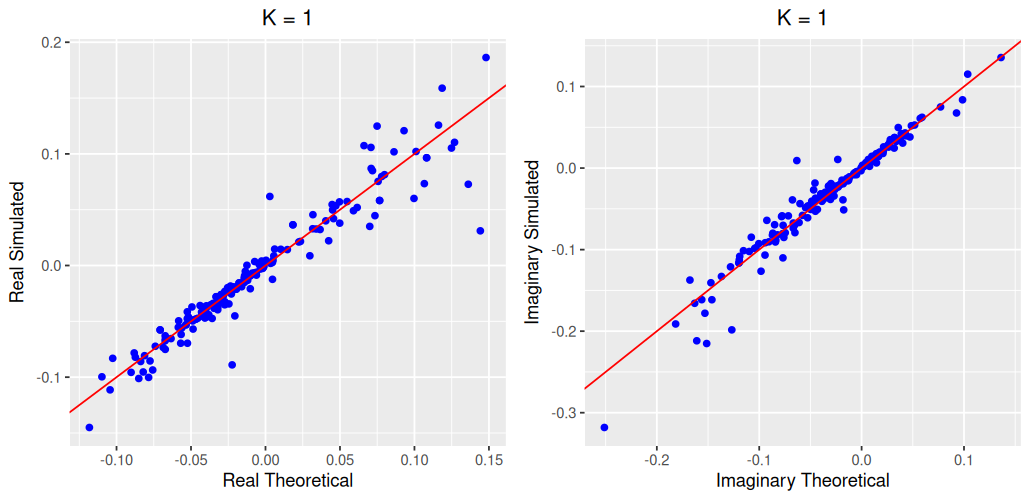}
     \caption{K=1:c=1}
     \label{fig1}
 \end{subfigure}
 \medskip
 \begin{subfigure}{0.8\textwidth}
\includegraphics[width=\textwidth,height=6cm]{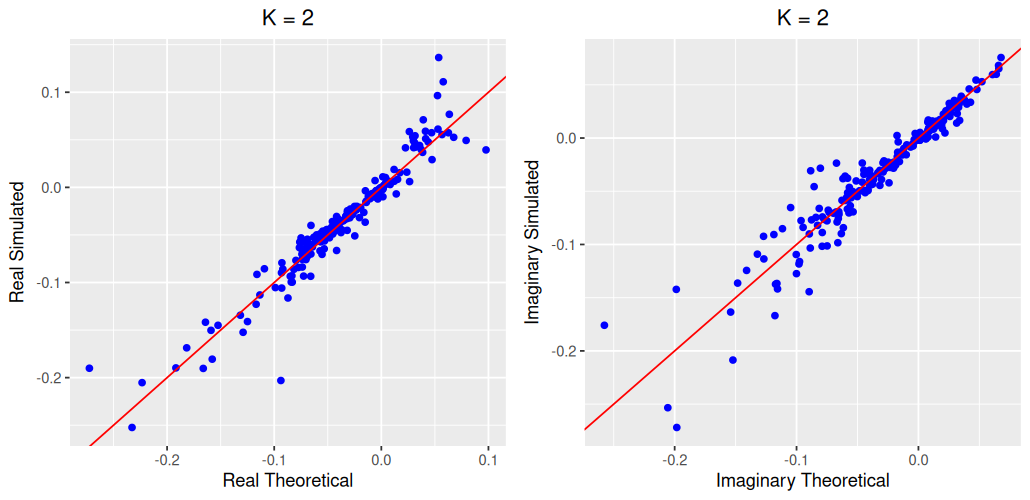}
     \caption{K=2:c=1}
     \label{fig2}
 \end{subfigure}
\medskip
 \begin{subfigure}{0.8\textwidth}
\includegraphics[width=\textwidth,height=6cm]{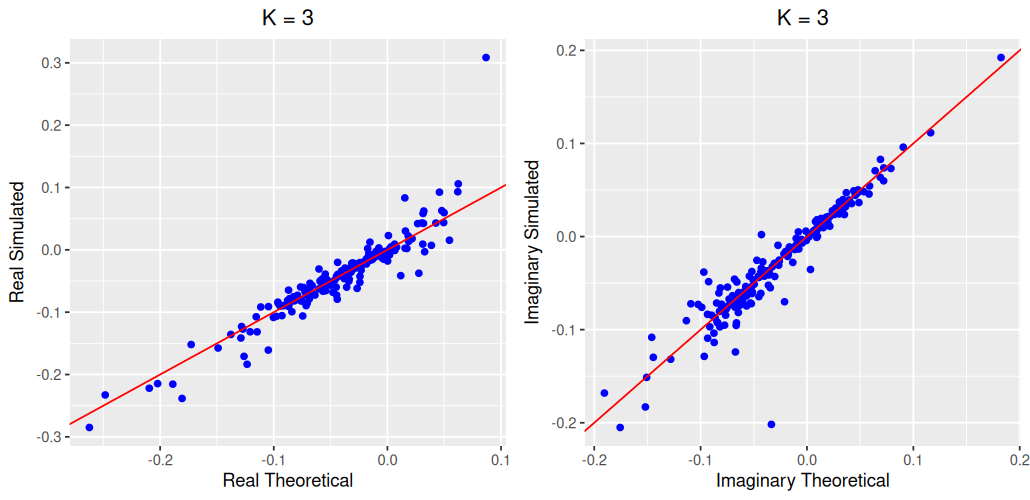}
     \caption{K=3:c=1}
     \label{fig3}
 \end{subfigure}
 \caption{Scatterplot of Simulated vs. Theoretical values (real and imaginary parts in the left and right panels respectively) of $Cov(M_n(z_1), M_n(z_2))$ calculated at 200 pairs of $z_1, z_2$ for $K=1,2,3$ and $c=1$ with standard Gaussian innovations.}
 \label{fig1_to_6}
\end{figure}

\begin{figure}
\centering
 \begin{subfigure}{0.75\textwidth}
\includegraphics[width=\textwidth,height=6cm]{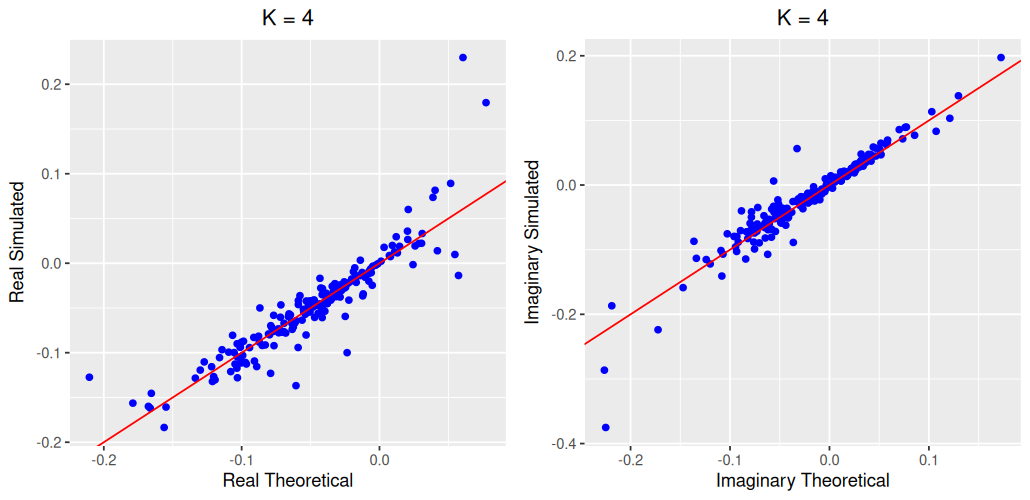}
     \caption{K=4:c=1}
     \label{fig4}
 \end{subfigure}
 \medskip
   \begin{subfigure}{0.75\textwidth}     \includegraphics[width=\textwidth,height=6cm]{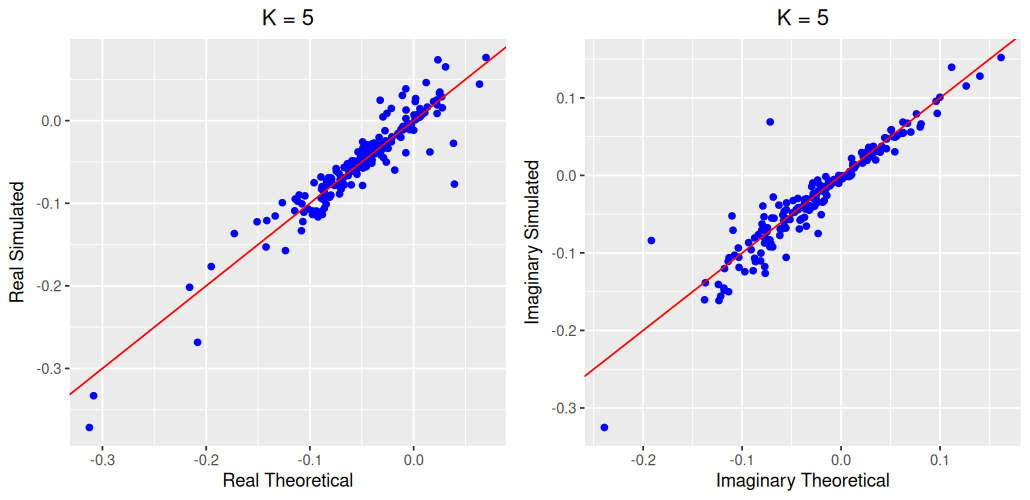}
     \caption{K=5:c=1}
     \label{fig5}
 \end{subfigure}
 \medskip
 \begin{subfigure}{0.75\textwidth}
\includegraphics[width=\textwidth,height=6cm]{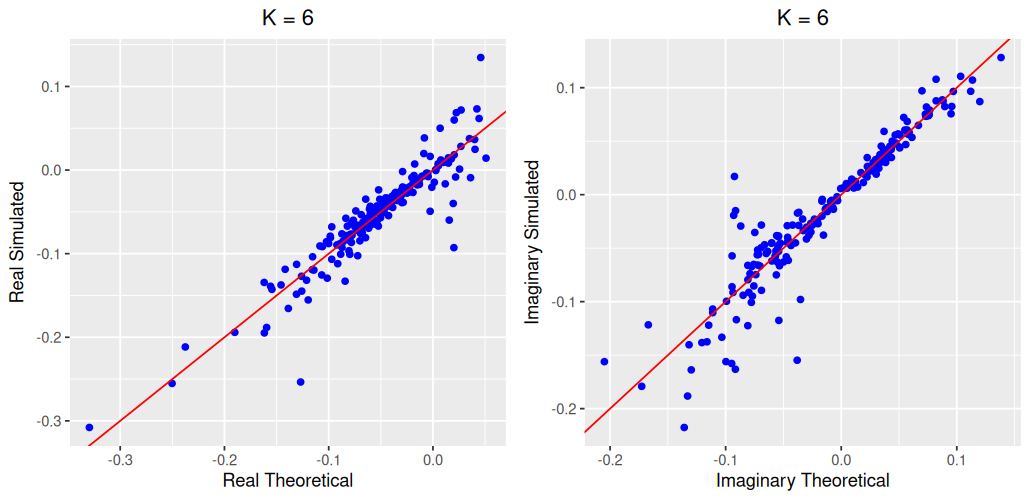}
     \caption{K=6:c=1}
     \label{fig6}
 \end{subfigure}
 \caption{Scatterplot of Simulated vs. Theoretical values (real and imaginary parts in the left and right panels respectively) of $Cov(M_n(z_1), M_n(z_2))$ calculated at 200 pairs of $z_1, z_2$ for $K=4,5,6$ and $c=1$ with standard Gaussian innovations.}
  \label{fig7_to_12}
\end{figure}

Recall the definitions of $U_{2,2}$ (\ref{defining_population_U_matrices}) and $W_2$ (\ref{defining_population_W_matrices}). From Lemma \ref{lemma:hn_Ehn_hn0_uniform_results} and DCT, we have
$$\|U^0(z_1, z_2) - U_{2,2}(z_1, z_2)\|_{op} = o(1), \quad \|W^0(z_1, z_2) - W_2(z_1,z_2)\|_{op} = o(1).$$  
By Lemma \ref{lemma:transfer_of_derivative}, we finally infer that
\begin{align*}
    \operatorname{Cov}(M(z_1), M(z_2)) &= \underset{n \rightarrow \infty}{\lim} \operatorname{Cov}(M_{n1}(z_1), M_{n1}(z_2))\\
&= -\frac{\partial^2}{\partial z_2\partial z_1}\operatorname{tr}\bigg(\log\bigl(I - cU_{2,2}(z_1,z_2)W_2(z_1,z_2)\bigr)\bigg).
\end{align*}
This concludes the proof of the covariance part (\ref{covariance_function_real}) of Lemma \ref{lemma:GaussianProcess}. A few numerical simulations (see Figures \ref{fig1_to_6}, \ref{fig7_to_12}) have been added for validation.
\newpage
\subsection{Numerical Validation}

\section{Deterministic/ Bias Part}\label{sec:BiasTerm}
In this section, we determine the limit of the deterministic component, i.e., $M_{n2}(z)$ derived in (\ref{decomposing_Mn}). The basic goal of this section is to demonstrate a decomposition of the form
$$\mathbb{E}s_n(z) = s_n^0(z) + \frac{1}{p}m_2(z) + o\bigg(\frac{1}{n}\bigg),$$
which leads to an expression for $m_2(z)$ involved in the mean function (\ref{mean_Xf}). Note that unlike the covariance function which was sufficient to establish only for the horizontal arms, we derive the limiting bias expression for the contour $\mathcal{C}_n$. 

Recall the definitions of $\psi_j(z)$, $b_j^0(z)$ $\textbf{A}_j$. Then, similar to $R(z)$ which serves as a deterministic equivalent of $Q(z)$, we define the finite sample analog of $Q(z)$ as follows:
\begin{align}\label{defining_R0z}
    R^0(z) = \bigg(\frac{1}{n}\sum_{j=1}^n b_j^0(z)\textbf{A}_j - zI_p\bigg)^{-1}.
\end{align}

Observe that using properties (\ref{finite_sample_properties}), we can also express $R^0(z)$ as follows:
\begin{align}
    R^0(z) = \bigg(\sum_{r=1}^K \bigg(\frac{1}{n}\sum_{j=1}^nb_j^0\theta_{rj}\bigg)A_r - zI_p\bigg)^{-1} = \bigg(-zI_p - \sum_{r=1}^K zg_{rn}^0(z)A_r\bigg)^{-1}.
\end{align}

We start with the following decomposition of $M_{n2}(z)$.
\begin{align}\label{decomposing_Mn2}
  M_{n2}(z) &= p(\mathbb{E}s_n(z) - s_n^0(z))\\ \notag
  &= \operatorname{tr}(\mathbb{E}Q(z) - R(z)) + \operatorname{tr}(R(z) - R^0(z))\\ \notag
  &:= d_{0n}(z) + e_{0n}(z).
\end{align}
Analogously, we also define the following: 
\begin{align}
    d_{rn}(z) :=\operatorname{tr}(\mathbb{E}Q(z)-R(z))A_r, \quad e_{rn}(z) :=  \operatorname{tr}(R(z)-R^0(z))A_r, \quad r \in [K].
\end{align}
Before moving on to the analysis of $d_{0n}, e_{0n}, \{d_{rn}\}_{r=1}^K, \{e_{rn}\}_{r=1}^K$, an important step is to show that $\|R(z)\|$ (see \ref{essential_quantities1}) is uniformly bounded for $z \in \mathcal{C}_n$. (proof in Section \ref{proof_Rz_uniform_bound})
\begin{lemma}\label{lemma:norm_Rz_uniform_bound}
    In addition to Lemma \ref{lemma:Bounds_Important_Quantities}, the following holds:
    $\underset{n;z \in \mathcal{C}_n}{\sup}\|R(z)\| < \infty$.
\end{lemma}

\subsection{Analysis of \texorpdfstring{$d_{0n}(z)$}{} and \texorpdfstring{$d_{rn}(z)$}{}}
Although $M_{n2}(z) = d_{0n}(z) + e_{0n}(z)$, the analysis of $e_{0n}(z)$ involves $\textbf{d}_n(z) = \{d_{rn}(z)\}_{r=1}^K$ as we shall see in the next subsection. This comes at little extra cost since the expressions for $d_{rn}(z)$ can be developed by arguments similar to the ones to be used for $d_{0n}(z)$.
We start by expanding $d_{0n}(z)$ using the resolvent identity (\ref{Resolvent_Identity}) as follows:
\begin{align*}
    \mathbb{E}Q(z) - R(z)
    =& -R(z)\mathbb{E}\bigg(\sum_{j=1}^nx_jx_j^*Q(z) - \frac{1}{n}\sum_{j=1}^n\psi_j(z)\textbf{A}_jQ(z)\bigg)\\
    =& -R(z)\mathbb{E}\bigg(\sum_{j=1}^n\beta_j(z)x_jx_j^*Q_j(z) - \frac{1}{n}\sum_{j=1}^n\psi_j(z)\textbf{A}_jQ(z)\bigg).
\end{align*}
Masking the dependence of $\beta_j(\cdot), \psi_j(\cdot), b_j(\cdot),\epsilon_j(\cdot), \rho_j(\cdot), R(\cdot), Q(\cdot), Q_j(\cdot)$ on $z$ for simplicity, we have
\begin{align*}
    d_{0n}(z)=& \operatorname{tr}(\mathbb{E}Q(z) - R(z))\\
    =&-\sum_{j=1}^n \mathbb{E}\beta_jx_j^*Q_jRx_j + \frac{1}{n}\sum_{j=1}^n\psi_j\operatorname{tr}(\textbf{A}_jQR)\\
    =& -\sum_{j=1}^n \mathbb{E}\beta_j\bigl(x_j^*Q_jRx_j - n^{-1}\mathbb{E}\operatorname{tr}(\textbf{A}_jQ_jR)\bigr)\\
    &-\frac{1}{n}\sum_{j=1}^n\mathbb{E}\beta_j\mathbb{E}\operatorname{tr}(\textbf{A}_j(Q_j-Q)R)\bigr)\\
    &-\frac{1}{n}\sum_{j=1}^n \mathbb{E}(\beta_j - \psi_j)\mathbb{E}\operatorname{tr}(\textbf{A}_jQR)\\
    &:= \mathcal{I}_1+\mathcal{I}_2+\mathcal{I}_3,
\end{align*}
where $\rho_j(z)$ is defined in (\ref{defining_rhoj}).Using (\ref{beta_b_relationship}), $\mathcal{I}_1$ can be written as
\begin{align*}
    \mathcal{I}_1 &= \sum_{j=1}^n b_j^2\mathbb{E}\rho_j(x_j^*Q_jRx_j -  n^{-1}\mathbb{E}\operatorname{tr}(\textbf{A}_jQ_jR)) \\
   & - \sum_{j=1}^n b_j^2 \mathbb{E}\beta_j\rho_j^2\bigl(x_j^*Q_jRx_j - n^{-1}\operatorname{tr}(\textbf{A}_jQ_jR)\bigr)\\
    &- \frac{1}{n}\sum_{j=1}^n b_j^2 \mathbb{E}\beta_j\rho_j^2\bigl(\operatorname{tr}(\textbf{A}_jQ_jR) - \mathbb{E}\operatorname{tr}(\textbf{A}_jQ_jR))\bigr)\\
    &= \mathcal{I}_{11} + \mathcal{I}_{12} + \mathcal{I}_{13}.
\end{align*}
By Lemmas \ref{lemma:bounds_tightness_section} and \ref{lemma:bounds_bias_section}, we have
$$\mathcal{I}_{12} \leq C\sum_{j=1}^n\mathbb{E}^\frac{1}{4}|\beta_j|^4\mathbb{E}^\frac{1}{4}|\rho_j|^8\mathbb{E}^\frac{1}{2}\bigg(x_j^*Q_jRx_j - n^{-1}\operatorname{tr}(\textbf{A}_jQ_jR)\bigg)^2 \leq \frac{C}{\sqrt{n}} \rightarrow 0,$$
and 
$$\mathcal{I}_{13} \leq \frac{C}{n}\sum_{j=1}^n\mathbb{E}^\frac{1}{4}|\beta_j|^4\mathbb{E}^\frac{1}{4}|\rho_j|^8\mathbb{E}^\frac{1}{2}\bigl|\operatorname{tr}(\textbf{A}_jQ_jR) - \mathbb{E}\operatorname{tr}(\textbf{A}_jQ_jR)\bigr|^2 \leq \frac{C}{n} \rightarrow 0.$$

We again decompose $\mathcal{I}_{11}$ as follows:
\begin{align*}
    \mathcal{I}_{11} &= \sum_{j=1}^n b_j^2\mathbb{E}\epsilon_j (x_j^*Q_jRx_j - n^{-1}\operatorname{tr}(\textbf{A}_jQ_jR)) + \frac{1}{n^2}\sum_{j=1}^n b_j^2\operatorname{Cov}(\operatorname{tr}\textbf{A}_jQ_j, \operatorname{tr}\textbf{A}_jQ_jR)\\
    &= \mathcal{I}_{111} + \mathcal{I}_{112}.
\end{align*}

We again have $|\mathcal{I}_{112}| \leq C/n$. Therefore, we have
\begin{align}\label{script_I1_final}
    \mathcal{I}_1 = \sum_{j=1}^n b_j^2\mathbb{E}\epsilon_j (x_j^*Q_jRx_j - n^{-1}\operatorname{tr}(\textbf{A}_jQ_jR)) + o(1).
\end{align}
Now we decompose $\mathcal{I}_2$. Using (\ref{Leave_one_out_property}), we have
\begin{align*}
  \mathcal{I}_2 =&-\frac{1}{n}\sum_{j=1}^n\mathbb{E}\beta_j\mathbb{E}\operatorname{tr}(\textbf{A}_j(Q_j-Q)R)\bigr)\\
  &= -\frac{1}{n}\sum_{j=1}^n b_j^2 \mathbb{E}x_j^*Q_jR\textbf{A}_jQ_jx_j + \frac{1}{n}\sum_{j=1}^n b_j^2\mathbb{E}\beta_j\rho_jx_j^*Q_jR\textbf{A}_jQ_jx_j+ \frac{1}{n}\sum_{j=1}^n b_j\mathbb{E} \beta_j\rho_j\mathbb{E}x_j^*Q_jR\textbf{A}_jQ_jx_j\\
  &:= \mathcal{I}_{21}+\mathcal{I}_{22}+\mathcal{I}_{23}.
\end{align*}
From Lemma \ref{lemma:bounds_tightness_section}, we have $|\mathcal{I}_{22}+\mathcal{I}_{23}| \leq C/\sqrt{n}$. Therefore,
\begin{align}\label{script_I2_final}
    \mathcal{I}_2 = -\frac{1}{n^2}\sum_{j=1}^n b_j^2 \mathbb{E}\operatorname{tr}R\textbf{A}_jQ_j\textbf{A}_jQ_j + o(1).
\end{align}

Recall the expression for $(\mathbb{E}\beta_j -\psi_j)$ from (\ref{exp_beta_j_minus_psi_j}). This leads to the following expression for $\mathcal{I}_3:$
\begin{align}\label{script_I3_final}
    \mathcal{I}_3 = -\frac{1}{n}\sum_{j=1}^n b_j^3\mathbb{E}\epsilon_j^2\mathbb{E}\operatorname{tr}R\textbf{A}_jQ + \frac{1}{n^3}\sum_{j=1}^n b_j^2\psi_j\bigl(\mathbb{E}\operatorname{tr}Q_j\textbf{A}_jQ_j\textbf{A}_j\bigr)\bigl(\mathbb{E}\operatorname{tr}R\textbf{A}_jQ\bigr) + o(1).
\end{align}
Combining (\ref{script_I1_final}), (\ref{script_I2_final}) and (\ref{script_I3_final}), we have
\begin{align}
    d_{0n}(z) &= \sum_{j=1}^n b_j^2\mathbb{E}\epsilon_j (x_j^*Q_jRx_j - n^{-1}\operatorname{tr}(\textbf{A}_jQ_jR))\\ \notag
    &-\frac{1}{n^2}\sum_{j=1}^n b_j^2 \mathbb{E}\operatorname{tr}R\textbf{A}_jQ_j\textbf{A}_jQ_j\\ \notag
    &-\frac{1}{n}\sum_{j=1}^n b_j^3\mathbb{E}\epsilon_j^2\mathbb{E}\operatorname{tr}R\textbf{A}_jQ\\ \notag
    &+ \frac{1}{n^3}\sum_{j=1}^n b_j^2\psi_j\bigl(\mathbb{E}\operatorname{tr}Q_j\textbf{A}_jQ_j\textbf{A}_j\bigr)\bigl(\mathbb{E}\operatorname{tr}R\textbf{A}_jQ\bigr) + o(1).
\end{align}
Also by Lemma \ref{lemma:exp_betaj_psij_bj_difference}, we have  $\sup_{j,z}|b_j - \psi_j| = O(n^{-1})$. Therefore, we have
\begin{align}\label{dnz_consolidated_expression}
    d_{0n}(z) &= \sum_{j=1}^n \psi_j^2\mathbb{E}\epsilon_j (x_j^*Q_jRx_j - n^{-1}\operatorname{tr}(\textbf{A}_jQ_jR))\\ \notag
    &-\frac{1}{n^2}\sum_{j=1}^n \psi_j^2 \mathbb{E}\operatorname{tr}R\textbf{A}_jQ_j\textbf{A}_jQ_j\\ \notag
    &-\frac{1}{n}\sum_{j=1}^n \psi_j^3\mathbb{E}\epsilon_j^2\operatorname{tr}R\textbf{A}_jR\\ \notag
    &+ \frac{1}{n^3}\sum_{j=1}^n \psi_j^3\bigl(\mathbb{E}\operatorname{tr}Q_j\textbf{A}_jQ_j\textbf{A}_j\bigr)\bigl(\operatorname{tr}R\textbf{A}_jR\bigr) + o(1),
\end{align}
where we used the fact from Lemma \ref{lemma:tracial_expansion_bias} that
\begin{align}\label{average_tracial_deviation}
    \frac{1}{n}\bigl|\mathbb{E}\operatorname{tr}R\textbf{A}_jQ - \operatorname{tr}R\textbf{A}_jR\bigr| = o(1).
\end{align}
\color{black}
When the innovation matrices are real-valued, an application of Lemma \ref{lemma:covariance_quadratic_forms} gives us
$$\mathbb{E}\epsilon_j (x_j^*Q_jRx_j - n^{-1}\operatorname{tr}(\textbf{A}_jQ_jR)) = \frac{2}{n^2}\sum_{j=1}^n \psi_j^2 \mathbb{E}\operatorname{tr}R\textbf{A}_jQ_j\textbf{A}_jQ_j + o(1),$$ and 
$$\mathbb{E}\epsilon_j^2 = \frac{2}{n^2}\operatorname{tr}Q_j\textbf{A}_jQ_j\textbf{A}_j + o(1).$$
We find similar expressions in the complex case except that the factor of $2$ gets replaced by $1$.
Therefore, it turns out that in the complex case we have
\begin{align}\label{dnz_complex_case}
d_{0n}(z) = o(1)
\end{align}
and for the real case, we have
\begin{align}\label{dnz_real_case}
    d_{0n}(z)&= \frac{1}{n^2}\sum_{j=1}^n \psi_j^2 \mathbb{E}\operatorname{tr}R\textbf{A}_jQ\textbf{A}_jQ -\frac{1}{n^3}\sum_{j=1}^n \psi_j^3\bigl(\mathbb{E}\operatorname{tr}Q\textbf{A}_jQ\textbf{A}_j\bigr)\bigl(\operatorname{tr}R\textbf{A}_jR\bigr) + o(1)\\ \notag
    &=c_n\sum_{r,s=1}^K\bigg(\frac{1}{n}\sum_{j=1}^n \theta_{rj}\theta_{sj}\psi_j^2\bigg)\bigg(\frac{1}{p}\mathbb{E}\operatorname{tr}RA_rQA_sQ\bigg)\\\notag
    &-c_n^2\sum_{r,s,t=1}^K\bigg(\frac{1}{n}\sum_{j=1}^n \theta_{rj}\theta_{sj}\theta_{tj}\psi_j^3\bigg) \bigg(\frac{1}{p}\mathbb{E}\operatorname{tr}QA_rQA_s\bigg)\bigg(\frac{1}{p}\operatorname{tr}RA_tR\bigg) + o(1).
\end{align}

\color{black}
Following arguments similar to those that led to (\ref{dnz_complex_case}) and (\ref{dnz_consolidated_expression}), we have $d_{rn}(z) = o(1)$ in the complex case and for the real case, 
\begin{align}\label{drnz_real_case}
    d_{rn}(z) &= \mathbb{E}(Q(z) - R(z))A_r\\ \notag
    &= \frac{1}{n^2}\sum_{j=1}^n \psi_j^2 \mathbb{E}\operatorname{tr}R\textbf{A}_jQ\textbf{A}_jQA_r -\frac{1}{n^3}\sum_{j=1}^n \psi_j^3\bigl(\mathbb{E}\operatorname{tr}Q\textbf{A}_jQ\textbf{A}_j\bigr)\bigl(\mathbb{E}\operatorname{tr}R\textbf{A}_jRA_r\bigr) + o(1)\\ \notag
    &= c_n\sum_{s,t=1}^K \bigg(\frac{1}{n}\sum_{j=1}^n \theta_{sj}\theta_{tj}\psi_j^2\bigg)\bigg(\frac{1}{p}\mathbb{E}\operatorname{tr}RA_sQA_tQA_r\bigg)\\\notag
    &-c_n^2\sum_{s,t,u=1}^K\bigg(\frac{1}{n}\sum_{j=1}^n \theta_{sj}\theta_{tj}\theta_{uj}\psi_j^3\bigg) \bigg(\frac{1}{p}\mathbb{E}\operatorname{tr}QA_sQA_t\bigg)\bigg(\frac{1}{p}\mathbb{E}\operatorname{tr}RA_uRA_r\bigg) + o(1).
\end{align}
Before proceeding we remark that (\ref{dnz_real_case}) and (\ref{drnz_real_case}) comprise of expressions of three kinds as follows:
\begin{enumerate}
    \item[1] \textbf{Type-W:} These are averages of deflators which depend upon the $j$-th column of the data matrix. 
    \item[2] \textbf{Type-U:} These are average traces of products of deterministic approximations of resolvents ($R(z)$) and $A_r$ type matrices.
    \item[3] \textbf{Type-V:} These are average traces of products of resolvents ($Q(z)$), deterministic approximations of resolvents ($R(z)$)  and $A_r$ type matrices.
\end{enumerate}

We handle the \textbf{type-W} objects in (\ref{dnz_real_case}) and (\ref{drnz_real_case}) as follows. For $r,s,t \in [K]$, we have
\begin{align}\label{defining_W2^0_W3^0}
   &\frac{1}{n}\sum_{j=1}^n \psi_j^2\theta_{rj}\theta_{sj} = \int \frac{\theta_r\theta_sdG_n(\boldsymbol{\theta})}{(1 + c_n \boldsymbol{\theta}^\top \textbf{h}_n^0(z))^2} + o(1) := W_2^0(z)[r,s] + o(1),\\ \notag
   &\frac{1}{n}\sum_{j=1}^n\psi_j^3\theta_{rj}\theta_{sj}\theta_{tj} = \int \frac{\theta_r\theta_s\theta_tdG_n(\boldsymbol{\theta})}{(1 + c_n \boldsymbol{\theta}^\top \textbf{h}_n^0(z))^3} + o(1) := W_3^0(z)[r,s,t] + o(1),
\end{align}
where the expressions in the middle follow from Lemma \ref{lemma:hn_Ehn_hn0_uniform_results}. $W_3^0(z)$ is reshaped as a $K^2 \times K$ object for multiplicative compatibility as we shall see. For the \textbf{type-U} objects, we introduce the following notations:
\begin{align*}
    & U_{1,2}^{(n)}(z)[r] :=  \frac{1}{p}\operatorname{tr}R(z)A_rR(z)=\int \frac{\lambda_r dH_n(\boldsymbol{\lambda})}{(-z + \frac{1}{n}\sum_{j=1}^n\psi_j\boldsymbol{\theta}_j^\top \boldsymbol{\lambda})^2},\\
    &U_{2,2}^{(n)}(z)[r,s] :=  \frac{1}{p}\operatorname{tr}R(z)A_rR(z)A_s = \int \frac{\lambda_r\lambda_s dH_n(\boldsymbol{\lambda})}{(-z + \frac{1}{n}\sum_{j=1}^n\psi_j\boldsymbol{\theta}_j^\top \boldsymbol{\lambda})^2},\\
    &U_{2,3}^{(n)}(z)[r,s] :=  \frac{1}{p}\operatorname{tr}R(z)A_rR(z)A_sR(z) =\int \frac{\lambda_r\lambda_s dH_n(\boldsymbol{\lambda})}{(-z + \frac{1}{n}\sum_{j=1}^n\psi_j\boldsymbol{\theta}_j^\top \boldsymbol{\lambda})^3},\\
   &U_{3,3}^{(n)}(z)[r,s,t] :=  \frac{1}{p}\operatorname{tr}R(z)A_rR(z)A_sR(z)A_t =\int \frac{\lambda_r\lambda_s\lambda_t dH_n(\boldsymbol{\lambda})}{(-z + \frac{1}{n}\sum_{j=1}^n\psi_j\boldsymbol{\theta}_j^\top \boldsymbol{\lambda})^3}.
\end{align*}
\begin{remark}
    The integral expressions arise due to the pairwise commutativity of $\{A_r\}_{r=1}^K$.
\end{remark}

Finally, for $r,s,t \in [K]$, the \textbf{type-V} objects are denoted as follows:
\begin{align}\label{defining_the_V_objects}
    &V_{1,2}^{(n)}(z)[r] = \frac{1}{p}\mathbb{E}\operatorname{tr}R(z)A_rQ(z),\\ \notag
    &V_{2,2}^{(n)}(z)[r,s] = \frac{1}{p}\mathbb{E}\operatorname{tr}Q(z)A_rQ(z)A_s,\\ \notag
    &V_{2,3}^{(n)}(z)[r,s] = \frac{1}{p}\mathbb{E}\operatorname{tr}R(z)A_rQ(z)A_sQ(z),\\ \notag
    &V_{3,3}^{(n)}(z)[r,s,t] = \frac{1}{p}\mathbb{E}\operatorname{tr}R(z)A_rQ(z)A_sQ(z)A_t.
\end{align}

Continuing from (\ref{dnz_real_case}), we use this terminology to  express $d_{0n}(z)$ in the real case as follows:
\begin{align}\label{dnz_intermsof_V_W}
    d_{0n}(z) &= c_n\sum_{r,s=1}^K W_2^0[r,s]V_{2,3}^{(n)}[r,s]- c_n^2\sum_{r,s,t=1}^KW_3^0[r,s,t]V_{2,2}^{(n)}[r,s]U_{1,2}^{(n)}[t] + o(1)\\\notag
    &= c_n\operatorname{tr}W_2^0V_{2,3}^{(n)} - c_n^2\operatorname{vec}(V_{2,2}^{(n)})^\top W_3^0U_{1,2}^{(n)} + o(1).
\end{align}
Similarly from (\ref{drnz_real_case}), for any $r \in [K]$, $d_{rn}(z)$  in the real case can be expressed as follows:
\begin{align}\label{drnz_intermsof_V_W}
    d_{rn}(z) &= c_n\sum_{s,t=1}^K W_2^0[s,t]V_{3,3}^{(n)}[r,s,t]- c_n^2\sum_{s,t,u=1}^KW_3^0[s,t,u]V_{2,2}^{(n)}[s,t]U_{2,2}^{(n)}[u,r] + o(1)\\ \notag
    &=c_n \operatorname{tr}W_2^0V_{3,3}^{(n)}[r,\cdot, \cdot] - c_n^2\operatorname{vec}(V_{2,2}^{(n)})^\top W_3^0U_{2,2}^{(n)}[\cdot,r] + o(1),
\end{align}
where $V_{3,3}^{(n)}[r,\cdot,\cdot]$ refers to the $K\times K$ matrix keeping the first coordinate fixed at $r$, $U_{2,2}^{(n)}[\cdot,r]$ refers to the $r-$th column of $U_{2,2}$ and $W_3^0(\cdot, \cdot, \cdot)$ is \textit{reshaped} as a $K^2 \times K$ matrix.

So all that remains is to derive expressions to approximate $V_{1,2}, V_{2,2}, V_{2,3}$ and $V_{3,3}$ in terms of \textbf{type-U} and \textbf{type-W} objects. The next lemma will be critical in determining the remaining two terms.
\begin{lemma}\label{lemma:tracial_expansion_bias}
Fix $r \in [K]$. Under Assumption, for any deterministic matrix $M$ with bounded spectral norm, we have
$$\frac{1}{n}\operatorname{tr}R\textbf{A}_jQ = \frac{1}{n}R\textbf{A}_jR + o(1), \quad \text{and,}$$
\begin{align*}
    \operatorname{tr}(Q(z)A_rQ(z)M) &= \operatorname{tr}R(z)A_rR(z)M + a(z)\\ &+ \sum_{s,t=1}^K\bigg(\frac{1}{n}\sum_{j=1}^n\psi_j^2\theta_{sj}\theta_{tj}\bigg) \bigl(\mathbb{E}\operatorname{tr}Q(z)A_rQ(z)A_s\bigr)\bigl(\frac{1}{n}\mathbb{E}\operatorname{tr}R(z)MR(z)A_t\bigr),
\end{align*} 
where $\mathbb{E}|a(z)| = O(n^{\frac{1}{2}})$.
\end{lemma}

Recall the population versions of $U_{1,2}^{(n)}, U_{2,2}^{(n)}$,$U_{2,3}^{(n)}$,$U_{3,3}^{(n)}$ from Subsection \ref{sec:Glossary}. Similar to the proof of Lemma \ref{lemma:Uj_U0_U_element_wise_difference}, we can show that 
\begin{align*}
  &\|U_{1,2}^{(n)}(z) - U_{1,2}(z)\| = o(1),\\  
  &\|U_{2,2}^{(n)}(z) - U_{2,2}(z)\| = o(1),\\  
  &\|U_{2,3}^{(n)}(z) - U_{2,3}(z)\| = o(1).
\end{align*}
Moreover, $|U_{3,3}^{(n)}[r,s,t] - U_{3,3}[r,s,t]| = o(1)$ for any $r,s,t \in [K]$. As an immediate consequence of Lemma \ref{lemma:tracial_expansion_bias} for $r_0,s_0,t_0 \in [K]$, we have the following:
\begin{align}\label{relationship_V_U_W}
    &V_{2,2}^{(n)}[r_0,s_0] = U_{2,2}^{(n)}[r_0, s_0] + c_n\sum_{r,s = 1}^K V_{2,2}^{(n)}[r_0, r]W_2^0[r,s]U_{2,2}^{(n)}[s, s_0] + a(z),\\
    &V_{2,3}^{(n)}[r_0,s_0] = U_{2,3}^{(n)}[r_0,s_0] + c_n\sum_{r,s=1}^K V_{2,2}^{(n)}[r_0,r]W_2^0[r,s]U_{2,3}^{(n)}[s,s_0] + a(z),\\
    &V_{3,3}^{(n)}[r_0,s_0,t_0] = U_{3,3}^{(n)}[r_0,s_0,t_0] + c_n\sum_{r,s,t=1}^K V_{2,2}^{(n)}[r_0,r]W_2^0[r,s]U_{3,3}^{(n)}[s, s_0, t_0] + a(z),
\end{align}
where $\mathbb{E}|a(z)| = O(n^{\frac{1}{2}})$. Note that from the first equation in (\ref{relationship_V_U_W}), for large $n$, we have the following:
\begin{align*}
    V_{2,2}^{(n)} &= U_{2,2}^{(n)} + V_{2,2}^{(n)}W_2^0U_{2,2}^{(n)} + A(z) = U_{2,2}^{(n)}(I - c_nW_2^0U_{2,2}^{(n)})^{-1} + A(z)\\
    &= (I - c_nU_{2,2}^{(n)}W_2^0)^{-1}U_{2,2}^{(n)} + A(z).
\end{align*}
Recall that the existence of the inverse was established in Lemma \ref{lemma:spectral_radius_le_1}. Therefore, we can \textit{plugin} $V_{2,2}$ into the remaining two equations to completely determine $V_{2,3}$ and $V_{3,3}$ upto a remainder term. Using this in (\ref{dnz_intermsof_V_W}) and (\ref{drnz_intermsof_V_W}), we see that 
\begin{align}\label{d0n_to_d0_drn_to_dr_convergence}
  d_{0n}(z) \longrightarrow d_0(z), \quad d_{rn}(z) \longrightarrow d_r(z),  
\end{align}
where $d(z), d_r(z)$ are stated in (\ref{defining_d0_dr}).

\subsection{Analysis of \texorpdfstring{$e_{0n}(z)$}{}}
\begin{figure}[!htbp]
\centering
 \begin{subfigure}{0.75\textwidth}
 \centering
\includegraphics[width=\textwidth,height=6cm]{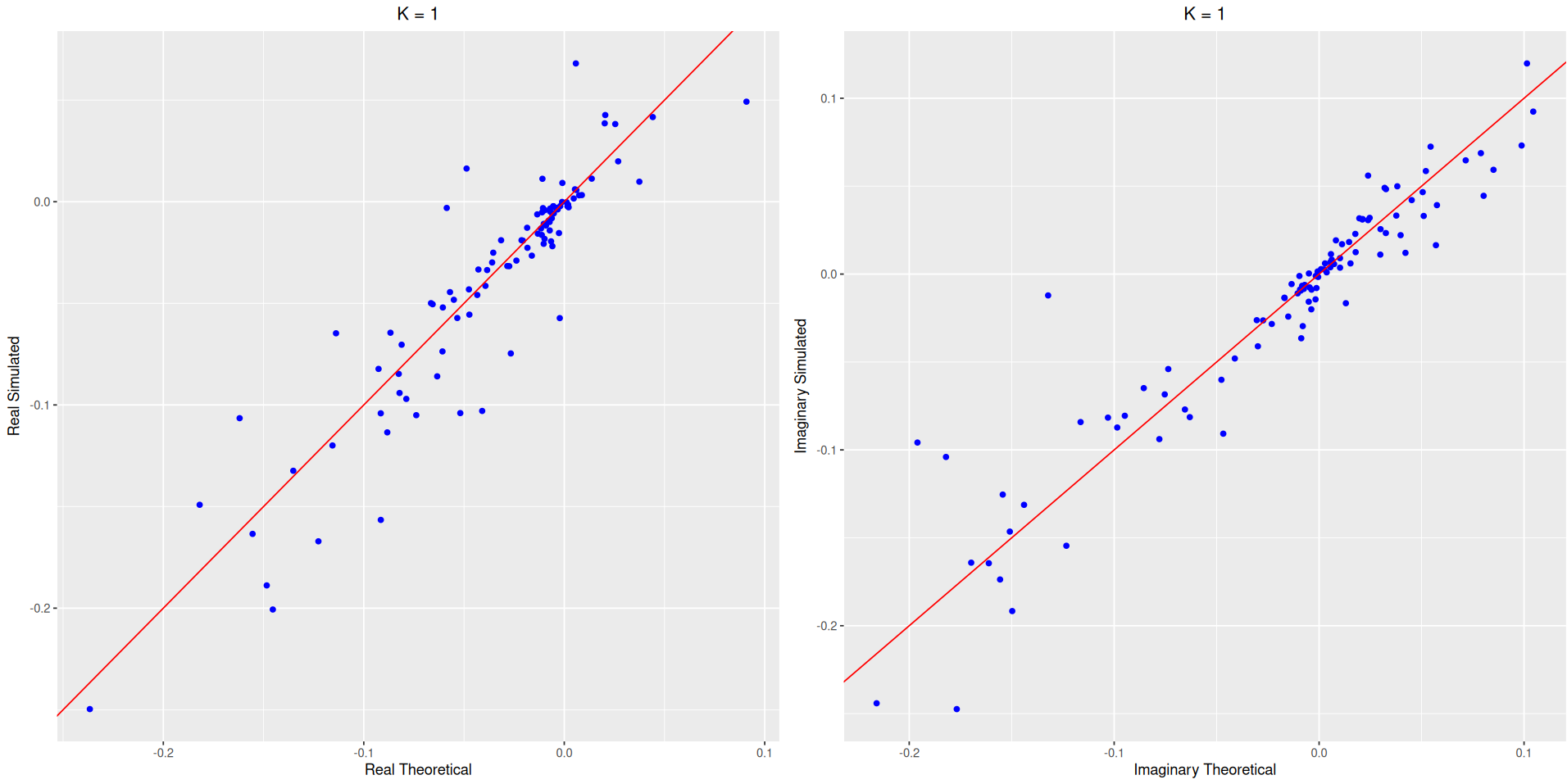}
     \caption{K=1:c=1}
     \label{fig7}
 \end{subfigure}
 \medskip
 \begin{subfigure}{0.75\textwidth}
\includegraphics[width=\textwidth,height=6cm]{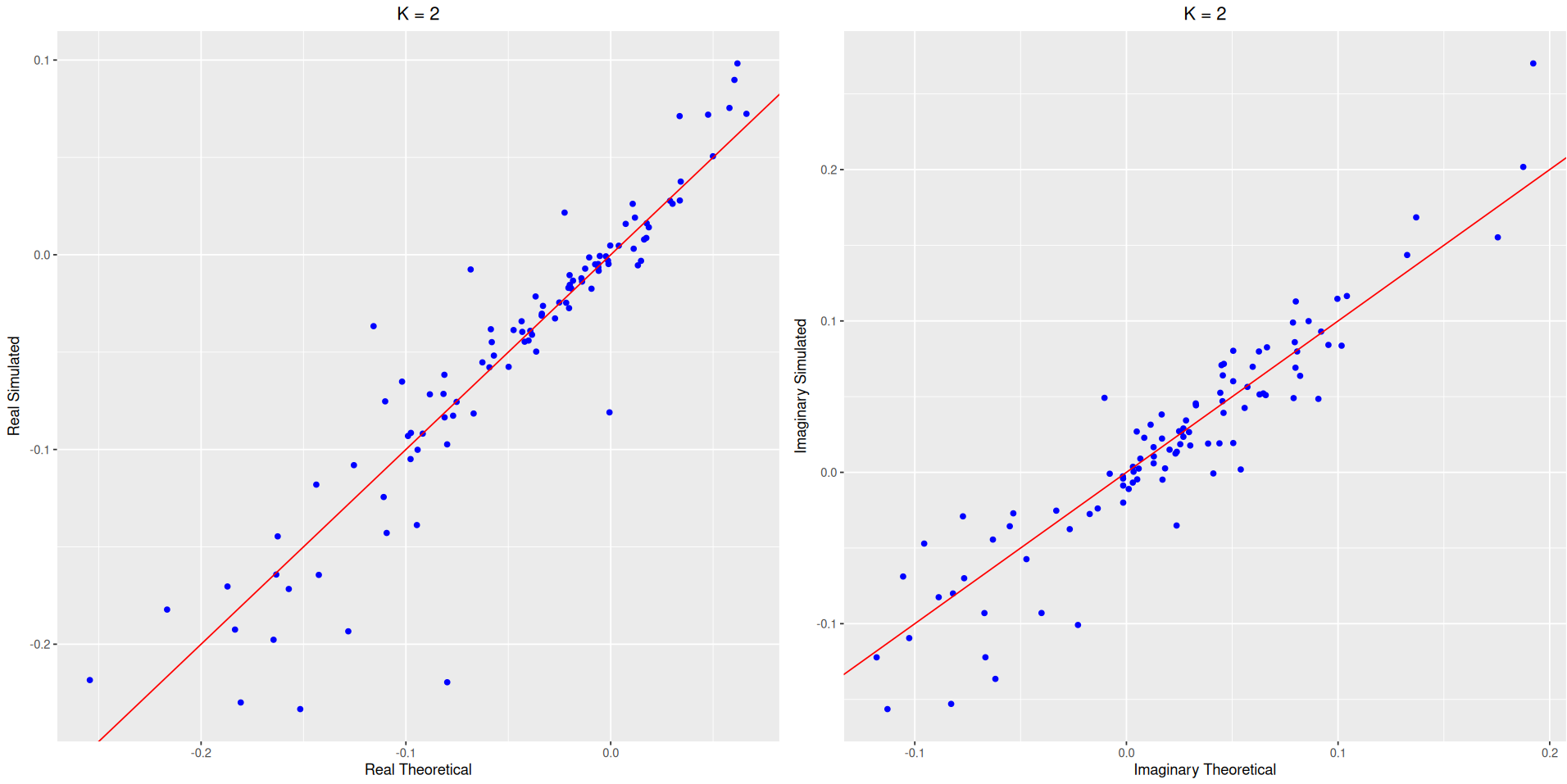}
     \caption{K=2:c=1}
     \label{fig8}
 \end{subfigure}
\medskip
 \begin{subfigure}{0.75\textwidth}
\includegraphics[width=\textwidth,height=6cm]{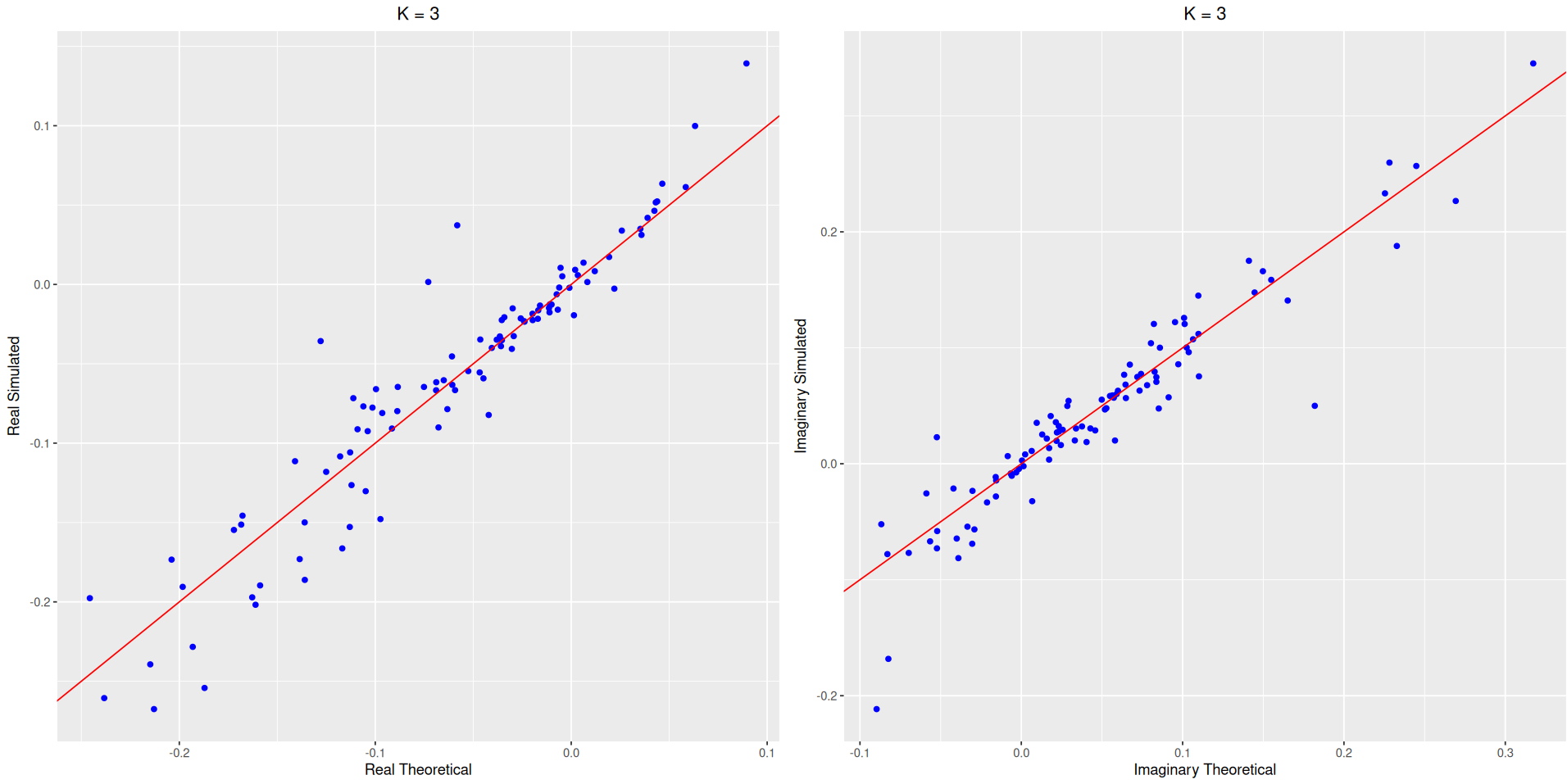}
     \caption{K=3:c=1}
     \label{fig9}
 \end{subfigure}
 \caption{Scatterplot of Simulated vs. Theoretical values (real and imaginary parts in the left and right panels respectively) of $M_{n2}(z)$ calculated at 100 values of $z$ for $K=1,2,3$ and $c=1$ with standard Gaussian innovations.}
 \label{fig1_to_6_bias}
\end{figure}

\begin{figure}
\centering
 \begin{subfigure}{0.75\textwidth}
\includegraphics[width=\textwidth,height=6cm]{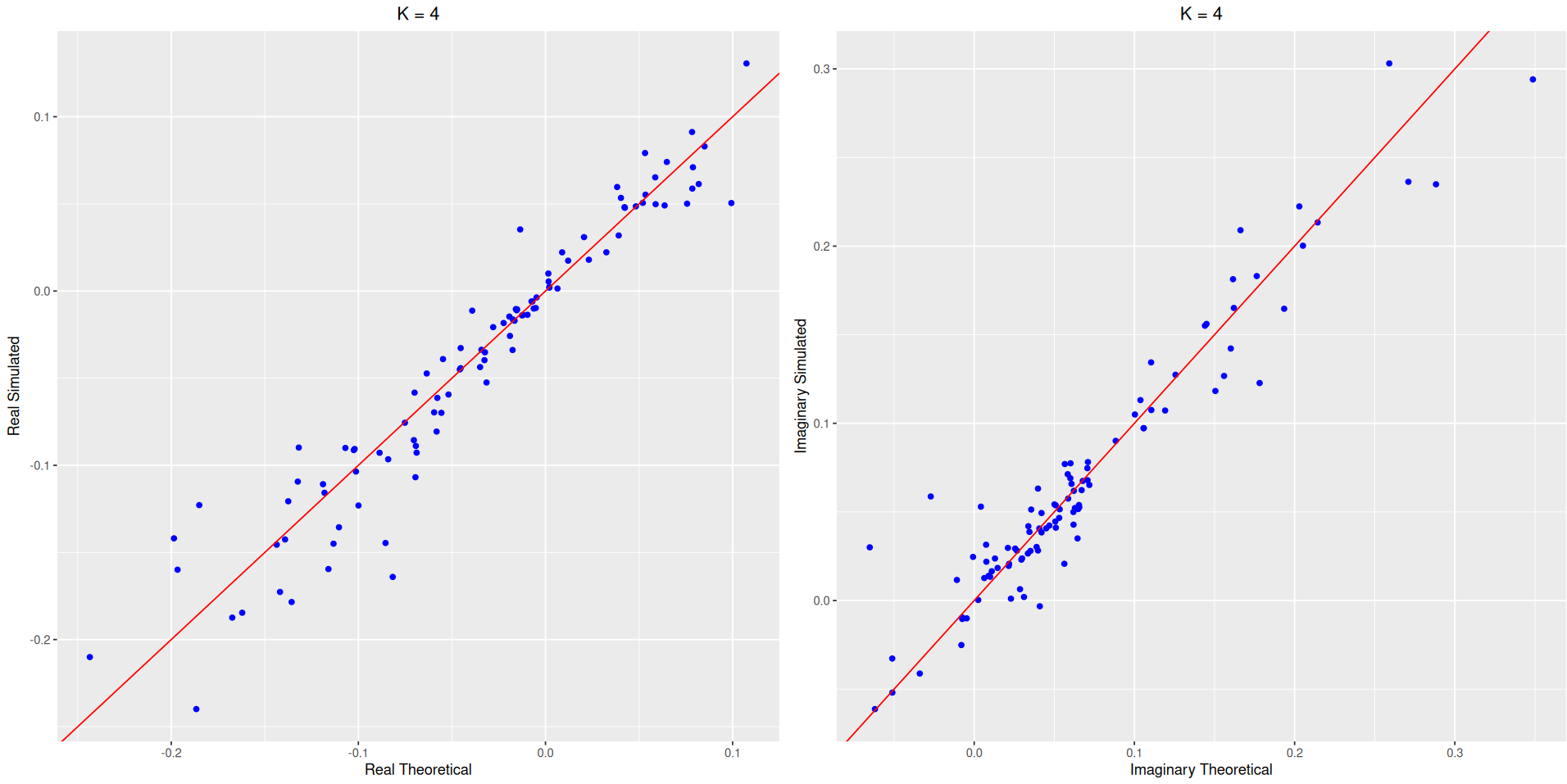}
     \caption{K=4:c=1}
     \label{fig10}
 \end{subfigure}
 \medskip
   \begin{subfigure}{0.75\textwidth}     \includegraphics[width=\textwidth,height=6cm]{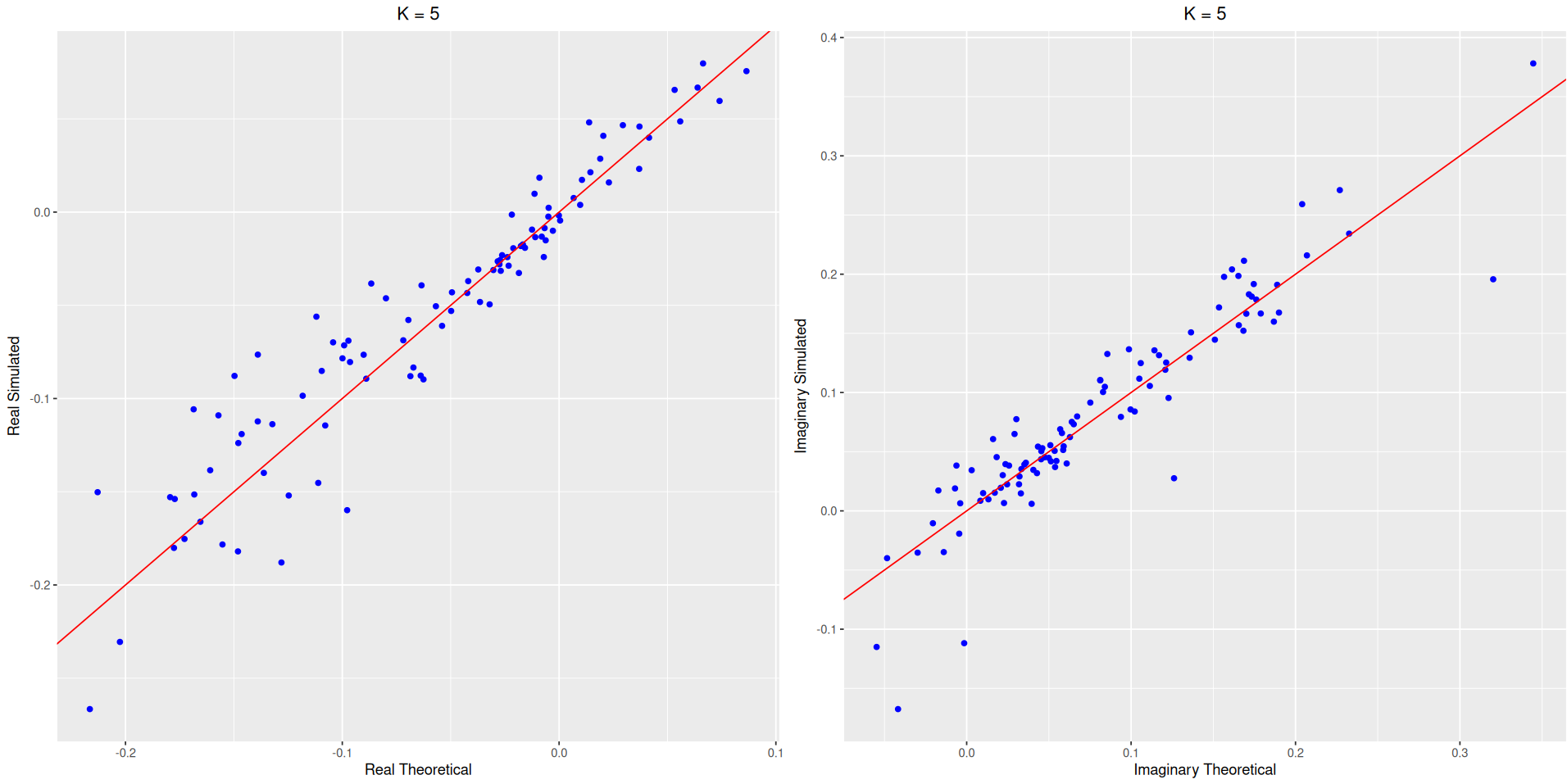}
     \caption{K=5:c=1}
     \label{fig11}
 \end{subfigure}
 \medskip
 \begin{subfigure}{0.75\textwidth}
\includegraphics[width=\textwidth,height=6cm]{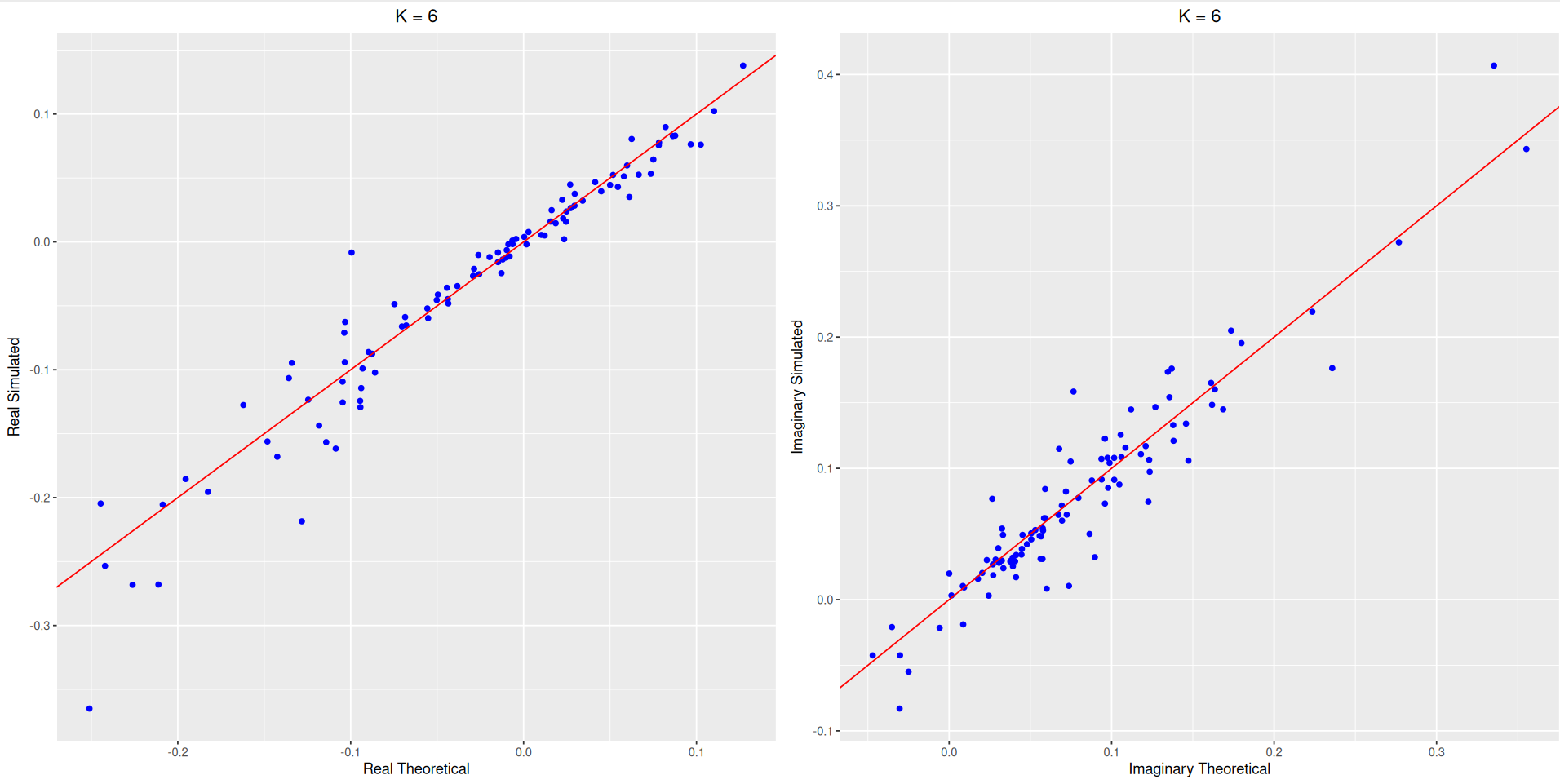}
     \caption{K=6:c=1}
     \label{fig12}
 \end{subfigure}
 \caption{Scatterplot of Simulated vs. Theoretical values (real and imaginary parts in the left and right panels respectively) of $M_{n2}(z)$ calculated at 100 values of $z$ for $K=4,5,6$ and $c=1$ with standard Gaussian innovations.}
  \label{fig7_to_12_bias}
\end{figure}

The analysis of $e_{0n}(z)$ depends on $d_{rn}(z)$ established in the previous subsection. Using the resolvent identity (\ref{Resolvent_Identity}), we have 
\begin{align}\label{decomposing_Mn2_1}
    e_{0n}(z)&=\operatorname{tr}(R(z) - R^0(z))\\\notag
    &= \operatorname{tr}\bigg(-zI_p + \frac{1}{n}\sum_{j=1}^n\psi_j\textbf{A}_j\bigg)^{-1} + \operatorname{tr}\bigg(zI_p + \sum_{r=1}^Kzg_{rn}^0A_r\bigg)^{-1}\\\notag
    &= - \sum_{r=1}^K \bigl(zg_{rn}^0(z) + \frac{1}{n}\sum_{j=1}^n\psi_j\theta_{rj}\bigr)\operatorname{tr}R(z)A_rR^0(z).
    \end{align}
Using properties of finite sample solutions from (\ref{finite_sample_properties}), we have
\begin{align*}
    zg_{rn}^0(z) + \frac{1}{n}\sum_{j=1}^n\psi_j\theta_{rj}
    &= \frac{1}{n}\sum_{j=1}^n\theta_{rj}(\psi_j - b_j^0)\\ 
    &=\frac{1}{n}\sum_{j=1}^n \bigg(\frac{\theta_{rj}}{1 + c_n\boldsymbol{\theta}_j^\top \mathbb{E}\textbf{h}_n(z)} - \frac{\theta_{rj}}{1 + c_n\boldsymbol{\theta}_j^\top \textbf{h}_n^0(z)}\bigg)\\
    &= -c_n\sum_{s=1}^K (\mathbb{E}h_{sn}(z)-h_{sn}^0(z))\bigg(\frac{1}{n}\sum_{j=1}^n\theta_{rj}\theta_{sj}\psi_j(z)b_j^0(z)\bigg).
\end{align*}
From Lemma \ref{lemma:bounds_tightness_section}, we have
\begin{align*}
    \sup_{n;j\in [n],z \in \mathcal{C}_n}|\psi_j(z) - b_j^0(z)| = o(1).
\end{align*}
Therefore, 
\begin{align*}
    \frac{1}{n}\sum_{j=1}^n\theta_{rj}\theta_{sj}\psi_j^2= \frac{1}{n}\sum_{j=1}^n\theta_{rj}\theta_{sj}\psi_j(z)b_j^0(z) + o(1).
\end{align*}
Denote $\boldsymbol{\mu}_n = \{\mu_{rn}\}_{r=1}^K$, where $\mu_{sn}(z) = p(\mathbb{E}h_{sn}(z) - h_{sn}^0(z))$. Combining everything into (\ref{decomposing_Mn2_1}), we have
\begin{align}\label{enz_final}
    e_{0n}(z) &=c_n\sum_{r=1}^K \sum_{s=1}^K \mu_{sn}(z)\bigg(\frac{1}{n}\sum_{j=1}^n\theta_{rj}\theta_{sj}\psi_j^2(z)\bigg)\bigg(\frac{1}{p}\operatorname{tr}R(z)A_rR^0(z)\bigg) + o(1)\\ \notag
    &= c_n\sum_{r,s=1}^K \mu_{sn}(z)W_2[r,s]U_{1,2}[r] + o(1)\\ \notag
    &= c_n \boldsymbol{\mu}_n^\top W_2U_{1,2} + o(1).
\end{align}
The expression will be complete once we solve for $\boldsymbol{\mu}_n$. Note that similar to the expansion for $e_{0n}(z)$ in (\ref{decomposing_Mn2_1}), we also have
\begin{align*}
    &p(\mathbb{E}h_{rn}(z) - h_{rn}^0(z))\\
    =&\, \operatorname{tr}A_r(\mathbb{E}Q(z) - R^0(z))\\
    =&\, \operatorname{tr}A_r(\mathbb{E}(Q - R) + \operatorname{tr}A_r(R - R^0)\\
    =&\, d_{rn}(z) + e_{rn}(z)\\
    =&\, d_{rn}(z) + c_n\sum_{s=1}^K \sum_{t=1}^K \mu_{tn}(z)\bigg(\frac{1}{n}\sum_{j=1}^n\theta_{sj}\theta_{tj}\psi_j(z)b_j^0(z)\bigg)\bigg(\frac{1}{p}\operatorname{tr}R(z)A_sR^0(z)A_r\bigg)\\
    =&\, d_{rn}(z) + c_n\sum_{s,t=1}^K U_{2,2}[r,s]W_2[s,t]\mu_{tn}(z) + o(1).
\end{align*}
 Then, we have the following system of equations:
\begin{align*}
  &\boldsymbol{\mu}_n = \textbf{d}_n + c_nU_{2,2}W_2\boldsymbol{\mu}_n + o(1) \implies \boldsymbol{\mu}_n = (I - c_nU_{2,2}W_2)^{-1}\textbf{d}_n + o(1).  
\end{align*}
Since $\textbf{d}_n$ has already been identified in (\ref{drnz_intermsof_V_W}), the expression for $e_{0n}(z)$ in (\ref{enz_final}) is complete. Having already established that $\textbf{d}_n = o(1)$ in the complex case, we immediately have $\boldsymbol{\mu}_n = o(1)$ entry-wise and therefore $$M_{n2}(z) = o(1).$$
In the real case, we have 
$$M_{n2}(z) = d_{0n}(z) + c_n\textbf{d}_n^\top (I 
- c_nW_2U_{2,2})^{-1} W_2U_{1,2} + o(1).$$
Using (\ref{d0n_to_d0_drn_to_dr_convergence}), we have finally proved that $M_{n2}(z) \to m_2(z)$.  A few numerical simulations (see Figures \ref{fig1_to_6_bias}, \ref{fig7_to_12_bias}) have been added for validation.
\newpage
\subsection{Numerical Validation of the bias term}

\newpage
\printbibliography

@article{BaiLiPan2019,
    author = {Bai,Z. and Li, H. and Pan, G.},
    title = {Central limit theorem for linear spectral
statistics of large dimensional separable
sample covariance matrices},
    journal = {Bernoulli},
    volume = {25},
    pages = {1838--1869},
    year = {2019}
}

@article{GeroHill03,
    author = {Geronimo, J.S. and Hill, T.P.},
    title = {Necessary \& sufficient condition that the limit of {Stieltjes} transforms is a {Stieltjes} transform},
    journal = {J. Approx. Theory},
    volume = {121}, 
    issue = {1},
    pages = {54--60},
    year = {2003}
}

@article{PaulSilverstein2009,
    author = {Paul, D. and Silverstein, J.W.},
    title = {No eigenvalues outside the support of the limiting empirical spectral distribution of a separable covariance matrix},
    journal = {J. Multivar. Anal.},
    volume = {100},
    number = {1}, 
    pages = {37-57},
    year = {2009}
}

@phdthesis{Lixin06,
    author = {Zhang, L.},
    title = {Spectral Analysis of Large Dimensional Random matrices},
    school = {National University of Singapore},
    year = {2006}
}

@article{SepCovar,
    author = {Couillet, R. and Hachem, W.},
    title = {Analysis of the limiting spectral measure of large random matrices of the separable covariance type},
    journal = {Random Matrices: Theory Appl},
    volume = {},
    pages = {},
    year = {2014}
}

@article{Lindeberg,
    author = {Lindeberg, J. W.},
    title = {Eine neue herleitung des exponentialgesetzes in der
wahrscheinlichkeitsrechnung},
    journal = {Math. Z.},
    year = {1922}
}

@article{NajimYao,
    author = {Najim, J. and Yao, J.},
    title = {Gaussian Fluctuations For Linear Spectral Statistics of Large Random Covariance Matrices},
    journal = {The Annals of Applied Probability},
    year = {2016}
}

@article{BaiSilv04,
    author = {Bai, Z. and Silverstein, J.W.},
    title = {CLT For Linear Spectral Statistics of Large-Dimensional Sample Covariance Matrices},
    journal = {The Annals of Applied Probability},
    year = {2004}
}

@article{JHDP26,
    author = {Hazarika, J. and Paul, D.},
    title = {LSD of sample covariances of superposition of matrices with separable covariance structure},
    journal = {Statistical Learning and Data Science},
    year = {2026}
}

@article{XieJohnstone,
    author = {Xie, R. and Johnstone, I.},
    title = {CLT for Linear Spectral Statistics in High-Dimensional
Random Effects Models},
    journal = {arXiv preprint},
    year = {2024}
}

@article{BaiYin,
    author = {Bai, Z.D. and Yin, Y.Q.},
    title = {Limit of the smallest eigenvalue of a large-dimensional sample covariance matrix.},
    journal = {Annals of Probability},
    year = {1993}
}

@article{BaiYinKri,
    author = {Yin, Y.Q. and Bai, Z.D. and Krishnaiah, P.R.},
    title = {On the limit of the largest eigenvalue of the largedimensional sample covariance matrix.},
    journal = {Probab. Theory Related Fields},
    year = {1988}
}

@article{WangPanJiang,
    author = {Wang, R. and Pan, G. and Jiang, D.},
    title = {Phase Transition of Spectral Fluctuations in Large Gram Matrices with a Variance Profile: A Unified Framework for Sparse CLTs},
    journal = {arXiv preprint},
    year = {2026}
}
\appendix

\section{A Few Established Results}
\textbf{Resolvent Identity:} For invertible matrices $A$ and $B$, 
    \begin{align}\label{Resolvent_Identity}
        A^{-1} - B^{-1} = A^{-1}(B - A)B^{-1} = B^{-1}(B-A)A^{-1}.
    \end{align}
\textbf{Tracial Inequalities:}
\begin{align}\label{tr_ABCD}
    \operatorname{tr}(ABCD) \leq \|A\|\,\|C\| \bigg(\operatorname{tr}(BB^*) \operatorname{tr}(DD^*)\bigg)^\frac{1}{2}.
\end{align}

\begin{lemma}\label{Lemma:Abel}
\textbf{Generalization of Abel's Lemma:} Let $\{\textbf{a}_n\}_{n=1}^\infty, \{\textbf{b}_n\}_{n=1}^\infty$ be sequences of matrices. Define the partial sums $A_n = \sum_{k=1}^n \textbf{a}_k, A_0 = \textbf{0}$. Then,
$$\sum_{k=1}^n \textbf{a}_k\textbf{b}_k  = A_n\textbf{b}_n - \sum_{k=1}^{n-1} A_k(\textbf{b}_{k+1}-\textbf{b}_k).$$
\end{lemma}
\begin{lemma}\label{lemma:oneRankTracePerturbation}
    (Lemma 0.11 of \cite{BaiLiPan2019}) For a symmetric matrix $M$ and $z \in \mathbb{C}_+$, the following holds: 
    \begin{align*}
        |\operatorname{tr}M(Q(z) - Q_j(z)| \leq \frac{\|M\|}{\Im(z)}.
    \end{align*}
\end{lemma}

\begin{lemma}\label{lemma:covariance_quadratic_forms}
Let $A,B$ be deterministic $p \times p$ matrices and $x \in \mathbb{C}^p$ be a random vector with i.i.d. entries with $\mathbb{E}x_1 = 0, \mathbb{E}|x_1|^2=\sigma^2, \mathcal{V} = \mathbb{E}x_1^2$. Define $\kappa = \mathbb{E}|x_1|^4 - |\mathcal{V}|^2 - 2\sigma^4.$ Then,
    \begin{align*}
    \mathbb{E}(x^*Ax - \operatorname{tr}A)(x^*Bx - \operatorname{tr}B)
    =\operatorname{tr}AB + |\mathcal{V}|^2\operatorname{tr}AB^\top  + \kappa \sum_i A_{ii}B_{ii}.
\end{align*}
\end{lemma}
\begin{proof}
We have 
\begin{align}\label{covariance_QF_expansion1}
    \mathbb{E}(x^*Ax - \operatorname{tr}A)(x^*Bx - \operatorname{tr}B) = \mathbb{E}(x^*Ax)(x^*Bx) - 2\sigma^4\operatorname{tr}A\operatorname{tr}B+ \operatorname{tr}A\operatorname{tr}B.
\end{align}
Denote $D = \sum_{i=1}^p A_{ii} B_{ii}$. The first term on the RHS can be expanded as follows.
\begin{align*}
&\mathbb{E}(x^*Ax)(x^*Bx)\\
=&\mathbb{E} \sum_{i,j,k,l}\bar{x}_ix_j\bar{x}_kx_lA_{ij}B_{kl}\\
    =& \sum_{i=j=k=l} \mathbb{E}|x_1|^4A_{ii}B_{ii}+ \mathcal{}\sum_{i=j\neq k=l}\sigma^4A_{ii}B_{kk} + \sum_{i=k\neq j=l}|\mathcal{V}|^2A_{ij}B_{ij} + \sum_{i=l\neq j=k} \sigma^4 A_{ij}B_{ji}\\
    =& \mathbb{E}|x_1|^4\sum_{i=1}^p A_{ii}B_{ii} + \sigma^4\sum_{i\neq j} A_{ii}B_{jj} + |\mathcal{V}|^2\sum_{i\neq j}A_{ij}B_{ij} + \sigma^4\sum_{i\neq j}A_{ij}B_{ji}\\
    =&\mathbb{E}|x_1|^4D + \sigma^4(\operatorname{tr}A\operatorname{tr}B - D) + |\mathcal{V}|^2(\operatorname{tr}AB^\top  - D) +\sigma^4(\operatorname{tr}AB - D)\\
    =& \kappa D + |\mathcal{V}|^2\operatorname{tr}AB^\top  + 2\sigma^4\operatorname{tr}A\operatorname{tr}B.
\end{align*}
\end{proof}

\begin{lemma}\label{lemma:tightness_condition}
    \textbf{Tightness of Random Processes:} Suppose $X^{(n)}$ are continuous random processes over $\mathbb{C}$. If there exist constants $\alpha, \beta,c>0$ with $\alpha\beta>1$ such that
    $$\mathbb{E}|X^{(n)}(z_1) - X^{(n)}(z_2)|^\alpha \leq c|z_1 - z_2|^{1+\beta}, \forall \, z_1,z_2 \in \mathbb{C}, \forall \, n,$$
    then $\{X^{(n)}\}$ is tight.
\end{lemma}

\begin{lemma}\label{lemma:Burkholder}
(Lemma B.26 in \cite{BaiSilv04}). Let $A = (a_{jk})$ be an $n \times n$ nonrandom matrix and $x = (x_1,\ldots,x_n)^\top $ be a random vector of independent entries. Assume that $\mathbb{E}x_j = 0$, $\mathbb{E}|x_j|^2=1$ and $\mathbb{E}|x_j|^l \leq \nu_l$. Then for $p \geq 1$, there exists a constant $C_p$ depending only on $p$ such that
    $$\mathbb{E}|x^*Ax - \operatorname{tr}A|^p \leq C_p[(\nu_4\operatorname{tr}AA^*)^\frac{p}{2} + \nu_{2p}\operatorname{tr}(AA^*)^\frac{p}{2}]$$
\end{lemma}

\begin{theorem}\label{theorem:martingale_CLT}
(Lemma 2.4 in \cite{BaiSilv04}). Suppose for each $n$, $Y_{n1},Y_{n2}, \ldots, Y_{nr_n}$ is a real martingale difference sequence with respect to the increasing $\sigma$-field $\{\mathcal{F}_{nj}\}$ having second moments. If as $n \rightarrow \infty$,
\begin{enumerate}
    \item[(i)] $$\sum_{j=1}^{r_n} \mathbb{E}(Y_{nj}^2|\mathcal{F}_{n,j-1}) \xrightarrow{i.p.} \sigma^2,$$
    where $\sigma^2$ is a positive constant, and for each $\epsilon \geq 0$,
    \item[(ii)] $$\sum_{j=1}^{r_n} \mathbb{E}(Y_{nj}^2I(|Y_{nj}| \geq \epsilon) \rightarrow 0,$$
\end{enumerate}
    then $$\sum_{j=1}^{r_n}Y_{nj} \xrightarrow{D} N(0, \sigma^2)$$
\end{theorem}

\begin{lemma}\label{lemma:Vitali_convergence}
(Lemma 2.3 of \cite{BaiSilv04}) Let $f_1, f_2, \ldots$ be analytic in $D$, a connected open set of $\mathbb{C}$ satisfying $f_n(z)| \leq M$ for every $n$ and $z \in D$, and $f_n(z)$ converges as $n \rightarrow  \infty$ for each $z$ in a subset of $D$ having a limit point of $D$. Then, there exists a function $f$, analytic in $D$ for which $f_n(z)\rightarrow f(z)$ and $f_n\prime(z)\rightarrow f^\prime(z)$ for all $z \in D$. Moreover, on any set bounded by a contour interior to $D$ , the convergence is uniform and $\{f_n^\prime(z)\}$ is uniformly bounded by $2M/\epsilon$, where $\epsilon$ is the distance between the contour and the boundary of $D$.
\end{lemma}

\begin{lemma}\label{lemma:Hoffman_Wielandt}
   \textbf{Hoffman-Wielandt Inequality:}For matrices $A,B \in \mathbb{C}^{p\times n}$, we have
    \begin{align*}
        |\sum_{i=1}^p \lambda_i(AA^*) - \lambda_i(BB^*)| \leq 2\sqrt{p}\|A - B\|_F\bigl(\lambda_{max}(AA^*) + \lambda_{max}(BB^*)\bigr)^\frac{1}{2}.
    \end{align*}
\end{lemma}

\subsection{Glossary} \label{sec:Glossary} 
 This a complete list of the ingredients for the mean and the covariance functions stated in Theorem \ref{mainTheorem_CLT}. $\textbf{h}(z_i), \textbf{g}(z_i)$ are the unique solutions of the matrix Dyson system from (\ref{h_main_equation}) and (\ref{g_main_equation}) respectively corresponding to $c, H$ and $G$. $W_2, W_3$ are maps that serve as \textit{weights} for certain tracial quantities. The motivation for the subscript should be clear from the definition. 
\begin{align}\label{defining_population_W_matrices}
    &W_2: \mathbb{C} \times \mathbb{C} \mapsto \mathbb{C}^{K \times K}, \quad W_2(z_1, z_2)[r,s] = \int \frac{\theta_r\theta_s \,dG(\boldsymbol{\theta})}{\prod_{i=1}^2(1 + c\boldsymbol{\theta}^\top\textbf{h}(z_i))}\\ \notag
    &W_3: \mathbb{C} \mapsto \mathbb{C}^{K \times K \times K}, \quad W_3(z)[r,s,t] = \int \frac{\theta_r\theta_s\theta_t\,dG(\boldsymbol{\theta})}{(1 + c\boldsymbol{\theta}^\top\textbf{h}(z))^3}.
\end{align}
$U_{1,2}, U_{2,2}, U_{2,3}, U_{3,3}$ are maps which are average traces of products of p.s.d. matrices ($A_r$'s) and resolvent-type objects (i.e. $Q(z)$ or its deterministic estimate). 
\begin{align}\label{defining_population_U_matrices}
        &U_{2,2}: \mathbb{C} \times \mathbb{C} \mapsto \mathbb{C}^{K \times K}, \quad U_{2,2}(z_1, z_2)[r,s] := \int \frac{\lambda_r\lambda_sdH(\boldsymbol{\lambda})}{\prod_{i=1}^2z_i(1 + \boldsymbol{\lambda}^\top\textbf{g}(z_i))}\\ \notag
        &U_{1,2}: \mathbb{C} \mapsto \mathbb{C}^{K}, \quad U_{1,2}(z)[r] = \int \frac{\lambda_rdH(\boldsymbol{\lambda})}{\bigl(-z(1 + \boldsymbol{\lambda}^\top \textbf{g}(z))\bigr)^2}\\ \notag
        &U_{2,3}: \mathbb{C} \mapsto \mathbb{C}^{K \times K}, \quad U_{2,3}(z)[r,s] := \int \frac{\lambda_r\lambda_sdH(\boldsymbol{\lambda})}{\bigl(-z(1 + \boldsymbol{\lambda}^\top\textbf{g}(z))\bigr)^3}\\ \notag
        &U_{3,3}: \mathbb{C} \mapsto \mathbb{C}^{K \times K \times K}, \quad U_{3,3}(z)[r,s,t] := \int \frac{\lambda_r\lambda_s\lambda_tdH(\boldsymbol{\lambda})}{\bigl(-z(1 + \boldsymbol{\lambda}^\top\textbf{g}(z))\bigr)^3}.
    \end{align}
\begin{remark}
    The second index in the subscript of $U$ indicates the number of \textit{resolvent-type} objects inside the trace operator while the first index indicates the number of $A_r$ matrices involved. The discrepancy in the number of complex arguments in these maps (i.e. only $z$ vs both $z_1$ and $z_2$) is to indicate whether they are used only in the bias function or both bias and covariance functions.
\end{remark}
Now, we define $V_{2,2}, V_{2,3}, V_{3,3}$ which are similar to their \textit{$U$-counterparts}. Note that once $V_{2,2}$, which is a function of $U_{2,2}$ and $W_2$ is derived, it gets plugged into $V_{2,3}$ and $V_{3,3}$. 
\begin{align}\label{defining_population_V_matrices1}
    &V_{2,2}(z_1, z_2) = U_{2,2}(z_1, z_2)\bigl(I - cW_2(z_1, z_2)U_{2,2}(z_1, z_2)\bigr)^{-1}.
\end{align}
We now define $V_{2,3} \in \mathbb{C}^{K\times K}$ and $V_{3,3} \in \mathbb{C}^{K\times K \times K}$ as functions of a single argument ($z$). Writing $U_{2,3}, U_{3,3},W_{2,2}, V_{2,2}$ instead of $U_{2,3}(z), U_{3,3}(z),W_{2,2}(z,z), V_{2,2}(z,z)$ respectively, we also define the quantities for $r_0,s_0,t_0\in [K]$ as follows:
\begin{align}\label{defining_population_V_matrices2}
    &V_{2,3}(z)[r_0,s_0] = U_{2,3}[r_0,s_0] + c\sum_{r,s=1}^K V_{2,2}[r_0,r]W_2[r,s]U_{2,3}[s,s_0],\\ \notag
    &V_{3,3}(z)[r_0,s_0,t_0] = U_{3,3}[r_0,s_0,t_0] + c\sum_{r,s,t=1}^K V_{2,2}[r_0,r]W_2[r,s]U_{3,3}[s, s_0, t_0].
\end{align}
Now, define $d_0(z)$ and $\textbf{d}(z) = \{d_r(z)\}_{r=1}^K$ as follows:
\begin{align}\label{defining_d0_dr}
        &d_0(z) = c\operatorname{tr}W_2(z,z)V_{2,3}(z) - c^2\operatorname{vec}(V_{2,2}(z,z))^\top W_3(z)U_{1,2}(z),\\ \notag
        &d_r(z) = c\operatorname{tr}(W_2(z,z)V_{3,3}(z)[r,\cdot, \cdot]) - c^2\operatorname{vec}(V_{2,2}(z,z))^\top W_3(z)U_{2,2}(z)[,r].
\end{align}
Here, $W_3$ is \textit{reshaped} into a $K^2\times K$ matrix for compatibility. Finally, the \textbf{mean function} is
\begin{align}\label{defining_population_mean_function}
  m_2(z) &= d_0(z) - \textbf{d}^\top(z)\bigl(I - cW_2(z,z)U_{2,2}(z,z)\bigr)^{-1}W_2(z,z)U_{1,2}(z)\\ \notag
  &= d_0(z) - U_{1,2}(z)^\top W_2(z,z) \bigl(I - cU_{2,2}(z,z)W_2(z, z)\bigr)^{-1} \textbf{d}(z),
\end{align}
and the \textbf{covariance kernel} is
\begin{align}\label{defining_population_covariance_function}
  \operatorname{ker}(z_1, z_2) = -\frac{\partial^2}{\partial z_2 \partial  z_1}\operatorname{tr} \log \bigl(I - cU_{2,2}(z_1,z_2)W_2(z_1, z_2)\bigr).
\end{align}

\subsection{Separable case (\texorpdfstring{$K=1$}{}): Expression matches with Bai, Li and Pan (2019)}
\begin{proof}
Denote the sample covariance matrix of \cite{BaiLiPan2019} is $\widetilde{S}_n := \frac{1}{p}XX^*,$
and the difference in normalization factor implies that
$$S_n = \frac{p}{n}\widetilde{S}_n = c_n\widetilde{S}_n.$$
As a result, the support of the LSD grows by a factor of $c$. We will write $\textbf{h}_n(z), \textbf{h}(z), \textbf{g}_n(z), \textbf{g}(z)$ as $h_n(z), h(z),g_n(z),g(z)$ respectively since $K = 1$. 
We will now express the relationships between $h_n, g_n (h, g)$ and $g_{2n}, g_{1n}(g_2, g_1)$. Observe that 
$$h_n(z)= \frac{1}{p}\operatorname{tr}(c_n\widetilde{S}_n -zI_p)^{-1}A_n = \frac{1}{c_n}g_{2n}(\frac{z}{c_n}).$$
As a result, $h(z) = c^{-1}g_2(z/c)$. Similarly, we also deduce that $g(z) = g_1(z/c)$. 
Denote the analog of the random process $M_n(z)$ arising out of $\widetilde{S}_n$ as $\widetilde{M}_n(z)$. Then, it is clear that 
\begin{align*}
    &\operatorname{Cov}(M_n(z_1),M_n(z_2))\\ =&\, c_n^{-2}\operatorname{Cov}(\widetilde{M}_n(z_1/c_n),\widetilde{M}_n(z_2/c_n)) \longrightarrow - c^{-2}\frac{\partial^2}{\partial z_1 \partial z_2} \log(1 - d(z_1/c, z_2/c)).
\end{align*}
Note that, we only have to match the expression on the right apart from the factor of $c^{-2}$ since the support of the LSD has grown by a factor of $c$. Using the value of the expression $d(z_1, z_2)$ from Theorem 2.1 of \cite{BaiLiPan2019}, we obtain
be\begin{align*}
        d(z_1/c,z_2/c) &= \frac{c^2}{z_1z_2}\frac{z_1g_1(z_1/c) - z_2g_1(z_2/c)}{c(g_2(z_1/c)-g_2(z_2/c))}\frac{z_1g_2(z_1/c)-z_2g_2(z_2/c)}{c(g_1(z_1/c) - g_1(z_2/c))}\\
        &= \frac{1}{z_1z_2} \frac{z_1g(z_1) - z_2g(z_2)}{h(z_1)-h(z_2)}\frac{z_1h(z_1)-z_2h(z_2)}{g(z_1)-g(z_2)}\\
   &= \frac{1}{z_1z_2}\frac{\int\frac{dG(\theta)}{1+c\theta h(z_1)} - \int\frac{dG(\theta)}{1+c\theta h(z_2)}}{h(z_1) - h(z_2)}\frac{\int\frac{dH(\lambda)}{1+\lambda g(z_1)} - \int\frac{dG(\lambda)}{1+\lambda g(z_2)}}{g(z_1) - g(z_2)} \\
   &=c\int \frac{\lambda^2}{z_1z_2(1+\lambda g(z_1))(1+\lambda g(z_2))} \int\frac{\theta^2 dG(\theta)}{(1 + c\theta h(z_1))(1 + c\theta h(z_2))}\\ 
   &= cU_{2,2}(z_1,z_2)W_2(z_1,z_2).
\end{align*}
Now observe that 
\begin{align*}
    \int_0^{d(z_1,z_2)} \frac{dz}{1-z} = -\log(1-d(z_1,z_2)).
\end{align*}
Therefore the expression for covariance matches exactly with \cite{BaiLiPan2019} for $K=1$.

Denote the bias component of BLP (analog of our $m_2(z)$) as $\widetilde{m}_2(z)$. Then, it is clear that 
$$\widetilde{m}_2(z) = cm_2(cz).$$

First, we show that $d_{n1}(z)$ of \cite{BaiLiPan2019} converges to $cd_0(cz)$ of our manuscript. Define $\widetilde{c} := n/p$. Note that $\widetilde{c}$ is the analog of $c$ in \cite{BaiLiPan2019}. We now observe that for $K=1$, we have
\begin{align*}
    (I - cU_{2,2}(cz,cz)W_2(cz,cz))^{-1} &= \bigg(1 - c\int \frac{\theta^2 dG(\theta)}{(1 + c\theta h(cz))^2}\int \frac{\lambda^2 dH(\lambda)}{(-cz(1 + \lambda g(z)))^2}\bigg)^{-1}\\
    &= \bigg(1 - \frac{1}{c}\int \frac{\theta^2 dG(\theta)}{(1 + \theta g_2(z))^2}\int \frac{\lambda^2 dH(\lambda)}{(-z(1 + \lambda g_1(z)))^2}\bigg)^{-1}\\
    &= \bigg(1 - \widetilde{c}z^{-2}\int \frac{\theta^2 dG(\theta)}{(1 + \theta g_2(z))^2}\int \frac{\lambda^2 dH(\lambda)}{(1 + \lambda g_1(z))^2}\bigg)^{-1},
\end{align*}
which is an expression that appears in the definitions of $d_1(z), d_2(z)$ in \cite{BaiLiPan2019}. Observe that $V_{2,3}(z) = U_{2,3}(z)(1 - cU_{2,2}W_2(z))^{-1}$. As a result, we have
\begin{align*}
    d_0(cz) &= cW_2V_{2,3}(cz) - c^2V_{2,2}W_3U_{1,2}(cz)\\
    &= (cW_2U_{2,3}(cz) - c^2W_3U_{1,2}U_{2,2}(cz))(1 - cU_{2,2}W_2(cz))^{-1}\\
    &= c^{-2}(1 - cU_{2,2}W_2(cz))^{-1} \times \bigg(\int \frac{\theta^2 dG(\theta)}{(1 + \theta g_2(z))^2}\int \frac{\lambda^2 dH(\lambda)}{(-z(1 + \lambda g_1(z))^3}  \\
    &\,\,\,-\int \frac{\theta^3 dG(\theta)}{(1 + \theta g_2(z))^3}\int \frac{\lambda dH(\lambda)}{(-z(1 + \lambda g_1(z))^2}\int \frac{\lambda^2 dH(\lambda)}{(-z(1 + \lambda g_1(z))^2}\bigg).
\end{align*}
Therefore, $cd_0(cz)$ matches exactly with $d_1(z)$ of \cite{BaiLiPan2019}. 

Finally, working with the expression for $\widetilde{M}_{n2}(z)$ in page 7 of the supplement of \cite{BaiLiPan2019}, we find that
\begin{align*}
    \widetilde{M}_{n2}(z) =& d_{n1}(z) - n(\mathbb{E}g_{2n}(z) - g_{2n}^0(z))\frac{1}{n}\sum_{j=1}^n \frac{s_j\psi_j}{1 + s_jg_{2n}^0(z)}\\
    &\times \int \frac{x dx}{(\frac{1}{n}\sum_{j=1}^ns_j\psi_j(z)x - z)(z + zg_{1n}^0(z)x)} +o(1)\\
    &\longrightarrow cd_0(cz) - cn(\mathbb{E}h_n(cz)-h_n^0(cz))W_2(cz)c^2U_{1,2}(cz)\\
    &=c(d_0(cz) - c\mu(cz)W_2(cz)U_{1,2}(cz)) = cm_2(cz).
\end{align*}
\end{proof}

\subsection{Proof of (\ref{eq:ratio_moment})}\label{proof:eq:ratio_moment}
It suffices to show the following results:
\begin{align*}
    \mu_{1n} = \mathbb{E}(p^{-1}\operatorname{tr}S_n) = \alpha_1\beta_1, \quad 
    \mu_{2n} = \mathbb{E}(p^{-1}\operatorname{tr}S_n^2) = \alpha_2\beta_1^2 + c_n\alpha_1^2\beta_2 + \frac{1}{p}c_n\alpha_2\beta_2.
\end{align*}

The first one is a trivial result and the second one follows from Wick's formula (assuming $Z$ is a real matrix) as shown below:
\begin{align*}
    pn^2\mu_{2n} &= \mathbb{E}\operatorname{tr}AZBZ^\top AZBZ^\top\\
    &= \sum_{i,j,k,l,a,b,c,d} \mathbb{E}A_{ij}Z_{ja}B_{ab}Z_{kb}A_{kl}Z_{lc}B_{cd}Z_{id}\\
    &= \sum_{i,j,k,l,a,b,c,d} A_{ij}A_{kl}B_{ab}B_{cd}\mathbb{E}Z_{ja}Z_{kb}Z_{lc}Z_{id}\\
    &= \sum_{i,j,a,c}A_{ij}A_{ji}B_{aa}B_{cc}+ \sum_{i,j,a,b}A_{ij}^2B_{ab}^2  + \sum_{i,k,a,b}A_{ii}A_{kk}B_{ab}B_{ba}\\
    &= \operatorname{tr}(A^2)(\operatorname{tr}B)^2 + \operatorname{tr}(A^2)\operatorname{tr}(B^2) + (\operatorname{tr}A)^2\operatorname{tr}(B^2)\\
    &= pn^2\alpha_{2n} \beta_{1n}^2 + pn\alpha_{2n}\beta_{2n} + p^2n\alpha_{1n}^2\beta_{2n},
\end{align*}
where $\alpha_{jn}, \beta_{jn}$ are the empirical analogues of $\alpha_j, \beta_j$ respectively. Then, observe that for large $n$, we have
\begin{align*}
    \mu_{2n} \longrightarrow \alpha_2\beta_1^2 + c\alpha_1^2\beta_2.
\end{align*}

\section{Impact of Truncation, Centralization and Scaling}\label{sec:Truncation_Centralization_rescaling}
Here, the indices $i,j,r$ belong to $[p], [n], [K]$ respectively. Denote\begin{align*}
    a_{ij}^{(r)} := \mathbb{E}|z_{ij}^{(r)}|^4\mathbbm{1}_{\{|z_{ij}^{(r)}|>\delta_n\sqrt{n}\}},
\end{align*} 
From (\ref{L4}), for any $\delta>0$, we have 
\begin{align}\label{lindeberg_4_condition_delta}
    \frac{1}{n^2\delta^4}\sum_{i,j,r}a_{ij}^{(r)} \longrightarrow 0.
\end{align}
By \textit{standard} arguments, we can choose a sequence $\{\delta_n\}_{n=1}^\infty \downarrow 0$ such that the above holds even when $\delta$ has been replaced by $\delta_n$, i.e.
\begin{align}\label{lindeberg_4_condition_delta_n}
    \frac{1}{Kn^2\delta_n^4}\sum_{i,j,r}\mathbb{E}\bigg(|z_{ij}^{(r)}|^4\mathbbm{1}_{\{|z_{ij}^{(r)}|>\delta_n\sqrt{n}\}}\bigg) \longrightarrow 0.
\end{align}

In particular, $\{\delta_n\}$ can be chosen to satisfy $n\delta_n^4 \rightarrow \infty$. We thus have
\begin{align}\label{L4_condition}
    \sum_{i,j} a_{ij}^{(r)} = o(n^2\delta_n^4).
\end{align}

\textbf{Truncation:}
Recall that $z_{ij}^{(r)}$ is the ij-th entry of the matrix $Z_r$.  We now define the truncated versions of $Z_r$ as $\widehat{Z}_r = (\widehat{z}_{ij}^{(r)})$, where
$$\widehat{z}_{ij}^{(r)} = z_{ij}^{(r)}\mathbbm{1}_{\{|z_{ij}^{(r)}| \leq \delta_n\sqrt{n}\}}.$$
The truncated version of $S_n$ is thus
$$\widehat{S}_n = \frac{1}{n}\widehat{X}_n\widehat{X}_n^*, \quad \widehat{X}_n = \sum_{r=1}^KA_r^\frac{1}{2}\widehat{Z}_rB_r^\frac{1}{2}.$$
Let us denote the event where all innovation entries are below $\delta_n\sqrt{n}$ by 
\begin{align*}
    \mathcal{G}_n :=\bigcap_{i,j,r} \bigl\{|z_{ij}^{(r)}| \leq \delta_n\sqrt{n}\bigr\}.
\end{align*}
As a consequence, by a simple union bound and from (\ref{L4_condition}), we have
\begin{align*}
\mathbb{P}(\mathcal{G}_n^c) &\leq 
\sum_{i,j,r}\mathbb{P}\bigl(|z_{ij}^{(r)}| > \delta_n\sqrt{n}\bigr) \leq \frac{1}{\delta_n^4n^2}\sum_{i,j,r}\mathbb{E}|z_{ij}^{(r)}|^4\mathbbm{1}_{\{|z_{ij}^{(r)}|>\delta_n\sqrt{n}\}} =o(1).
\end{align*}
Let $\xi_n(\cdot)$ and $\widehat{\xi}_n(\cdot)$ represent the empirical processes arising from $S_n$ and $\widehat{S}_n$ respectively. On $\mathcal{G}_n$, we have $S_n = \widehat{S}_n$ and thus $\xi_n(\cdot) = \widehat{\xi}_n(\cdot)$. As a result, for analytic $f$, we get
$$\int f(\lambda)d\xi_n(\lambda) = \int f(\lambda)d\widehat{\xi}_n(\lambda) + o_P(1).$$
\textbf{Centering and Rescaling:}
Similar to $\widehat{Z}_r$ and $\widehat{S}_n$, we now define $\widetilde{Z}_r$ as $$\widetilde{Z}_r = (\widehat{Z}_r - \mathbb{E}\widehat{Z}_r) \circ \Gamma_r, \quad i \in [p], j \in [n], r \in [K]$$ where $\circ$ represents element-wise product and 
$$[\Gamma_r]_{ij} = 1/\sigma_{ij}^{(r)}, \quad (\sigma_{ij}^{(r)})^2 = \mathbb{E}|\widehat{z}_{ij}^{(r)} - \mathbb{E}\widehat{z}_{ij}^{(r)}|^2.$$
Then, we define $\widetilde{S}_n$ as
$$\widetilde{S}_n = \frac{1}{n}\widetilde{X}_n\widetilde{X}_n^*, \quad \widetilde{X}_n = \sum_{r=1}^KA_r^\frac{1}{2}\widetilde{Z}_rB_r^\frac{1}{2}.$$
From \cite{BaiYinKri}, for any $\epsilon>0$, we have
$$\mathbb{P}\bigg(\limsup_n \lambda_{\max}(\widehat{S}_n) > K\tau(1+\sqrt{c})^2+\epsilon, \quad \text{i.o.}\bigg) = 0$$
and the same bound holds also for 
$\limsup_n \lambda_{\max}(\widetilde{S}_n)$. Let $\widetilde{\xi}_n(\cdot)$ denote the analog of $\xi_n(\cdot)$ with $S_n$ replaced by $\widetilde{S}_n$. Let $C_j$ be a bound on $\|f_j^\prime\|$ in the relevant region. Using (a version of) the Hoffman-Wielandt inequality (Lemma \ref{lemma:Hoffman_Wielandt}), we get
\begin{align}\label{hat_tilde_integral_error}
    \bigg|\int f_j(\lambda)d\widehat{\xi}_n(\lambda) - \int f_j(\lambda)d\widetilde{\xi}(\lambda)\bigg| &= \absmod{\sum_{j=1}^pf_j(\lambda_i(\widehat{S}_n)) - f_j(\lambda_i(\widetilde{S}_n))}\\ \notag
    &\leq\, C_j \sum_{i=1}^p |\lambda_i(\widehat{S}_n) - \lambda_i(\widetilde{S}_n)|\\ \notag
    &\leq C_j\bigg(2\operatorname{tr}(\widehat{S}_n + \widetilde{S}_n)\bigg)^\frac{1}{2}\bigg(\frac{1}{n}\operatorname{tr}(\widehat{X}_n-\widetilde{X}_n)(\widehat{X}_n-\widetilde{X}_n)^*\bigg)^\frac{1}{2}\\ \notag
    &\leq \sqrt{2c_n}C_j\bigg(\lambda_{\max}(\widehat{S}_n) + \lambda_{\max}(\widetilde{S}_n)\bigg)^\frac{1}{2}\|\widehat{X}_n-\widetilde{X}_n\|_F.
\end{align}

Note that in order to show that 
$$\bigg|\int f_j(\lambda)d\widehat{\xi}_n(\lambda) - \int f_j(\lambda)d\widetilde{\xi}(\lambda)\bigg| \overset{p}{\longrightarrow} 0,$$
we will use Markov's inequality. To establish that the L.H.S. of (\ref{hat_tilde_integral_error}) converges to $0$ in expectation, we will use Lemma \ref{lemma:infinitely_often_rigorous}.
Since 
$$\lambda_{max}(\widehat{S}_n) \leq \operatorname{tr}(\widehat{S}_n) = O(n^2\delta_n^2),$$ 
and similarly for $\lambda_{max}(\widetilde{S}_n)$, it thus suffices to show that 
$$\mathbb{E}\|\widehat{X}_n-\widetilde{X}_n\|_F^2 \longrightarrow 0.$$

Let $M_r := \mathbbm{1}_p\mathbbm{1}_n^\top - \Gamma_r$. Then, we have
\begin{align*}
    \widehat{Z}_r - \widetilde{Z}_r &= \widehat{Z}_r - (\widehat{Z}_r - \mathbb{E}\widehat{Z}_r)\circ \Gamma_r = \widehat{Z}_r \circ M_r + \mathbb{E}\widehat{Z}_r\circ \Gamma_r\\
    \implies \widehat{X}_n - \widetilde{X}_n &= \sum_{r=1}^K A_r^\frac{1}{2}\bigl(\widehat{Z}_r \circ M_r + \mathbb{E}\widehat{Z}_r)\circ \Gamma_r\bigr)B_r^\frac{1}{2}
    \end{align*}

As a consequence, we have
\begin{align}\label{Xhat_Xtilde_frobenius}
\|\widehat{X}_n - \widetilde{X}_n\|_F &\leq  \sum_{r=1}^K\bigl\|A_r^\frac{1}{2} \bigl(\widehat{Z}_r \circ M_r\bigr)B_r^\frac{1}{2}\bigr\|_F + \sum_{r=1}^K\bigl\|A_r^\frac{1}{2} \bigl(\mathbb{E}\widehat{Z}_r\circ \Gamma_r)\bigr)B_r^\frac{1}{2}\bigr\|_F\\ \notag
    &\leq \tau \sum_{r=1}^K \bigg(\bigr\|\widehat{Z}_r \circ M_r\bigr\|_F + \bigr\|\mathbb{E}\widehat{Z}_r \circ \Gamma_r\bigr\|_F\bigg).
\end{align}

To bound the second term above, we need uniform bounds on $\mu_{ij}^{(r)} := \mathbb{E}\widehat{z}_{ij}^{(r)}$. Observe that
\begin{align*}
    |\mu_{ij}^{(r)}| \leq \frac{1}{\delta_n^3n^\frac{3}{2}}\mathbb{E}|z_{ij}^{(r)}|^4\mathbbm{1}_{\{|z_{ij}^{(r)}|>\delta_n\sqrt{n}\}} = \frac{a_{ij}^{(r)}}{\delta_n^3n^\frac{3}{2}}.
\end{align*}
Note that all for large $n$, we have the following inequality:
\begin{align}\label{crude_inequality}
    (a_{ij}^{(r)})^2 = \bigg(\mathbb{E}|z_{ij}^{(r)}|^4\mathbbm{1}_{\{|z_{ij}^{(r)}| > \delta_n\sqrt{n}\}}\bigg)^2 \leq M_4a_{ij}^{(r)},
\end{align}
we obtain
\begin{align}\label{l2_avergae_mu_ijr}
    \sum_{i,j}|\mu_{ij}^{(r)}|^2 \leq \frac{1}{\delta_n^6n^3}\sum_{i,j}\bigl(a_{ij}^{(r)}\bigr)^2 \leq \frac{M_4}{\delta_n^6n^3}\sum_{i,j} a_{ij}^{(r)} = o(\delta_n^{-2}n^{-1}) = o(n^{-\frac{1}{2}}).
\end{align}
where we used the fact that $n\delta_n^4 \rightarrow \infty$ implies $\delta_n^{-2} = o(\sqrt{n})$. To finish the argument for the second term of (\ref{Xhat_Xtilde_frobenius}), we need some uniform bounds on $\sigma_{ij}^{(r)}$. Starting from,
\begin{align*}
    (\sigma_{ij}^{(r)})^2 &= \mathbb{E}|\widehat{z}_{ij}^{(r)}|^2 - |\mathbb{E}\widehat{z}_{ij}^{(r)}|^2,
\end{align*}
we find that
\begin{align}\label{uniform_bound_sigma_ijr}
    |1 - (\sigma_{ij}^{(r)})^2| &\leq \mathbb{E}|z_{ij}^{(r)}|^2\mathbbm{1}_{\{|z_{ij}^{(r)}| > \delta_n\sqrt{n}\}} + |\mu_{ij}^{(r)}|^2 \\ \notag
    &\leq \frac{a_{ij}^{(r)}}{\delta_n^2n} + \frac{(a_{ij}^{(r)})^2}{\delta_n^6n^3}\\ \notag
    &\leq \frac{a_{ij}^{(r)}}{\delta_n^2n}\bigg(1 + \frac{M_4}{\delta_n^4n^2}\bigg) = O\bigg(\frac{a_{ij}^{(r)}}{\delta_n^2n}\bigg).
\end{align}
Thus, we have
\begin{align*}
    \max_{i,j,r} |1 - (\sigma_{ij}^{(r)})^2| \leq \max_{i,j,r}\frac{Ca_{ij}^{(r)}}{\delta_n^2n} \leq \frac{CM_4}{\delta_n^2n}\longrightarrow 0.
\end{align*}
Therefore, for all large $n$ we have $(\sigma_{ij}^{(r)})^2 \in [\frac{1}{2}, 2]$. Using (\ref{l2_avergae_mu_ijr}), the second term of (\ref{Xhat_Xtilde_frobenius}) is thus 
\begin{align*}
    \|\mathbb{E}\widehat{Z}_r\circ \Gamma_r\|_F^2 = \sum \frac{|\mu_{ij}^{(r)}|^2}{(\sigma_{ij}^{(r)})^2} \leq 2\sum_{i,j}|\mu_{ij}^{(r)}|^2 = o(n^{-\frac{1}{2}}).
\end{align*}

Now, we will control the first term of (\ref{Xhat_Xtilde_frobenius}). Observe that for $x \approx 1$, we have $|1 - 1/x| \leq C(1 - x^2)$. Thus,
$$[M_r]_{ij} = 1 - 1/\sigma_{ij} \leq C(1 - (\sigma_{ij}^{(r)})^2).$$ 
Using the fact that $$|\mu_{ij}^{(r)}|^2 = \mathbb{E}|\widehat{z}_{ij}^{(r)}|^2 = (\sigma_{ij}^{(r)})^2 + |\mu_{ij}^{(r)}|^2 = O(1),$$ we use (\ref{crude_inequality}) and (\ref{uniform_bound_sigma_ijr}) to get
\begin{align*}
    \mathbb{E}\bigg(\|M_r\circ \widehat{Z}_r\|_F^2\bigg) \leq M_4\sum_{i,j} \bigl(1 - (\sigma_{ij}^{(r)})^2\bigr)^2 \times |\mu_{ij}^{(r)}|^2 \leq \frac{C}{\delta_n^4n^2}\sum_{i,j}a_{ij}^{(r)} = o(1).
\end{align*}
As a consequence, we have 
\begin{align*}
    \mathbb{E}\|\widehat{X}_n - \widetilde{X}_n\|^2 \leq 4K^2\tau^2 \sum_{r=1}^K \bigg(\mathbb{E}\|M_r\circ \widehat{Z}_r\|_F^2 + \|\widehat{Z}_r\circ \Gamma_r\|_F^2\bigg) = o(1).
\end{align*}
This concludes the proof.

\begin{lemma}\label{lemma:infinitely_often_rigorous}
Suppose $\{X_n\}$ satisfies
\begin{enumerate}
    \item $X_n = O(n^k)$ for some $k$,
    \item $X_n \leq C_0$ a.s.,
    \item $\mathbb{P}(X_n > C_0 + \epsilon) = n^{-l}$ for any positive $l$.
\end{enumerate} 
and additionally, $\mathbb{E}|Y_n|^2 \to 0$. Then,
    $$\mathbb{E}X_nY_n \to 0.$$
\end{lemma}

\begin{proof}
For arbitrary $\epsilon>0$, we have
    $\mathbb{P}(X_n > C_0 + \epsilon) = n^{-l}$ for any positive $l$. Let $\mathcal{G} = \{X_n \leq C_0 + \epsilon\}$. Using $|X_n| \leq C_1n^k$ and Cauchy-Schwarz inequality, we have
    \begin{align*}
        \mathbb{E}|X_nY_n| \leq (C_0 + \epsilon)\mathbb{E}|Y_n|\mathbbm{1}_\mathcal{G} + C_1n^k\sqrt{\mathbb{E}|Y_n^2|} \sqrt{\mathbb{P}(\mathcal{G}^c)}.
    \end{align*}
    The result is immediate from the above inequality.
\end{proof}

\section{Extension to non-Gaussian innovations}\label{sec:NonGaussian}
\textbf{Note on notation:} Throughout this section and Appendix \ref{sec:appendixD}, $\{e_i\}_{i=1}^p$ and $\{\widetilde{e}_j\}_{j=1}^n$ represent canonical bases of $\mathbb{C}^p$ and $\mathbb{C}^n$ respectively.

Without loss of generality, we will assume that the $Z_r$ matrices are purely real, the entries are truncated at $\delta_n\sqrt{n}$, have zero mean and have unit variance. We start by constructing an analogous version of $S_n$ based on i.i.d. standard Gaussian innovations. Denote
\begin{align}
    T_n:=\frac{1}{n}Y_nY_n^\top, \quad Y_n :=\sum_{r=1}^K A_r^\frac{1}{2}W_rB_r^\frac{1}{2},
\end{align}
where the entries of $W_r=(w_{ij}^{(r)}) \in \mathbb{R}^{p\times n}$ are independent real Gaussian random variables such that
$$\mathbb{E} w_{ij}^{(r)}=0, \quad \mathbb{E} (w_{ij}^{(r)})^2=1, \quad i \in [p], j \in [n], r \in [K].$$
We also assume that $\{Z_r\}_{r=1}^K$ and $\{W_r\}_{r=1}^K$ are independent random matrices. Let $f(\lambda)$ be a smooth function with Fourier transform given by 
\begin{align}\label{fourier_transform}
    \widehat{f}(t) = \frac{1}{2\pi} \int_{-\infty}^\infty f(\lambda)e^{-\mathbbm{i}t\lambda}d\lambda.
\end{align}
Our goal for this section is to prove that the distributions of $\operatorname{tr}f(S_n)$ and  $\operatorname{tr}f(T_n)$ get arbitrarily close to each other for large $n$. To this end, we define the following matrices as functions of $\theta\in[0,\pi/2]$:
$$D_{n,r}(\theta):= Z_r\sin\theta + W_r\cos\theta,$$
and
\begin{align}\label{defining_En_Fn_theta}
    \mathbf{C}_n(\theta)=\frac{1}{n} E_n(\theta)E_n^{\top}(\theta), \quad E_n(\theta) = \sum_{r=1}^K A_r^\frac{1}{2}D_{n,r}(\theta)B_r^\frac{1}{2}.
\end{align}
\begin{remark}
    This particular interpolation scheme preserves the variance of the innovations. It also has a technical advantage in that many remainder terms vanish in (\ref{decomposing_del_Kn_theta_3}), which is partly due to the fact that $D_{n,r}(\theta)$ and its derivative with respect to $\theta$ are uncorrelated.
\end{remark}
We denote the $ij^{th}$ element of $D_{n,r}(\theta)$ by $d_{ijr}(\theta)$. Furthermore, let
\begin{align}\label{defining_In_Jn_Jn0_Kn_theta}
    &I_n(t,\theta)=e^{it\mathbf{C}_n(\theta)}, \quad J_n(\theta)=\operatorname{tr} f\bigl(\mathbf{C}_n(\theta)\bigr),\\ \notag
    &J_n^0(\theta)=J(\theta)-p\int f(x)\,dF^{c_n,H_{n},G_{n}}(x),\qquad
K_n(x,\theta)=\mathbb{E} e^{\mathbbm{i}xJ^0(\theta)}.
\end{align}
So, $J_n^0(\theta)$ is the new \textit{centered} empirical process corresponding to the sample covariance matrix $\mathbf{C}_n(\theta)$ and $K_n(\cdot, \theta)$ represents its characteristic function. Since $\mathbf{C}_n(\pi/2) = S_n$ and $\mathbf{C}_n(0)=T_n$, it is sufficient to show that 
$$|K_n(x, \pi/2) - K_n(x, 0)| \longrightarrow{} 0.$$
Observing that
\begin{align*}
    K_n(x, \pi/2) - K_n(x, 0) = \int_0^{\pi/2} \frac{\partial K_n(x, \theta)}{\partial \theta} d\theta,
\end{align*}
it therefore suffices to show that for any $x \in \mathbb{R}$, ${\partial K_n(x, \theta)}/{\partial \theta}$ uniformly converges to $0$ for $\theta \in [0, \pi/2]$. Along with $f(\cdot)$ and $\widehat{f}(\cdot)$, we introduce another related function $\widetilde{f}(\cdot)$ as follows:
\begin{align}\label{defining_f_tilde_Fn_theta}
    \widetilde{f}(\mathbf{C}_n) = \mathbbm{i}\int_{-\infty}^\infty u\widehat{f}(u)I_n(u)du.
\end{align}
For simplicity, we omit the argument $\theta$ from the notations of $D_{n,r}(\theta)$, $\mathbf{C}_n(\theta)$, $I_n(t,\theta)$ and denote them by $D_{n,r}$, $\mathbf{C}_n$, $I_n(t)$, respectively. 
 Denote $$\widetilde{d}_{ijr}(\theta) := \frac{\partial d_{ijr}(\theta)}{\partial\theta} = z_{ij}^{(r)}\cos\theta - w_{ij}^{(r)}\sin\theta.$$
We observe that $\widetilde{d}_{ijr}$ is uncorrelated with $d_{ijr}$. From Lemma \ref{lemma:derivative_Jn_theta}, we have 
\begin{align}\label{decomposing_del_Kn_theta}
    \frac{\partial K_n(x, \theta)}{\partial \theta} &=  \frac{\partial \mathbb{E}[e^{\mathbbm{i}xJ_n^0(\theta)}]}{\partial \theta}\\ \notag
    &= \mathbbm{i}x\,\mathbb{E} \bigg(e^{\mathbbm{i}xJ_n^0(\theta)}  \frac{\partial J_n(\theta)}{\partial \theta}\bigg)\\ \notag
    &= \mathbbm{i}x \sum_{i=1}^p\sum_{j=1}^n\sum_{r=1}^K\mathbb{E} \bigg(e^{\mathbbm{i}xJ_n^0(\theta)}  \frac{\partial J_n(\theta)}{\partial d_{ijr}}\frac{\partial d_{ijr}}{\partial \theta}\bigg)\\ \notag
    &= \frac{2\mathbbm{i}x}{n} \sum_{i=1}^p\sum_{j=1}^n\sum_{r=1}^K\mathbb{E}\, \bigg( \widetilde{d}_{ijr}[A_r^\frac{1}{2}\widetilde{f}(\mathbf{C}_n(\theta))E_nB_r^\frac{1}{2}]_{ij}e^{\mathbbm{i}xJ_n^0(\theta)}\bigg).
\end{align}
We will sequentially introduce a few more notations to simplify the expression within the expectation parentheses above. Note that $E_n(\theta)$ defined in (\ref{defining_En_Fn_theta}) is a sum of $K$ components. We introduce an analogous version where the $ij^{th}$ element of the $r^{th}$ component (i.e. $d_{ijr}(\theta)$) is replaced with $\psi$ as follows:
\begin{align}\label{defining_Enijr_omega}
    E_{n,ijr}(\psi, \theta) := \sum_{s\neq r}A_s^\frac{1}{2}D_{n,s}(\theta)B_s^\frac{1}{2} + A_r^\frac{1}{2}\bigg(D_{n,r}(\theta) + (\psi-d_{ijr}(\theta))e_i\widetilde{e}_j\bigg)B_r^\frac{1}{2}.
\end{align} 
Since we have fixed an arbitrary $\theta \in [0, \pi/2]$, we will henceforth ignore $\theta$ in the argument of $E_{n,ijr}$ and the following objects for simplicity.
Then analogous to (\ref{defining_En_Fn_theta}), we introduce the following: 
\begin{align}\label{defining_Fnijr_omega}
    \mathbf{C}_{n,ijr}(\psi)=\frac{1}{n} E_{n,ijr}(\psi)E_{n,ijr}^{\top}(\psi).
\end{align}
Furthermore, let
\begin{align}\label{defining_In_Jn_Jn0_Kn_omega_theta}
    &I_{n,ijr}(\psi, t)=e^{it\mathbf{C}_{n,ijr}(\psi)}, \quad J_{n,ijr}(\psi)=\operatorname{tr} f\bigl(\mathbf{C}_{n,ijr}(\psi)\bigr),\\ \notag
    &J_{n,ijr}^0(\psi)=J_{n,ijr}(\psi)-p\int f(x)\,dF^{c_n,H_{n},G_{n}}(x),\qquad
K_{n,ijr}(\psi, x)=\mathbb{E} e^{\mathbbm{i}xJ_{n}^0(\psi)}.
\end{align}
The analogous version of (\ref{defining_f_tilde_Fn_theta}) becomes:
\begin{align}\label{defining_f_tilde_Fn_omega_theta}
    \widetilde{f}(\mathbf{C}_{n,ijr}(\psi)) = \mathbbm{i}\int_{-\infty}^\infty u \widehat{f}(u) I_{n,ijr}(\psi, u) du.
\end{align}
Finally let
\begin{align}\label{defining_psi_ijr_omega}
    \psi_{ijr}(\psi) = [A_r^\frac{1}{2}\widetilde{f}(\mathbf{C}_{n,ijr}(\psi))E_nB_r^\frac{1}{2}]_{ij}e^{\mathbbm{i}xJ_{n}^0(\psi)}, \quad i \in [p], j \in [n], r \in [K].
\end{align}
Observe that plugging $d_{ijr}$ back into $\mathbf{C}_{n,ijr}(\psi)$ gives us the original objects, i.e.
$$\mathbf{C}_{n,ijr}(d_{ijr}) = \mathbf{C}_n, \quad \psi_{ijr}(d_{ijr}) = [A_r^\frac{1}{2}\widetilde{f}(\mathbf{C}_{n})E_nB_r^\frac{1}{2}]_{ij}e^{\mathbbm{i}xJ_{n}^0}.$$
Therefore, we can express (\ref{decomposing_del_Kn_theta}) as follows:
\begin{align}\label{decomposing_del_Kn_theta_1}
    \frac{\partial J_n(x, \theta)}{\partial \theta} = \frac{2\mathbbm{i}x}{n} \sum_{i=1}^p\sum_{j=1}^n\sum_{r=1}^K\mathbb{E}\, \bigl( \widetilde{d}_{ijr}\psi_{ijr}(d_{ijr})\bigr).
\end{align}
By Taylor expansion of $\psi_{ijr}(\cdot)$ at $d_{ijr}$, we have
\begin{align}\label{Taylor_expansion}
    \psi_{ijr}(d_{ijr}) = \sum_{l=0}^3 \frac{1}{l!} d_{ijr}^l\psi_{ijr}^{(l)}(0) + \frac{1}{4!}d_{ijr}^4 \psi_{ijr}^{(4)}(\rho_{ijr}\, d_{ijr}), \quad \rho_{ijr} \in (0,1).
\end{align}
Plugging this into (\ref{decomposing_del_Kn_theta_1}) and using the independence of $\bigl( d_{ijr}, \widetilde{d}_{ijr}\bigr)$ and $\psi_{ijr}(0)$, we get
\begin{align}\label{decomposing_del_Kn_theta_2}
     \frac{\partial J_n(x, \theta)}{\partial \theta} = \frac{2\mathbbm{i}x}{n} \sum_{l=0}^3 \frac{1}{l!}\sum_{i,j,r}\mathbb{E} \bigl(\widetilde{d}_{ijr}d_{ijr}^l\bigr) \, \mathbb{E}\,(\psi_{ijr}^{(l)}(0)) + \frac{2\mathbbm{i}x}{4!n}\sum_{i,j,r} \mathbb{E}\bigl(\widetilde{d}_{ijr}(\theta) d_{ijr}^4 \psi_{ijr}^{(4)}(\rho_{ijr}d_{ijr})\bigr).
\end{align}

By using the properties of $z_{ij}^{(r)}$ and $w_{ij}^{(r)}$, we have the following results:
$$\mathbb{E}\widetilde{d}_{ijr}d_{ijr}^0 = 0, \quad \mathbb{E}\widetilde{d}_{ijr}d_{ijr}^1 = 0,
\mathbb{E}\widetilde{d}_{ijr}d_{ijr}^2= \mathbb{E}d_{ijr}^3\sin^2\theta\cos\theta, \mathbb{E}\widetilde{d}_{ijr} d_{ijr}^3 = o(1)\sin^3\theta\cos\theta.$$

Therefore, (\ref{decomposing_del_Kn_theta_2}) reduces to
\begin{align}\label{decomposing_del_Kn_theta_3}
    \frac{\partial J_n(x, \theta)}{\partial \theta} &= \frac{\mathbbm{i}x}{n} \sum_{i,j,r}\mathbb{E}d_{ijr}^3\sin^2\theta\cos\theta\,\mathbb{E}\psi_{ijr}^{(2)}(0)+ \frac{\mathbbm{i}x}{12n} \mathbb{E}\bigl(\widetilde{d}_{ijr}(\theta) d_{ijr}^4 \psi_{ijr}^{(4)}(\rho_{ijr}d_{ijr})\bigr)\\ \notag
    &= \mathcal{I}_1 + \mathcal{I}_2.
\end{align}

\subsection{Analysis of \texorpdfstring{$\psi_{ijr}^{(2)}(0)$}{} and \texorpdfstring{$\mathcal{I}_1$}{}}
We start by observing that a detailed analysis of $\mathcal{I}_1$ requires a control on $\psi_{ijr}^{(2)}(0)$ since $\mathbb{E}d_{ijr}^3$ are uniformly bounded. We fix $i \in [p],j\in [n],r\in [K]$ arbitrarily and then analyze $\psi_{ijr}^{(2)}(0)$. Note that
\begin{align}\label{defining_h_omega}
    \psi_{ijr}(\psi)=h(\psi)e^{\mathbbm{i}xJ_{n,ijr}^0(\psi)}, \text{ where  } h(\psi) := [A_r^\frac{1}{2}\widetilde{f}(\mathbf{C}_{n,ijr}(\psi))E_nB_r^\frac{1}{2}]_{ij}.
\end{align}
So, we need the derivatives of $h(\psi)$ and $J_{n,ijr}(\psi)$ at $0$. Going forward, the $k^{th}$ derivatives of $\widetilde{f}(\mathbf{C}_{n,ijr}(\psi))$, $I_{n,ijr}(\psi, t)$ and $J_{n,ijr}(\psi)$ shall be denoted as:
\begin{align*}
  &\nabla_\psi^{(k)}\widetilde{f}(\mathbf{C}_{n,ijr}(\psi_0)):= \frac{\partial^k \widetilde{f}(\mathbf{C}_{n,ijr}(\psi))}{\partial \psi^k}\bigg|_{\psi=\psi_0},\\  
  &\nabla_\psi^{(k)}I_{n,ijr}(\psi_0, t)):= \frac{\partial^k I_{n,ijr}(\psi, t))}{\partial \psi^k}\bigg|_{\psi=\psi_0},\\
  &\nabla_\psi^{(k)}J_{n,ijr}(\psi_0)):= \frac{\partial^k J_{n,ijr}(\psi))}{\partial \psi^k}\bigg|_{\psi=\psi_0}.
\end{align*}
\textbf{Notation:} Since the indices $i,j,r$ are fixed throughout this subsection, we will write $E_{n,ijr}(\psi)$, $\mathbf{C}_{n,ijr}(\psi)$, $I_{n,ijr}(\psi, t)$, $J_{n,ijr}(\psi)$ and $J_{n,ijr}^0(\psi)$ as $E_n(\psi)$, $\mathbf{C}_n(\psi)$, $I_n(\psi, t)$, $J_n(\psi)$ and $J_n^0(\psi)$ respectively.

From (\ref{defining_h_omega}), we have
\begin{align*}
    \psi_{ijr}^{(1)}(\psi) &= h^{(1)}(\psi)e^{\mathbbm{i}xJ_{n}^0(\psi)} + \mathbbm{i}xh(\psi)e^{\mathbbm{i}xJ_{n}^0(\psi)}\nabla_\psi^{(1)} J_{n}(\psi)\\
    &= h^{(1)}(\psi)e^{\mathbbm{i}xJ_{n}^0(\psi)} + \frac{2\mathbbm{i}x}{n}h^2(\psi)e^{\mathbbm{i}xJ_{n}^0(\psi)},
    \end{align*}
where the last equality follows from the fact that by Lemma \ref{lemma:derivative_Jn_theta}, we have $$\nabla_\psi^{(1)} J_n(0) = \frac{2h(0)}{n}, \quad \nabla_\psi^{(k)} J_n(0) = \frac{2h^{(k-1)}(0)}{n}.$$
As a result, we have
\begin{align}\label{psi_ijr_2nd_derivative}
    \psi_{ijr}^{(2)}(\psi)
    &= e^{\mathbbm{i}xJ_{n}^0(\psi)}\bigg(h^{(2)}(\psi) + \mathbbm{i}xh^{(1)}(\psi)\nabla_\psi^{(1)} J_n(\psi)\\ \notag
    &+ \frac{2\mathbbm{i}x}{n}\bigl(2h(\psi)h^{(1)}(\psi) +  \mathbbm{i}xh^2(\psi)\nabla_\psi^{(1)} J_n(\psi)\bigr)\bigg)\\ \notag
    \implies \psi_{ijr}^{(2)}(0) &= e^{\mathbbm{i}xJ_n^0(0)}\bigg(h^{(2)}(0) + \frac{6\mathbbm{i}x}{n}h(0)h^{(1)}(0) - \frac{4x^2}{n^2}h^3(0)\bigg).
\end{align}
Now, all that we need are the derivatives of $h(\psi)$ at $0$. Upon direct differentiation, we obtain the following:
\begin{align}\label{derivative_h_d_ijr}
    \implies h^{(1)}(0) &= \bigg[A_r^\frac{1}{2}\nabla_\psi^{(1)}\widetilde{f}(\mathbf{C}_n(0))E_n(0)B_r^\frac{1}{2}\bigg]_{ij} + \bigg[A_r^\frac{1}{2}\widetilde{f}(\mathbf{C}_n(0))A_r^\frac{1}{2}e_i\tilde{e}_j^\top B_r^\frac{1}{2}B_r^\frac{1}{2}\bigg]_{ij}\\ \notag
     &= \bigg[A_r^\frac{1}{2}\nabla_\psi^{(1)}\widetilde{f}(\mathbf{C}_n(0))E_n(0)B_r^\frac{1}{2}\bigg]_{ij} + [A_r^\frac{1}{2}\widetilde{f}(\mathbf{C}_n(0))A_r^\frac{1}{2}]_{ii} [B_r]_{jj}\\ \notag
     \implies h^{(2)}(0) &=  \bigg[A_r^\frac{1}{2}\nabla_\psi^{(2)}\widetilde{f}(\mathbf{C}_n(0))E_n(0)B_r^\frac{1}{2}\bigg]_{ij} + 2\bigg[A_r^\frac{1}{2}\nabla_\psi^{(1)}\widetilde{f}(\mathbf{C}_n(0))A_r^\frac{1}{2}\bigg]_{ii} [B_r]_{jj}.
\end{align}
Combining everything into (\ref{psi_ijr_2nd_derivative}), we get
    \begin{align}
    \psi_{ijr}^{(2)}(0)
    &=e^{\mathbbm{i}xJ_n^0(0)}\bigg(\bigg[A_r^\frac{1}{2}\nabla_\psi^{(2)}\widetilde{f}(\mathbf{C}_n(0))E_n(0)B_r^\frac{1}{2}\bigg]_{ij}\\ \notag
    &+ 2\bigg[A_r^\frac{1}{2}\nabla_\psi^{(1)}\widetilde{f}(\mathbf{C}_n(0))A_r^\frac{1}{2}\bigg]_{ii} [B_r]_{jj} \\ \notag
    &+ \frac{6\mathbbm{i}x}{n}[A_r^\frac{1}{2}\widetilde{f}(\mathbf{C}_n(0))E_n(0)B_r^\frac{1}{2}]_{ij}\bigg[A_r^\frac{1}{2}\nabla_\psi^{(1)}\widetilde{f}(\mathbf{C}_n(0))E_n(0)B_r^\frac{1}{2}\bigg]_{ij}\\ \notag
    &+ \frac{6\mathbbm{i}x}{n}[A_r^\frac{1}{2}\widetilde{f}(\mathbf{C}_n(0))E_n(0)B_r^\frac{1}{2}]_{ij}\bigg[A_r^\frac{1}{2}\nabla_\psi^{(1)}\widetilde{f}(\mathbf{C}_n(0))A_r^\frac{1}{2}\bigg]_{ii} [B_r]_{jj} \\ \notag
    &- \frac{4x^2}{n^2}[A_r^\frac{1}{2}\widetilde{f}(\mathbf{C}_n(0))E_n(0)B_r^\frac{1}{2}]_{ij}^3\bigg)\\ \notag
    &:= \mathcal{I}_{ijr}^{(1)}+\mathcal{I}_{ijr}^{(2)}+\mathcal{I}_{ijr}^{(3)}+\mathcal{I}_{ijr}^{(4)}+\mathcal{I}_{ijr}^{(5)}.
\end{align}

\subsubsection{Analysis of \texorpdfstring{$\mathcal{I}_{ijr}^{(l)}\}_{l=1}^5$}{}}
We will handle only $\mathcal{I}_{ijr}^{(1)}$ since the others can be similarly derived.  From (\ref{defining_f_tilde_Fn_theta}), we have
\begin{align} \label{decimposing_I_ijr_1}
\mathcal{I}_{ijr}^{(1)}&= \bigg[A_r^\frac{1}{2}\nabla_\psi^{(2)}\widetilde{f}(\textbf{C}_n(0))E_n(0)B_r^\frac{1}{2}\bigg]_{ij}e^{\mathbbm{i}xJ_n^0(0)}\\ \notag
     &= \mathbbm{i}\bigg[A_r^\frac{1}{2}\bigg(\int_{-\infty}^\infty u\widehat{f}(u)\nabla_\psi^{(2)} I_n(0, u)du\bigg)E_n(0)B_r^
     \frac{1}{2}\bigg]_{ij}e^{\mathbbm{i}xJ_n^0(0)}.
\end{align}
Using \textit{Duhamel's formula} (Lemma \ref{lemma:Duhamel}), we simplify $\nabla_\psi^{(2)} I_n(0, u)$ as follows:
\begin{align*}
    \nabla_\psi^{(1)} I_n(0, u) &= \mathbbm{i}\int_0^u I_n(0, v) \nabla_\psi^{(1)} \mathbf{C}_n(0) I_n(0, u-v)\,dv
    \end{align*}
Differentiating again, we get    
\begin{align*}
&\nabla_\psi^{(2)} I_n(0, u)\\
=& \mathbbm{i}\int_0^u \nabla_\psi^{(1)} I_n(0,v) \nabla_\psi^{(1)} \mathbf{C}_n(0) I_n(0,u-v)\,dv \\
&+ \mathbbm{i}\int_0^u I_n(0,v)\nabla_\psi^{(2)} \mathbf{C}_n(0)I_n(0, u-v)\,dv \\ 
&+ \mathbbm{i}\int_0^u I_n(0,v) \nabla_\psi^{(1)} \mathbf{C}_n(0) \nabla_\psi^{(1)} I_n(0,u-v)\,dv\\
&= \mathbbm{i}\int_0^u I_n(0,v)\nabla_\psi^{(2)} \mathbf{C}_n(0)I_n(0, u-v)\,dv \\
&-\int_0^u\int_0^v I_n(0,w)\nabla_\psi^{(1)} \mathbf{C}_n(0) I_n(0,v-w)\nabla_\psi^{(1)} \mathbf{C}_n(0) I_n(0,u-v)\,dw\,dv \\
& -\int_0^u\int_0^{u-v} I_n(0,v)\nabla_\psi^{(1)}\mathbf{C}_n(0)I_n(0,w)\nabla_\psi \mathbf{C}_n(0)I_n(0,u-v-w)\,dw\,dv\\
&:= \mathcal{T}_{ijr}^{(1)}+2\mathcal{T}_{ijr}^{(2)},
\end{align*}
where we note that the last and the second last terms in the expansion of the second partial derivative are equal. Therefore (\ref{decimposing_I_ijr_1}) splits as follows:
\begin{align}\label{decimposing_I_ijr_2}
    \mathcal{I}_{ijr}^{(1)} = \mathbbm{i}\int_{-\infty}^\infty u\widehat{f}(u)\bigg[A_r^\frac{1}{2}(\mathcal{T}_{ijr}^{(1)} + 2\mathcal{T}_{ijr}^{(2)})E_n(0)B_r^\frac{1}{2}]_{ij}e^{\mathbbm{i}xJ_n^0(0)}.
\end{align}
The first term of (\ref{decimposing_I_ijr_2}) can be expressed as follows:
\begin{align}\label{first_term_I_ijr_1}
    &\mathbbm{i}\int_{-\infty}^\infty u\widehat{f}(u)\bigg[A_r^\frac{1}{2}\mathcal{T}_{ijr}^{(1)}E_n(0)B_r^\frac{1}{2}\bigg]_{ij}e^{\mathbbm{i}xJ_n^0(0)}du\\ \notag
    =& - \int_{-\infty}^\infty u\widehat{f}(u) \bigg[A_r^\frac{1}{2}\bigg(\int_0^u I_n(0,v)\nabla_\psi^{(2)} \mathbf{C}_n(0)I_n(0,u-v)\,dv\bigg) E_n(0)B_r^\frac{1}{2}\bigg]_{ij}e^{\mathbbm{i}xJ_n^0(0)}du\\ \notag
    =& -\frac{2[B_r]_{jj}}{n} \int_{-\infty}^\infty u\widehat{f}(u) \int_0^u  [A_r^\frac{1}{2}I_n(0,v)A_r^\frac{1}{2}]_{ii}[A_r^\frac{1}{2}I_n(0,u-v)E_n(0)B_r^\frac{1}{2}]_{ij}e^{\mathbbm{i}xJ_n^0(0)}\,dv\,du\\ \notag
    =& -\frac{2[B_r]_{jj}}{n} \int_{-\infty}^\infty u\widehat{f}(u) [A_r^\frac{1}{2}I_n(0)A_r^\frac{1}{2}]_{ii}*[A_r^\frac{1}{2}I_n(0)E_n(0)B_r^\frac{1}{2}]_{ij}(u) e^{\mathbbm{i}xJ_n^0(0)}\,du,
\end{align}
where we used the fact that by Lemma \ref{lemma:derivative_Fn_theta}, we have
\begin{align*}
    \nabla_\psi^{(2)} \mathbf{C}_n(\psi) = \frac{2}{n}A_r^\frac{1}{2}e_i\tilde{e}_j^\top B_r\tilde{e}_je_i^\top A_r^\frac{1}{2},
\end{align*}
and the convolution $(*)$ of two functions is defined as $f*g(x) = \in_0^tf(s)g(t-s)ds.$ We now simplify the second term of (\ref{decimposing_I_ijr_2}), i.e.
$$2\mathbbm{i}\int_{-\infty}^\infty u\widehat{f}(u)\bigg[A_r^\frac{1}{2}\mathcal{T}_{ijr}^{(2)}E_n(0)B_r^\frac{1}{2}\bigg]_{ij}e^{\mathbbm{i}xJ_n^0(0)}\,du.$$
From Lemma \ref{lemma:derivative_Fn_theta}, we have
\begin{align*}
    \nabla_\psi^{(1)} \textbf{C}_n(0)= \frac{1}{n}(M_{ijr} + M_{ijr}^\top), \quad \text{ where, } M_{ijr} = A_r^\frac{1}{2}e_i\tilde{e}_j^\top B_r^\frac{1}{2}E_n^\top.
\end{align*}
Then, the $ij^{th}$ term of the involved matrix in the second term of (\ref{decimposing_I_ijr_2}) simplifies as follows:
\begin{align*}
    &-\bigg[A_r^\frac{1}{2}\mathcal{T}_{ijr}^{(2)}E_nB_r^\frac{1}{2}\bigg]_{ij}\\
    =& \bigg[A_r^\frac{1}{2} \bigg(\int_0^u\int_0^v I_n(w)\nabla_\psi^{(1)} \mathbf{C}_n(0) I_n(0,v-w)\nabla_\psi^{(1)} \mathbf{C}_n(0) I_n(0, u-v)\,dw\,dv\bigg)E_n(0)B_r^\frac{1}{2}\bigg]_{ij}\\
    =& \frac{1}{n^2}\int_0^u \int_0^v \bigg[A_r^\frac{1}{2}I_n(0,w)(M_{ijr} + M_{ijr}^\top)I_n(0,v-w)(M_{ijr} + M_{ijr}^\top)I_n(0,u-v)E_n(0)B_r^\frac{1}{2}\bigg]_{ij}\,dw\,dv\\
    =& \frac{1}{n^2}\int_0^u \int_0^v \bigg[A_r^\frac{1}{2}I_n(0,w)M_{ijr}I_n(0,v-w)M_{ijr}I_n(0,u-v)E_n(0)B_r^\frac{1}{2}\bigg]_{ij}\,dw\,dv\\
    &+ \frac{1}{n^2}\int_0^u \int_0^v \bigg[A_r^\frac{1}{2}I_n(0,w)M_{ijr}I_n(0,v-w)M_{ijr}^\top I_n(0,u-v)E_n(0)B_r^\frac{1}{2}\bigg]_{ij}\,dw\,dv\\
    &+\frac{1}{n^2}\int_0^u \int_0^v \bigg[A_r^\frac{1}{2}I_n(0,w)M_{ijr}^\top I_n(0,v-w)M_{ijr}I_n(0,u-v)E_n(0)B_r^\frac{1}{2}\bigg]_{ij}\,dw\,dv\\
    &+\frac{1}{n^2}\int_0^u \int_0^v \bigg[A_r^\frac{1}{2}I_n(0,w)M_{ijr}^\top I_n(0,v-w)M_{ijr}^\top I_n(0,u-v)E_n(0)B_r^\frac{1}{2}\bigg]_{ij}\,dw\,dv\\
    =& \frac{3}{n^2}\int_0^u \int_0^v [A_r^\frac{1}{2}I_n(0,w)Ar^\frac{1}{2}]_{ii}[A_r^\frac{1}{2}I_n(0,v-w)E_n(0)B_r^\frac{1}{2}]_{ij}[B_r^\frac{1}{2}E_n^\top I_n(0,u-v)E_n(0)B_r^\frac{1}{2}]_{jj}\,dw\,dv\\
    &+ \frac{1}{n^2}\int_0^u \int_0^v [A_r^\frac{1}{2}I_n(0,w)E_n(0)B_r^\frac{1}{2}]_{ij}[A_r^\frac{1}{2}I_n(0,v-w)E_n(0)B_r^\frac{1}{2}]_{ij}[A_r^\frac{1}{2}I_n(0,u-v)E_n(0)B_r^\frac{1}{2}]_{ij}\,dw\,dv\\
    =& \frac{3}{n^2}[A_r^\frac{1}{2}I_n(0)Ar^\frac{1}{2}]_{ii}*[A_r^\frac{1}{2}I_n(0)E_n(0)B_r^\frac{1}{2}]_{ij}*[B_r^\frac{1}{2}E_n^\top(0) I_n(0)E_n(0)B_r^\frac{1}{2}]_{jj}(u)\\ &+ \frac{1}{n^2}[A_r^\frac{1}{2}I_n(0)E_n(0)B_r^\frac{1}{2}]_{ij}*[A_r^\frac{1}{2}I_n(0)E_n(0)B_r^\frac{1}{2}]_{ij}*[A_r^\frac{1}{2}I_n(0)E_n(0)B_r^\frac{1}{2}]_{ij}(u).
\end{align*}
Therefore, the contribution of the second term of (\ref{decimposing_I_ijr_2}) turns out to be:
\begin{align}\label{second_term_I_ijr_1}
    &2\mathbbm{i}\int_{-\infty}^\infty u\widehat{f}(u)\bigg[A_r^\frac{1}{2}\mathcal{T}_{ijr}^{(2)}E_n(0)B_r^\frac{1}{2}\bigg]_{ij}(u)e^{\mathbbm{i}xJ_n^0(0)}du\\ \notag
    =& -\frac{6\mathbbm{i}}{n^2} \int_{-\infty} ^\infty u\widehat{f}(u) [A_r^\frac{1}{2}I_n(0)A_r^\frac{1}{2}]_{ii}*[A_r^\frac{1}{2}I_n(0)E_n(0)B_r^\frac{1}{2}]_{ij}*[B_r^\frac{1}{2}E_n^\top(0) I_n(0) E_n(0)B_r^\frac{1}{2}]_{jj}(u)e^{\mathbbm{i}xJ_n^0(0)}\,du \\ \notag
    &-\frac{2\mathbbm{i}}{n^2} \int_{-\infty} ^\infty u\widehat{f}(u) [A_r^\frac{1}{2}I_n(0)E_n(0)B_r^\frac{1}{2}]_{ij}*[A_r^\frac{1}{2}I_n(0)E_n(0)B_r^\frac{1}{2}]_{ij}* [A_r^\frac{1}{2}I_n(0)E_n(0)B_r^\frac{1}{2}]_{ij}(u)e^{\mathbbm{i}xJ_n^0(0)}\,du.
\end{align}
By the justification given in sub-section (\ref{subsec:boundedSpectralNorm}), we find that the quantities below 
$$\|A_r^\frac{1}{2}I_n(0,u)A_r^\frac{1}{2}\|, \quad \|n^{-\frac{1}{2}}A_r^\frac{1}{2}I_n(0,u)E_n(0)B_r^\frac{1}{2}\|, \quad \|n^{-1}B_r^\frac{1}{2}E_n^\top(0) I_n(0, u) E_n(0)B_r^\frac{1}{2}\|$$ have bounded moments. Therefore, applying Lemma \ref{lemma:5_by_4} gives us
\begin{align*}
    \absmod{\frac{1}{n}\sum_{i,j,r} \mathbb{E}\mathcal{I}_{ijr}^{(1)}} \leq \frac{C}{n^{1/4}} \int_{-\infty}^\infty (|u|^2 + |u|^3)\, |\widehat{f}(u)| \, du \leq \frac{C}{n^{1/4}}.
\end{align*}
Note that we can justify $\widehat{f^{(3)}}, \widehat{f^{(2)}}\in L^1(\mathbb{R})$ as follows. We can swap f with $\widetilde{f}$ where $\widetilde{f} \in \mathcal{C}_c(\mathbb{R})$ since $\operatorname{tr}f(\textbf{C}_n)$ and $\operatorname{tr}\widetilde{f}(\textbf{C}_n)$ agree everywhere except on a super-polynomially small event as in \cite{NajimYao}. In that case $\int |u|^m\widehat{f}(u)du < \infty$ for any $m$. By similar arguments, we also get
\begin{align*}
    \absmod{\frac{1}{n}\sum_{i,j,r} \mathbb{E}\mathcal{I}_{ijr}^{(2)}+\mathcal{I}_{ijr}^{(3)}+\mathcal{I}_{ijr}^{(4)}+\mathcal{I}_{ijr}^{(5)}} \leq \frac{C}{n^{1/4}}.
\end{align*}
Therefore, $|\mathcal{I}_1| \longrightarrow 0$ as $n \rightarrow \infty$.
\subsection{Analysis of \texorpdfstring{$\mathcal{I}_2$}{}}
Recall $\mathcal{I}_2$ from (\ref{decomposing_del_Kn_theta_2}). Following the arguments in Section C of supplementary file of \cite{BaiLiPan2019}, we will show that 
$$|\mathcal{I}_2| \leq C\delta_n.$$
As in the analysis of $\mathcal{I}_1$, it is sufficient to control the term $\psi_{ijr}^{(4)}(\psi_0)$ where we set 
$$\psi_0 = \rho_{ijr}d_{ijr}.$$
Let $a = \mathbbm{i}x$. We have 
\begin{align}\label{psi_dijr_4th_derivative}
    \psi_{ijr}^{(4)}(\psi) = \sum_{k=0}^4 h^{(4-k)}(\psi) \frac{\partial^k}{\partial \psi^4} e^{aJ_{n,ijr}^0(\psi)}
\end{align}
\textbf{Notation:} Similar to the analysis of $\mathcal{I}_1$,  we will write $E_{n,ijr}(\psi)$, $\mathbf{C}_{n,ijr}(\psi)$, $I_{n,ijr}(\psi, t)$, $J_{n,ijr}(\psi)$ and $J_{n,ijr}^0(\psi)$ as $E_n(\psi)$, $\mathbf{C}_n(\psi)$, $I_n(\psi, t)$, $J_n(\psi)$ and $J_n^0(\psi)$ respectively.

The derivatives of the exponential term are
\begin{align*}
    &\nabla_\psi^{(1)}e^{aJ_n^0(\psi)} = a\nabla_\psi^{(1)}J_n(\psi) e^{aJ_n^0(\psi)}, \quad \nabla_\psi^{(2)}e^{aJ_n^0(\psi)} = \bigl(a\nabla_\psi^{(2)}J_n(\psi)+a^2(\nabla_\psi^{(1)}J_n(\psi))^2\bigr)e^{aJ_n^0(\psi)},\\
    &\nabla_\psi^{(3)}e^{aJ_n^0(\psi)} = \bigl(a\nabla_\psi^{(3)}J_n(\psi) + 3a^2\nabla_\psi^{(1)}J_n(\psi) \nabla_\psi^{(2)}J_n(\psi) + a^3(\nabla_\psi^{(1)}J_n(\psi))^3\bigr)e^{aJ_n^0(\psi)},\\
    &\nabla_\psi^{(4)}e^{aJ_n^0(\psi)} = \bigl(a\nabla_\psi^{(4)}J_n(\psi) + 4a^2\nabla_\psi^{(1)}J_n(\psi) \nabla_\psi^{(3)}J_n(\psi) + 6a^3(\nabla_\psi^{(1)}J_n(\psi))^2\nabla_\psi^{(2)}J_n(\psi) \\
    &\hspace{2cm} + 3a^2(\nabla_\psi^{(2)}J_n(\psi))^2 + a^4(\nabla_\psi^{(1)}J_n(\psi))^4\bigr)e^{aJ_n^0(\psi)}.
\end{align*}
By Lemma \ref{lemma:derivative_Jn_theta}, we have 
\begin{align*}
    \nabla_\psi^{(k)}J_n(\psi_0) = 2h^{(k-1)}(\psi_0)/n.
\end{align*}
The other necessary ingredients of (\ref{psi_dijr_4th_derivative}) are the derivatives of $h(\cdot)$. $h^{(1)}, h^{(2)}$ have already been derived in (\ref{derivative_h_d_ijr}). Differentiating (\ref{derivative_h_d_ijr}) repeatedly, we obtain
\begin{align*}
    h^{(3)}(\psi_0)&= \bigg[A_r^\frac{1}{2}\nabla_\psi^{(3)} \widetilde{f}(\textbf{C}_n(\psi_0))E_n(\psi_0)B_r^\frac{1}{2}\bigg]_{ij} + 3\bigg[A_r^\frac{1}{2}\nabla_\psi^{(2)}\widetilde{f}(\textbf{C}_n(\psi_0))A_r^\frac{1}{2}\bigg]_{ii} [B_r]_{jj},\\
    h^{(4)}(\psi_0) &= \bigg[A_r^\frac{1}{2}\nabla_\psi^{(4)} \widetilde{f}(\textbf{C}_n(\psi_0))E_n(\psi_0)B_r^\frac{1}{2}\bigg]_{ij} + 4\bigg[A_r^\frac{1}{2}\nabla_\psi^{(3)}\widetilde{f}(\textbf{C}_n(\psi_0))A_r^\frac{1}{2}\bigg]_{ii} [B_r]_{jj}.
\end{align*}
Combining everything, (\ref{psi_dijr_4th_derivative}) expands into the following expression:
\begin{align}
    \psi_{ijr}^{(4)}(\psi_0) &= e^{\mathbbm{i}xJ_n^0(\psi_0)}\bigg(\bigg[A_r^\frac{1}{2}\nabla_\psi^{(4)} \widetilde{f}(\textbf{C}_n(\psi_0))E_n(\psi_0)B_r^\frac{1}{2}\bigg]_{ij}\\\notag
    &+ 4\bigg[A_r^\frac{1}{2}\nabla_\psi^{(3)} \widetilde{f}(\textbf{C}_n(\psi_0))A_r^\frac{1}{2}\bigg]_{ii}[B_r]_{jj}\\ \notag
    &+ \frac{10\mathbbm{i}x}{n}\bigg[A_r^\frac{1}{2}\nabla_\psi^{(3)} \widetilde{f}(\textbf{C}_n(\psi_0))E_n(\psi_0)B_r^\frac{1}{2}\bigg]_{ij}\bigg[A_r^\frac{1}{2}\widetilde{f}(\textbf{C}_n(\psi_0))E_n(\psi_0)B_r^\frac{1}{2}\bigg]_{ij} \\ \notag
    &+ \frac{30\mathbbm{i}x}{n}\bigg[A_r^\frac{1}{2}\widetilde{f}(\textbf{C}_n(\psi_0))E_n(\psi_0)B_r^\frac{1}{2}\bigg]_{ij}\bigg[A_r^\frac{1}{2}\nabla_\psi^{(2)} \widetilde{f}(\textbf{C}_n(\psi_0))A_r^\frac{1}{2}\bigg]_{ii}[B_r]_{jj}\\ \notag
    &+ \frac{20\mathbbm{i}x}{n}\bigg[A_r^\frac{1}{2}\nabla_\psi^{(2)} \widetilde{f}(\textbf{C}_n(\psi_0))E_n(\psi_0)B_r^\frac{1}{2}\bigg]_{ij}\bigg[A_r^\frac{1}{2}\nabla_\psi^{(1)} \widetilde{f}(\textbf{C}_n(\psi_0))E_n(\psi_0)B_r^\frac{1}{2}\bigg]_{ij}\\ \notag
    &+ \frac{20\mathbbm{i}x}{n}\bigg[A_r^\frac{1}{2}\nabla_\psi^{(2)} \widetilde{f}(\textbf{C}_n(\psi_0))E_n(\psi_0)B_r^\frac{1}{2}\bigg]_{ij}\bigg[A_r^\frac{1}{2}\widetilde{f}(\textbf{C}_n(\psi_0))A_r^\frac{1}{2}\bigg]_{ii}[B_r]_{jj}\\ \notag
    &+\frac{40\mathbbm{i}x}{n}\bigg[A_r^\frac{1}{2}\nabla_\psi^{(1)} \widetilde{f}(\textbf{C}_n(\psi_0))E_n(\psi_0)B_r^\frac{1}{2}\bigg]_{ij}\bigg[A_r^\frac{1}{2}\nabla_\psi^{(1)} \widetilde{f}(\textbf{C}_n(\psi_0))A_r^\frac{1}{2}\bigg]_{ii}[B_r]_{jj}\\ \notag
    &+\frac{40\mathbbm{i}x}{n}\bigg[A_r^\frac{1}{2}\nabla_\psi^{(1)} \widetilde{f}(\textbf{C}_n(\psi_0))A_r^\frac{1}{2}\bigg]_{ii}\bigg[A_r^\frac{1}{2}\widetilde{f}(\textbf{C}_n(\psi_0))A_r^\frac{1}{2}\bigg]_{ii}[B_r]_{jj}^2\\ \notag
    &-\frac{40x^2}{n^2}\bigg[A_r^\frac{1}{2}\nabla_\psi^{(2)} \widetilde{f}(\textbf{C}_n(\psi_0))E_n(\psi_0)B_r^\frac{1}{2}\bigg]_{ij}\bigg[A_r^\frac{1}{2}\widetilde{f}(\textbf{C}_n(\psi_0))E_n(\psi_0)B_r^\frac{1}{2}\bigg]_{ij}^2\\ \notag
    &-\frac{60x^2}{n^2}\bigg[A_r^\frac{1}{2}\widetilde{f}(\textbf{C}_n(\psi_0))E_n(\psi_0)B_r^\frac{1}{2}\bigg]_{ij}\bigg[A_r^\frac{1}{2}\nabla_\psi^{(1)} \widetilde{f}(\textbf{C}_n(\psi_0))E_n(\psi_0)B_r^\frac{1}{2}\bigg]_{ij}\\\notag
    &- \frac{60x^2}{n^2}\bigg[A_r^\frac{1}{2}\widetilde{f}(\textbf{C}_n(\psi_0))E_n(\psi_0)B_r^\frac{1}{2}\bigg]_{ij}\bigg[A_r^\frac{1}{2}\widetilde{f}(\textbf{C}_n(\psi_0))A_r^\frac{1}{2}\bigg]_{ii}^2[B_r]_{jj}^2\\\notag
    &- \frac{120x^2}{n^2}\bigg[A_r^\frac{1}{2}\widetilde{f}(\textbf{C}_n(\psi_0))E_n(\psi_0)B_r^\frac{1}{2}\bigg]_{ij}\bigg[A_r^\frac{1}{2}\nabla_\psi^{(1)} \widetilde{f}(\textbf{C}_n(\psi_0))E_n(\psi_0)B_r^\frac{1}{2}\bigg]_{ij}\bigg[A_r^\frac{1}{2}\widetilde{f}(\textbf{C}_n(\psi_0))A_r^\frac{1}{2}\bigg]_{ii}[B_r]_{jj}\\\notag
    &-\frac{80\mathbbm{i}x^3}{n^3}\bigg[A_r^\frac{1}{2}\nabla_\psi^{(1)} \widetilde{f}(\textbf{C}_n(\psi_0))E_n(\psi_0)B_r^\frac{1}{2}\bigg]_{ij}\bigg[A_r^\frac{1}{2}\widetilde{f}(\textbf{C}_n(\psi_0))E_n(\psi_0)B_r^\frac{1}{2}\bigg]_{ij}^3\\\notag
    &-\frac{80\mathbbm{i}x^3}{n^3}\bigg[A_r^\frac{1}{2}\widetilde{f}(\textbf{C}_n(\psi_0))E_n(\psi_0)B_r^\frac{1}{2}\bigg]_{ij}^3\bigg[A_r^\frac{1}{2}\widetilde{f}(\textbf{C}_n(\psi_0))A_r^\frac{1}{2}\bigg]_{ii}[B_r]_{jj}\\
    &+ \frac{16x^4}{n^4}\bigg[A_r^\frac{1}{2}\widetilde{f}(\textbf{C}_n(\psi_0))E_n(\psi_0)B_r^\frac{1}{2}\bigg]_{ij}^5\bigg).
\end{align}
We analyze only the first term since the others can be similarly derived.Denote
\begin{align*}
    &\mathcal{X}_n(\psi_0, u) := A_r^\frac{1}{2}I_n(\psi_0, u)A_r^\frac{1}{2}, \quad 
    \mathcal{Y}_n(\psi_0, u) := A_r^\frac{1}{2}I_n(\psi_0, u)E_n(\psi_0)B_r^\frac{1}{2}\\
    &\mathcal{Z}_n(\psi_0) = B_r^\frac{1}{2}E_n(\psi_0)^\top I_n(\psi_0, u)E_n(\psi_0)B_r^\frac{1}{2}, \quad 
    \mathcal{E}_n(\psi_0) := e^{\mathbbm{i}xJ_n^0(\psi_0)}.
\end{align*}
Also for a function $f$, we will use $f^2$ to denote $f*f$. Similar to the steps that led to (\ref{second_term_I_ijr_1}), we find that
\begin{align*}
     &\bigg[A_r^\frac{1}{2}\nabla_\psi^{(4)} \widetilde{f}(\textbf{C}_n(\psi_0))E_n(\psi_0)B_r^\frac{1}{2}\bigg]_{ij}e^{\mathbbm{i}xJ_n^0(\psi_0)}\\      &=\,\mathcal{Q}_{ijr}^{(1)}+\mathcal{Q}_{ijr}^{(2)}+\mathcal{Q}_{ijr}^{(3)}+\mathcal{Q}_{ijr}^{(4)}+\mathcal{Q}_{ijr}^{(5)}+\mathcal{Q}_{ijr}^{(6)},
     \end{align*}
where
\begin{align*}
    \mathcal{Q}_{ijr}^{(1)}&= -\frac{24\mathbbm{i}}{n^2}\int_{-\infty}^\infty u\widehat{f}(u)[B_r]_{jj}^2[\mathcal{X}_n(\psi_0)]_{ii}^2*[\mathcal{Y}_n(\psi_0)]_{ij}(u)\mathcal{E}_n(\psi_0)du\\
    \mathcal{Q}_{ijr}^{(2)} &= \frac{144}{n^3}\int_{-\infty}^\infty u\widehat{f}(u)[B_r]_{jj}[\mathcal{X}_n(\psi_0)]_{ii}*[\mathcal{Y}_n(\psi_0)]_{ii}^3\mathcal{E}_n(\psi_0)du\\
    \mathcal{Q}_{ijr}^{(3)} &=  \frac{144}{n^3}\int_{-\infty}^\infty u\widehat{f}(u)[B_r]_{jj}[\mathcal{X}_n(\psi_0)]_{ii}^2*[\mathcal{Y}_n(\psi_0)]_{ii}*[\mathcal{Z}_n(\psi_0)]_{jj}(u)\mathcal{E}_n(\psi_0)du\\
    \mathcal{Q}_{ijr}^{(4)} &=  \frac{240\mathbbm{i}}{n^4}\int_{-\infty}^\infty u\widehat{f}(u)[\mathcal{X}_n(\psi_0)]_{ii}*[\mathcal{Y}_n(\psi_0)]_{ii}^3*[\mathcal{Z}_n(\psi_0)]_{jj}(u)\mathcal{E}_n(\psi_0)du\\    
    \mathcal{Q}_{ijr}^{(5)} &=  \frac{120\mathbbm{i}}{n^4}\int_{-\infty}^\infty u\widehat{f}(u)[\mathcal{X}_n(\psi_0)]_{ii}^2*[\mathcal{Y}_n(\psi_0)]_{ii}*[\mathcal{Z}_n(\psi_0)]_{jj}^2(u)\mathcal{E}_n(\psi_0)du\\
    \mathcal{Q}_{ijr}^{(6)} &=  \frac{24\mathbbm{i}}{n^4}\int_{-\infty}^\infty u\widehat{f}(u)[\mathcal{Y}_n(\psi_0)]_{ii}^5(u)\mathcal{E}_n(\psi_0)du.
\end{align*}
By arguments similar to those in (\ref{subsec:boundedSpectralNorm}), the norms of $\|A_r^\frac{1}{2}I_n(\psi_0,u)A_r^\frac{1}{2}\|$, $\|n^{-\frac{1}{2}}A_r^\frac{1}{2}I_n(\psi_0,u)E_n(\psi_0)B_r^\frac{1}{2}\|$, $\|n^{-1}B_r^\frac{1}{2}E_n^\top(\psi_0) I_n(\psi_0, u) E_n(\psi_0)B_r^\frac{1}{2}\|$ are bounded. Therefore, we have
\begin{align*}
    &\absmod{\frac{1}{n}\sum_{i=1}^p\sum_{j=1}^n\sum_{r=1}^K \mathbb{E}\widetilde{d}_{ijr}d_{ijr}^4\mathcal{Q}_{ijr}^{(1)}}\\
    \leq& \,\frac{C}{n^3}\sum_{i,j,r} \int_{-\infty}^\infty |u^3\widehat{f}(u)|\mathbb{E}\bigg(|\widetilde{d}_{ijr}d_{ijr}^4|\|A_r^\frac{1}{2}I_n(\psi_0)A_r^\frac{1}{2}\|^2 \,\|A_r^\frac{1}{2}I_n(\psi_0)E_n(\psi_0)B_r^\frac{1}{2}\|\bigg) du\\
    \leq& \,\frac{C}{n^3}\sum_{i,j,r} \int_{-\infty}^\infty |u^3\widehat{f}(u)|\,\mathbb{E}\bigg(|\widetilde{d}_{ijr}d_{ijr}^4| \, \|E_n(\psi_0)\|\bigg)du,\\
    \leq& \, \frac{C}{n^3}\sum_{i,j,r}\int_{-\infty}^\infty |u^3\widehat{f}(u)| \bigg(\mathbb{E}|\widetilde{d}_{ijr}d_{ijr}^4| \, \|E_n + (\rho_{ijr}-1)d_{ijr}e_i\widetilde{e}_j^\top\| \\
    \leq& \, \frac{C}{n^3}\sum_{i,j,r}\int_{-\infty}^\infty |u^3\widehat{f}(u)| \bigg(\mathbb{E}|\widetilde{d}_{ijr}d_{ijr}^4| \, \mathbb{E}\|E_n\| + |1-\rho_{ijr}|\mathbb{E}|\widetilde{d}_{ijr}d_{ijr}^5|\bigg)du\\
     \leq& \, \frac{C}{n^3}\sum_{i,j,r}\int_{-\infty}^\infty |u^3\widehat{f}(u)|\bigg(\delta_nn + \delta_n^2n\bigg)du \leq   C\delta_n \longrightarrow 0.
\end{align*}
Applying similar arguments as above, we get that
\begin{align*}
    \absmod{\frac{1}{n}\sum_{i,j,r}\mathbb{E}\widetilde{d}_{ijr}d_{ijr}^4\bigg(\mathcal{Q}_{ijr}^{(2)} + \mathcal{Q}_{ijr}^{(3)} +\mathcal{Q}_{ijr}^{(4)} +\mathcal{Q}_{ijr}^{(5)} +\mathcal{Q}_{ijr}^{(6)}\bigg)} \leq C\delta_n \int_{-\infty}^\infty (|u|^4 + |u|^5)|\widehat{f}(u)|du \leq C\delta_n,
\end{align*}
which leads to $|\mathcal{I}_2| \leq C\delta_n \longrightarrow 0$.

\subsection{Results used for generalization to non-Gaussian innovations}\label{sec:appendixD}

\begin{lemma}\label{lemma:Duhamel}
    (Duhamel formula:) Let A(t) be a matrix valued $C^\infty$ function (element wise) of $t \in \mathbb{R}$. Then,
    $$\frac{de^{A(t)}}{dt} = \int_0^1 e^{sA(t)}A^\prime(t)e^{(1-s)A(t)}ds.$$
\end{lemma}

\begin{lemma}\label{basic_matrix_derivative}
Suppose $X \in \mathbb{R}^{p \times n}$ is a real matrix.
$$\frac{\partial AXB}{x_{ij}} = A_{\cdot i}B_{j \cdot}.$$
\end{lemma}
\begin{proof}
Let $\{e_i\}_{i=1}^p$ and $\{\tilde{e}_j\}_{j=1}^n$ be canonical bases of $\mathbb{R}^p$ and $\mathbb{R}^n$ respectively. Note that $X = \sum_{i,j} x_{ij}e_i\tilde{e}_j^\top$. Then, it is clear that
    \begin{align*}
        \frac{\partial AXB}{\partial x_{ij}} = Ae_i\tilde{e}_j^\top B = A_{\cdot i} B_{j \cdot}.
    \end{align*}
\end{proof}

\begin{lemma}\label{lemma:derivative_Fn_theta}
Let $i \in [p], j \in [n], r \in [K].$ Recall the definition of $\textbf{C}_{n,ijr}$ from (\ref{defining_Fnijr_omega}). Then,
    $$\bigg[\frac{\partial \textbf{C}_{n,ijr}(\psi)}{\partial \psi}\bigg]_{ab} = \frac{1}{n} \bigg(\bigg[A_r^\frac{1}{2}\bigg]_{ai}\bigg[B_r^\frac{1}{2}E_{n,ijr}^\top(\psi)\bigg]_{jb} + \bigg[E_{n,ijr}(\psi)B_r^\frac{1}{2}\bigg]_{aj}\bigg[A_r^\frac{1}{2}\bigg]_{ib}\bigg).$$
\end{lemma}

\begin{proof}
Starting from
$$E_{n,ijr}(\psi) := \sum_{s\neq r}A_s^\frac{1}{2}D_{n,s}B_s^\frac{1}{2}+ A_r^\frac{1}{2}\bigg(D_{n,r} + (\psi-d_{ijr})e_i\widetilde{e}_j\bigg)B_r^\frac{1}{2},$$ 
we use the arguments of Lemma \ref{basic_matrix_derivative} to get
  $$\frac{\partial E_{n,ijr}(\psi)}{\partial \psi} = A_r^\frac{1}{2}e_i\tilde{e}_j^\top B_r^\frac{1}{2}.$$
Starting from (\ref{defining_En_Fn_theta}), we differentiate both sides and use the above to get 
    \begin{align*}
        \textbf{C}_{n,ijr}(\psi) =& \frac{1}{n}E_{n,ijr}(\psi)E_{n,ijr}^\top(\psi)\\
        \implies \frac{\partial \textbf{C}_{n,ijr}(\psi)}{\partial \psi} =& \frac{1}{n}\bigg(\frac{\partial E_{n,ijr}(\psi)}{\partial \psi}E_{n,ijr}^\top(\psi) + E_{n,ijr}(\psi) \frac{\partial E_{n,ijr}^\top(\psi)}{\partial \psi}\bigg)\\
         =& \frac{1}{n}\bigg(A_r^\frac{1}{2}e_i\widetilde{e}_j^\top B_r^\frac{1}{2}E_{n,ijr}^\top(\psi) + E_{n,ijr}(\psi) B_r^\frac{1}{2}\widetilde{e}_je_i^\top A_r^\frac{1}{2}\bigg)\\
         \implies \bigg[\frac{\partial \textbf{C}_{n,ijr}(\psi)}{\partial \psi}\bigg]_{ab} =& \frac{1}{n}\bigg(e_a^\top  A_r^\frac{1}{2}e_i \widetilde{e}_j^\top B_r^\frac{1}{2}E_{n,ijr}^\top(\psi) e_b +   e_a^\top E_{n,ijr}(\psi) B_r^\frac{1}{2}\widetilde{e}_je_i^\top A_r^\frac{1}{2} e_b\bigg)\\
         =& \frac{1}{n} \bigg(\bigg[A_r^\frac{1}{2}\bigg]_{ai}\bigg[B_r^\frac{1}{2}E_{n,ijr}^\top(\psi)\bigg]_{jb} + \bigg[E_{n,ijr}(\psi)B_r^\frac{1}{2}\bigg]_{aj}\bigg[A_r^\frac{1}{2}\bigg]_{ib}\bigg).
    \end{align*}
\end{proof}

\begin{lemma}\label{lemma:derivative_In_theta}
$$\bigg[\frac{\partial I_n(\psi, t)}{\partial d_{ijr}}\bigg]_{dl} = \frac{\mathbbm{i}}{n} [I_nA_r^\frac{1}{2}]_{di} * [B_r^\frac{1}{2}E_n^\top I_n]_{jl}(t) + \frac{\mathbbm{i}}{n}[I_nE_nB_r^\frac{1}{2}]_{dj} * [A_r^\frac{1}{2}I_n]_{il}(t),$$
where $f*g(t) = \int_0^tf(s)g(t-s)ds$.
\end{lemma}
\begin{proof}
    Using Lemma \ref{lemma:Duhamel} and \ref{lemma:derivative_Fn_theta}, we get
    \begin{align*}
        &\frac{\partial I_{n,ijr}(\psi, t)}{\partial \psi}
        = \sum_{a,b=1}^p  \frac{\partial I_{n,ijr}(\psi, t)}{\partial [\textbf{C}_{n,ijr}(\psi)]_{ab}} \frac{\partial [\textbf{C}_{n,ijr}(\psi)]_{ab}}{\partial \psi}\\
        =& \frac{\mathbbm{i}}{n}\sum_{a,b=1}^p\int_0^t e^{\mathbbm{i}s\textbf{C}_{n,ijr}(\psi)}e_ae_b^\top e^{\mathbbm{i}(t-s)\textbf{C}_{n,ijr}(\psi)}ds \times \bigg(\bigg[A_r^\frac{1}{2}\bigg]_{ai}\bigg[B_r^\frac{1}{2}E_{n,ijr}^\top(\psi)\bigg]_{jb} + \bigg[E_{n,ijr}(\psi)B_r^\frac{1}{2}\bigg]_{aj}\bigg[A_r^\frac{1}{2}\bigg]_{ib}\bigg).
    \end{align*}
    Therefore, it follows that
    \begin{align*}
        &\bigg[\frac{\partial I_{n, ijr}(\psi, t)}{\partial \psi}\bigg]_{dl}\\
        =& \int_0^t  \bigg(\sum_{a=1}^p [I_{n,ijr}(\psi, s)]_{da}[A_r^\frac{1}{2}]_{ai}\bigg)\bigg(\sum_{b=1}^p[B_r^\frac{1}{2}E_{n,ijr}^\top(\psi)]_{jb} [I_{n,ijr}(\psi, t-s)]_{bl}\bigg)ds \\
        &+ \int_0^t  \bigg(\sum_{a=1}^p[I_{n,ijr}(s)]_{da}[E_{n,ijr}(\psi)B_r^\frac{1}{2}]_{aj}\bigg)\bigg(\sum_{b=1}^p[A_r^\frac{1}{2}]_{ib} [I_{n,ijr}(\psi, t-s)]_{bl}\bigg)ds\\
        =&\frac{\mathbbm{i}}{n} [I_{n,ijr}(\psi)A_r^\frac{1}{2}]_{di} * [B_r^\frac{1}{2}E_{n,ijr}^\top(\psi) I_{n,ijr}(\psi)]_{jl}(t) \\
        +& \frac{\mathbbm{i}}{n}[I_{n,ijr}(\psi)E_{n,ijr}(\psi)B_r^\frac{1}{2}]_{dj} * [A_r^\frac{1}{2}I_{n,ijr}(\psi)]_{il}(t).
    \end{align*}
\end{proof}

\begin{lemma}\label{lemma:derivative_Jn_theta}
Let $i \in [p], j \in [n]$. Recall the definitions of $J_{n,ijr}(\psi)$ and   $\widetilde{f}(\textbf{C}_{n,ijr}(\psi))$ from Section \ref{sec:NonGaussian}. Then,
$$\frac{\partial J_{n,ijr}(\psi)}{\partial \psi}
= \frac{2}{n}[A_r^\frac{1}{2}\widetilde{f}(\textbf{C}_{n,ijr}(\psi))E_{n,ijr}(\psi)B_r^\frac{1}{2}]_{ij}.$$
\end{lemma}
\begin{proof}
Recall that $I_{n,ijr}(\psi, u) = e^{\mathbbm{i}u\textbf{C}_{n,ijr}(\psi)}$ and $J_{n,ijr}(\psi) = \operatorname{tr}f(\textbf{C}_{n,ijr}(\psi))$. In particular, $I_{n,ijr}(\psi, \cdot)$ satisfies
\begin{align}\label{property_In_theta}
    I_{n,ijr}(\psi, u)I_{n,ijr}(v) = e^{\mathbbm{i}u\textbf{C}_{n,ijr}(\psi)}e^{\mathbbm{i}u\textbf{C}_{n,ijr}(\psi)} = e^{\mathbbm{i}(u+v)\textbf{C}_{n,ijr}(\psi)} = I_{n,ijr}(u+v).
\end{align}
By the inverse Fourier transform, we have
\begin{align*}
    f(\textbf{C}_{n,ijr}(\theta)) &= \int_{-\infty}^\infty \widehat{f}(u)I_{n,ijr}(\psi, u)du\\
    \implies \frac{\partial J_{n,ijr}(\psi)}{\partial \psi} &= \int_{-\infty}^\infty \widehat{f}(u) \operatorname{tr}\bigg(\frac{\partial I_{n,ijr}(u)}{\partial \psi}\bigg)du = \int_{-\infty}^\infty \widehat{f}(u)\bigg(\sum_{d=1}^p \frac{\partial [I_{n,ijr}(\psi,t)]_{dd}}{\partial \psi}\bigg)du.
\end{align*}
Let $M(v)=I_{n,ijr}(\psi, v)A_r^\frac{1}{2}$ and $N(u) = B_r^\frac{1}{2}E_{n,ijr}^\top(\psi) I_{n,ijr}(\psi, u)$ for $u,v \in  \mathbb{R}$. Then, using Lemma \ref{lemma:derivative_In_theta}, we have the following simplification:
$$\frac{\partial [I_{n,ijr}(\psi, t)]_{dd}}{\partial \psi} = \frac{2\mathbbm{i}}{n} (M_{di}*N_{jd})(t).$$
Thus, we have the following simplification:
\begin{align*}
\frac{\partial J_{n,ijr}(\psi)}{\partial \psi}
&=
\frac{2\mathbbm{i}}{n}\int_{-\infty}^{\infty}\widehat{f}(u)
\sum_{d=1}^{p}
(M_{di}*N_{jd})(u)\,du \\
&= \frac{2\mathbbm{i}}{n}\int_{-\infty}^{\infty}\widehat{f}(u) \bigg(\int_0^u \sum_{d=1}^p M_{di}(v)N_{jd}(u-v)dv \bigg)du\\ 
&= \frac{2\mathbbm{i}}{n}\int_{-\infty}^{\infty}\widehat{f}(u) \bigg(\int_0^u [N(u-v)M(v)]_{ji}(u)dv \bigg)du.
\end{align*}
Using property (\ref{property_In_theta}), for fixed $u$, for any $0\leq v \leq u$, we have
\begin{align}\label{NM_ij}
    N(u-v)M(v) = B_r^\frac{1}{2}E_{n,ijr}^\top(\psi) I_{n,ijr}(\psi, u-v)I_{n,ijr}(\psi, v)A_r^\frac{1}{2} = B_r^\frac{1}{2}E_{n,ijr}^\top(\psi) I_{n,ijr}(\psi, u)A_r^\frac{1}{2}.
\end{align}
Therefore, using (\ref{NM_ij}) and (\ref{defining_f_tilde_Fn_theta}), it follows that
\begin{align*}
\frac{\partial J_{n,ijr}(\psi)}{\partial \psi}
&= \frac{2\mathbbm{i}}{n} \int_{-\infty}^\infty u\widehat{f}(u)[B_r^\frac{1}{2}E_{n,ijr}^\top(\psi) I_{n,ijr}(\psi, u)A_r^\frac{1}{2}]_{ji}du\\
&= \frac{2\mathbbm{i}}{n} \int_{-\infty}^\infty u\widehat{f}(u)[A_r^\frac{1}{2}I_{n,ijr}(\psi, u)E_{n,ijr}(\psi)B_r^\frac{1}{2}]_{ij}du\\
&= \frac{2}{n}[A_r^\frac{1}{2}\widetilde{f}(\textbf{C}_{n,ijr}(\psi))E_nB_r^\frac{1}{2}]_{ij}.
\end{align*}
\end{proof}

\begin{lemma}
    Let $\widehat{f}$ be defined as in (\ref{fourier_transform}). For $k \in \mathbb{N}$, we have $$\widehat{f^{(k)}}(u) = u^k\widehat{f}(u).$$  As a result, if $\widehat{f^{(k)}} \in L^1(\mathbb{R})$, we have
    $$\int_{-\infty}^\infty |u|^k\,|\widehat{f}(u)|\,du = \int_{-\infty}^\infty |\widehat{f^{(k)}}(u)|\,du < \infty.$$
\end{lemma}

\begin{lemma}\label{lemma:5_by_4}
Suppose $A,B,C$ are random matrices whose spectral norms have bounded moments. Then for some constant $C>0$,
$$\absmod{\sum_{i,j}\mathbb{E}A_{ii}B_{ij}C_{jj}} \leq Cn^{5/4}, \quad \mathbb{E} \absmod{\sum_{i,j}A_{ij}B_{ij}C_{ij}} \leq Cn^{5/4}.$$
\end{lemma}
\begin{proof}
    Due to the given bounds, we have the following results:
    \begin{enumerate}
        \item $\sum_{i=1}^n\mathbb{E}|A_{ii}|^2 \leq n\mathbb{E}\|A\|^2 < Cn$,
        \item $\sum_{j,k}|(B^*B)_{jk}|^2 = \|B^*B\|_F^2 \leq n\|B\|^4$, and
        \item $\sum_{j,k}\mathbb{E}|C_{jj}C_{kk}|^2 \leq Cn^2$.
    \end{enumerate}
    By Cauchy-Schwarz, we have
    \begin{align*}
        \absmod{\sum_{i,j}\mathbb{E}A_{ii}B_{ij}C_{jj}} &\leq \bigg(\sum_i \mathbb{E}|A_{ii}|^2\bigg)^\frac{1}{2} \bigg(\sum_i \mathbb{E}\bigg|\sum_j B_{ij}C_{jj}\bigg|^2\bigg)^\frac{1}{2}\\
        &\leq C\sqrt{n}\bigg(\sum_i \mathbb{E} \bigl(\sum_j B_{ij}C_{jj}\bigr)\overline{\bigl(\sum_k B_{ik}C_{kk}\bigr)}\bigg)^\frac{1}{2}\\
        & = C\sqrt{n}\bigg(\mathbb{E}\sum_{j,k} \bigl(\sum_i B_{ij}\overline{B}_{ik}\bigr)C_{jj}\overline{C}_{kk}\bigg)^\frac{1}{2}\\
        &=C\sqrt{n}\bigg(\mathbb{E}\sum_{j,k} (B^*B)_{jk}C_{jj}\overline{C}_{kk}\bigg)^\frac{1}{2}\\
        &\leq C\sqrt{n} \bigg(\sum_{j,k}\mathbb{E}|B^*B)_{jk}|^2\bigg)^\frac{1}{4} \bigg(\sum_{jk}\mathbb{E}|C_{jj}C_{kk}|^2\bigg)^\frac{1}{4}\\
        &\leq C\sqrt{n}\,n^\frac{1}{4}(n^2)^\frac{1}{4} = Cn^\frac{5}{4}.
    \end{align*}

This proves the first result. For the other, we note that 
\begin{enumerate}
    \item $\sum_{i,j}|A_{ij}|\,B_{ij}|\,|C_{ij}\|^4 \leq Cn^2$,
    \item $|A_{ij}B_{ij}C_{ij}| = (|A_{ij}|\,|B_{ij}|)^\frac{1}{2}(|A_{ij}|\,|B_{ij}|\,|C_{ij}|^2)^\frac{1}{2}$, and 
    \item $|A_{ij}|\,|B_{ij}|\,|C_{ij}|^2 = (|A_{ij}|\,|B_{ij}|)^\frac{1}{2} (|A_{ij}|\,|B_{ij}|\,|C_{ij}|^4)^\frac{1}{2}$.
\end{enumerate}
Applying Cauchy-Schwarz, we get
\begin{align*}
    \absmod{\sum_{i,j}\mathbb{E}A_{ij}B_{ij}C_{ij}} &\leq \bigg(\sum_{i,j}\mathbb{E}|A_{ij}|\,|B_{ij}|\bigg)^\frac{1}{2}\bigg(\sum_{i,j}\mathbb{E}|A_{ij}|\,|B_{ij}|\,|C_{ij}|^2\bigg)^\frac{1}{2}\\
    & \leq \bigg(\sum_{i,j}\mathbb{E}|A_{ij}|\,|B_{ij}|\bigg)^\frac{1}{2}\bigg(\sum_{i,j}\mathbb{E}|A_{ij}|\,|B_{ij}|\bigg)^\frac{1}{4} \bigg(\sum_{i,j}\mathbb{E}|A_{ij}|\,|B_{ij}|\,|C_{ij}|^4\bigg)^\frac{1}{4}\\
    & \leq \bigg(\sum_{i,j}\mathbb{E}|A_{ij}|\,|B_{ij}|\bigg)^\frac{3}{4}\bigg(\sum_{i,j}\mathbb{E}|A_{ij}|\,|B_{ij}|\,|C_{ij}|^4\bigg)^\frac{1}{4}\\
    &\leq C\sqrt{n}\bigg(\sum_{i,j}\mathbb{E}|A_{ij}|^2\bigg)^\frac{3}{8} \bigg(\sum_{i,j}\mathbb{E}|B_{ij}|^2\bigg)^\frac{3}{8} \leq Cn^{\frac{1}{2} + \frac{3}{8} + \frac{3}{8}} = Cn^\frac{5}{4}.
\end{align*}
\end{proof}

\subsection{Justification for bounds on spectral norms}\label{subsec:boundedSpectralNorm}
 Here, we revert to the true notations, i.e. $E_n$ refers to the matrix defined in (\ref{defining_En_Fn_theta}) and $E_{n,ijr}(\psi)$ refers to (\ref{defining_Enijr_omega}) and similarly for $I_n, I_{n,ijr}(\psi)$, etc. This is to emphasize the fact that though we have replaced one innovation entry $d_{ijr}(\theta)$ by $0$, we still have the desired control on the matrix norms and hence on the derivatives of the associated functions.
 We want to establish that the spectral norms of the following matrices
\begin{enumerate}
    \item $A_r^\frac{1}{2}I_{n,ijr}(0,u)A_r^\frac{1}{2}$,
    \item $n^{-\frac{1}{2}}A_r^\frac{1}{2}I_{n,ijr}(0,u)E_{n,ijr}(0)B_r^\frac{1}{2}$,
    \item $n^{-1}B_r^\frac{1}{2}E_n^\top(0) I_n(0, u) E_n(0)B_r^\frac{1}{2}$
\end{enumerate}
have bounded moments. Observe that
$$\|E_{n,ijr}(0)\| \leq \|E_n\|+\|(E_{n,ijr}(0) - E_n)/\sqrt{n}\| =\|E_n\| +  |d_{ijr}|\|A_r^\frac{1}{2}e_i\widetilde{e}_j^\top B_r^\frac{1}{2}\|.$$
As a result,
\begin{align*}
    \mathbb{E}\|n^{-\frac{1}{2}}E_{n,ijr}\|^k \leq C\bigg(\mathbb{E}\|n^{-\frac{1}{2}}E_n\|^k + C_1\mathbb{E}|n^{-\frac{1}{2}}d_{ijr}|^k\bigg).
\end{align*}
Clearly, $\mathbb{E}|n^{-\frac{1}{2}}d_{ijr}|^k \to 0$ for any $k \geq 1$ and $\|n^{-\frac{1}{2}}E_n\|$ have bounded moments, therefore it implies that 
$$\mathbb{E}\|n^{-\frac{1}{2}}E_{n,ijr}(0)\|^k \leq C_k.$$
Therefore, $\|\textbf{C}_{n,ijr}(0)\|$ is also bounded. This immediately implies that $\|I_{n,ijr}(0,u)\| = \|e^{\mathbbm{i}u\textbf{C}_{n,ijr}(0)}\| \leq e^{\|\textbf{C}_{n,ijr}(0)\|}$ is also bounded. 

\section{Proof of Tightness of \texorpdfstring{$M_{n1}(z)$}{}}\label{sec:Tightness}
Having developed all the necessary tools in Lemma \ref{lemma:bounds_tightness_section}, we will now show that 
\begin{align}\label{goal_tightness}
    \underset{n;z_1,z_2 \in \mathcal{C}_n}{\sup} \frac{\mathbb{E}|M_{n1}(z_1) - M_{n1}(z_2)|^2}{|z_1-z_2|^2} < \infty.
\end{align}
By the resolvent identity (\ref{Resolvent_Identity}), we have 
\begin{align*}
    s_{n}(z_1) - s_{n}(z_2) = \frac{1}{p}\operatorname{tr}(Q(z_1) - Q(z_2)) = \frac{1}{p}(z_1 - z_2)\operatorname{tr}Q(z_1)\operatorname{tr}Q(z_2).
\end{align*}
Retaining our earlier practice of using $Q^{(i)}, Q_j^{(i)} \beta_j^{(i)}, b_j^{(i)}$ instead of $Q(z_i), Q_j(z_i), \beta_j(z_i), b_j(z_i)$ as and when convenient, we have 
$$Q(z_i) = Q_j(z_i) - \beta_j(z_i)Q_j(z_i)x_jx_j^*Q_j(z_i), \quad i=1,2.$$ 
Using these we simplify the expression in (\ref{goal_tightness}) as follows:
\begin{align}\label{decomposing_lipschitz_ratio}
    &\frac{M_{n1}(z_1) - M_{n1}(z_2)}{z_1 - z_2}\\ \notag
    =&\, \frac{p(s_n(z_1) - s_n(z_2) - \mathbb{E}(s_n(z_1) - s_n(z_2)))}{z_1 - z_2}\\ \notag
    =& \operatorname{tr}Q(z_1)Q(z_2) - \mathbb{E}\operatorname{tr}Q(z_1)Q(z_2)\\ \notag
    =& \sum_{j=1}^n (\mathbb{E}_j - \mathbb{E}_{j-1})Q(z_1)Q(z_2)\\ \notag
    =& \sum_{j=1}^n \Delta_j\beta_j(z_1)\beta_j(z_2)(x_j^*Q_j(z_1)Q_j(z_2)x_j)^2\\ \notag
    & -\sum_{j=1}^n\Delta_j \beta_j(z_1)x_j^*Q_j^2(z_1)Q_j(z_2)x_j -\sum_{j=1}^n \Delta_j\beta_j(z_2)x_j^*Q_j(z_1)Q_j^2(z_2)x_j\\ \notag
    :=& \mathcal{P}_1 + \mathcal{P}_2 + \mathcal{P}_3,
\end{align}
where, we used the fact that $\Delta_jQ_j(z_1)Q_j(z_2)=0$.
It suffices to show that $\mathbb{E}|\mathcal{P}_1 + \mathcal{P}_2+ \mathcal{P}_3|^2$ is bounded. Similar to $\epsilon_j(z)$ defined in (\ref{essential_quantities}), we define the quantity: 
\begin{align}\label{defining_rhoj}
    \rho_{j}(z) &:= x_j^*Q_j(z)x_j - \mathbb{E}\operatorname{tr}\mathcal{A}_j = \epsilon_j(z) + \operatorname{tr}(\mathcal{A}_j) - \mathbb{E}\operatorname{tr}\mathcal{A}_j.
\end{align}
Recall the definitions of $\beta_j, \widetilde{\beta}, b_j$ from (\ref{essential_quantities})). We start by decomposing $\beta_j(z)$  using $b_j(z)$ as the deterministic first order term (instead of $\widetilde{\beta}(z)$ as was done in (\ref{beta_betaTilde_relationship})).
\begin{align}\label{beta_b_relationship}
    \beta_j(z) &= b_j(z) - \beta_j(z)b_j(z)\rho_j(z)
    = b_j(z) - b_j^2(z)\rho_j(z) + \beta_j(z)b_j^2(z)\rho_j^2(z).
\end{align}
From Lemma \ref{lemma:bounds_tightness_section}, $\sup_{j\in[n], z\in\mathcal{C}_n}|b_j(z)|$ is bounded. We now analyze the first term in the decomposition derived in (\ref{decomposing_lipschitz_ratio}). Denote $T_j = x_j^*Q_j^{(1)}Q_j^{(2)}x_j$ for simplicity. Applying the property (\ref{beta_b_relationship}) sequentially, we get the following:
\begin{align*}
    \mathcal{P}_1 =& \sum_{j=1}^n \Delta_j\beta_j(z_1)\beta_j(z_2)T_j^2\\
    =& \sum_{j=1}^n \Delta_j(b_j(z_1) - b_j(z_1)\beta_j(z_1)\rho_j(z_1))\beta_j(z_2)T_j^2\\
    =& \sum_{j=1}^n \Delta_jb_j(z_1)(b_j(z_2) - b_j(z_2)\beta_j(z_2)\rho_j(z_2))T_j^2\\
    &- \sum_{j=1}^n b_j(z_1)\Delta_j\beta_j(z_1)\beta_j(z_2)\rho_j(z_1)T_j^2\\
    =& \sum_{j=1}^n b_j(z_1)b_j(z_2)\Delta_jT_j^2 - \sum_{j=1}^n b_j(z_1)b_j(z_2)\Delta_j\beta_j(z_2)\rho_j(z_2)T_j^2\\
    &- \sum_{j=1}^n b_j(z_1)\Delta_j\beta_j(z_1)\beta_j(z_2)\rho_j(z_1)T_j^2\\
    :=& \mathcal{P}_{11} + \mathcal{P}_{12} + \mathcal{P}_{13}.
\end{align*}
Note that  $\mathbb{E}_{-j}T_j = n^{-1}\operatorname{tr}Q_j(z_1)Q_j(z_2)\mathcal{A}_j$. We have
\begin{align*}
  &|T_j + \mathbb{E}_{-j}T_j|^2 = |T_j - \mathbb{E}_{-j}T_j + 2\mathbb{E}_{-j}T_j|^2  \leq 2|T_j-\mathbb{E}_{-j}T_j|^2 + 8|\mathbb{E}_{-j}T_j|^2.
  \end{align*}
Observe that $\mathbb{E}|\mathbb{E}_{-j}T_j|^2$ is bounded by Lemma \ref{lemma:bounds_tightness_section}. Therefore, there exists absolute constants $C_1, C_2>0$ such that 
  \begin{align}\label{basic_algebra}
& \mathbb{E}|T_j + \mathbb{E}_{-j}T_j|^2 \leq \, C_1 + C_2|T_j - \mathbb{E}_{-j}T_j|^2\\\notag
\implies& \mathbb{E}|T_j^2 - \mathbb{E}_{-j}^2T_j|^2 \leq  C_1|T_j-\mathbb{E}_{-j}T_j|^2 + C_2 |T_j-\mathbb{E}{_-j}T_j|^4.
\end{align}
Using Lemma \ref{lemma:bounds_tightness_section} and (\ref{basic_algebra}), we then have
\begin{align*}
    \mathbb{E}|\mathcal{P}_{11}|^2 &=\mathbb{E}\absmod{\sum_{j=1}^n b_j(z_1)b_j(z_2)\Delta_jT_j^2}^2 \\
    &=\mathbb{E}\absmod{\sum_{j=1}^n b_j(z_1)b_j(z_2)\Delta_j\bigg(T_j^2 - \mathbb{E}_{-j}^2T_j\bigg)}^2 \\
    &\leq C\sum_{j=1}^n \mathbb{E}|T_j^2 - \mathbb{E}_{-j}T_j^2|^2\\
    &= C_1\sum_{j=1}^n\mathbb{E}|T_j - \mathbb{E}_{-j}T_j|^4 + C_2\sum_{j=1}^n\mathbb{E}|T_j - \mathbb{E}_{-j}T_j|^2\\
    &\leq \frac{C_1}{n} + C_2 \leq C.
\end{align*}
Now analyzing $\mathcal{P}_{12}$, we get
\begin{align*}
    \mathbb{E}|\mathcal{P}_{12}|^2 &= \mathbb{E}\absmod{\sum_{j=1}^n b_j(z_1)b_j(z_2)\Delta_j\beta_j(z_2)\rho_j(z_2)T_j^2}^2\\
    &\leq C_1\mathbb{E}\absmod{\sum_{j=1}^n \beta_j(z_2)\rho_j(z_2)(T_j - \mathbb{E}_{-j}T_j)^2}^2 + C_2\mathbb{E}|\beta_j(z_2)\rho_j(z_2)|^2 \\
    &\leq C_1\sum_{j=1}^n\mathbb{E}^\frac{1}{2}\absmod{\rho_j^2(z_2)(T_j - \mathbb{E}_{-j}T_j)^2}^4\mathbb{E}^\frac{1}{2}|\beta_j(z_2)|^4+ C_2\sum_{j=1}^n\mathbb{E}^\frac{1}{2}|\beta_j(z_2)|^4 \mathbb{E}^\frac{1}{2}|\rho_j(z_2)|^4\\
    &\leq \frac{C_1}{n^2} + C_2 \leq C.
\end{align*}
By similar arguments, we also get $|\mathbb{E}\mathcal{P}_{13}|^2 \leq C$. Therefore, we have $$\mathbb{E}|\mathcal{P}_1|^2 \leq C.$$
Due to the similarity of $\mathcal{P}_2$ and $\mathcal{P}_3$ in (\ref{decomposing_lipschitz_ratio}), we analyze only one of them below. Now let $T_j = x_j^*Q_j^2(z_1)Q_j(z_2)x_j$. Then, $\mathbb{E}_{-j}T_j = \mathcal{A}_jQ_j^2(z_1)Q_j(z_2)$. Using (\ref{beta_b_relationship}), we have
\begin{align*}
    \mathcal{P}_2 &= -\sum_{j=1}^n\Delta_j \beta_j(z_1)T_j\\
    &= -\sum_{j=1}^n \Delta_j(b_j(z_1) - b_j(z_1)\beta_j(z_1)\rho_j(z_1))T_j\\
    &= -\sum_{j=1}^n b_j(z_1)\Delta_jT_j + \sum_{j=1}^n b_j(z_1)\Delta_j \beta_j(z_1)\rho_j(z_1)T_j.
\end{align*}
This leads to the following bound.
\begin{align*}
    \mathbb{E}|\mathcal{P}_2|^2 &\leq C_1\mathbb{E}|\sum_{j=1}^n b_j(z_1) \Delta_jT_j|^2 + C_2 \mathbb{E}|\sum_{j=1}^n b_j(z_1)\Delta_j \beta_j(z_1)\rho_j(z_1)T_j|^2\\
    &\leq C_1 \sum_{j=1}^n\mathbb{E}|T_j - \mathbb{E}_{-j}T_j|^2 + C_2 \sum_{j=1}^n\mathbb{E}|\beta_j(z_1)\rho_j(z_1)T_j|^2\\
    &\leq C_1 +  C_2\sum_{j=1}^n \mathbb{E}|\beta_j(z_1)\rho_j(z_1)(T_j - \mathbb{E}_{-j}T_j)|^2 + C_3 \sum_{j=1}^n\mathbb{E}|\beta_j(z_1)\rho_j(z_1)|^2\leq C.
\end{align*}

Therefore, we conclude that 
\begin{align*}
    \underset{n;z_1,z_2 \in \mathcal{C}_n}{\sup} \frac{\mathbb{E}|M_{n1}(z_1) - M_{n1}(z_2)|^2}{|z_1-z_2|^2} \leq \underset{n;z_1,z_2 \in \mathcal{C}_n}{\sup} \mathbb{E}|\mathcal{P}_1+\mathcal{P}_2+\mathcal{P}_3|^2 \leq C.
\end{align*}
By Lemma \ref{lemma:tightness_condition}, this implies that $\widehat{M}_{n1}(z)$ defined in (\ref{defining_truncated_contour}) is tight.

\section{Results related to the Covariance kernel}
\subsection{Proof of Lemma \ref{lemma:Bounds_Important_Quantities}}\label{sec:proofLemma4.1}
\begin{proof}
    \textbf{Proof of (1):} 
    
Let $S_n = P^*\Lambda P$ be a spectral decomposition of $S_n$. Observe that for $x \in \mathbb{C}^p$, we have
\begin{align*}
    \Im(zx^*Q_jx) &= \frac{1}{2\mathbbm{i}}\bigg(x^*(S_{nj}/z - I_p)^{-1}x - x^*(S_{nj}/\bar{z} - I_p)^{-1}x\bigg)\\
    &= \frac{1}{2\mathbbm{i}}\sum_{i=1}^p \bigg(\frac{|y_i|^2}{\lambda_i/z-1} - \frac{|y_i|^2}{\lambda_i/\bar{z}-1}\bigg)\\
    &= \frac{1}{2\mathbbm{i}}\sum_{i=1}^p |y_i|^2 \frac{z(\lambda_i - \bar{z}) - \bar{z}(\lambda_i - z)}{|\lambda_i - z|^2} = 2v\sum_{i=1}^n \frac{|y_i|^2}{|\lambda_i - z|^2} \geq 0,
\end{align*}
where $y = Px$. Therefore, we have 
\begin{align}\label{betaBound}
    |\beta(z)| \leq \frac{|z|}{v}, \quad \frac{1}{1 + x_j^*Q_jx_j} \leq \frac{|z|}{v}.
\end{align}
Note that $h_{rn}(z) = \frac{1}{p}\operatorname{tr}{A_rQ_j}$ is the Stieltjes transform of a measure with positive support, thus implying that $\Im(z\operatorname{tr}(A_rQ_j)) \geq 0$. In particular, we have 
\begin{align}\label{betaTildeBound}
    \absmod{{z(1 + \frac{1}{n}\operatorname{tr}(A_rQ_j))}} \geq {\Im(z) + \frac{1}{n}\Im(z\operatorname{tr}(A_rQ_j))} \geq {v} \implies |\widetilde{\beta}_j(z)| \leq \frac{|z|}{v}.
\end{align}

By similar arguments, we have $\Im(z\mathbb{E}\operatorname{tr}A_rQ_j) \geq 0$ and $\Im(z\mathbb{E}\operatorname{tr}A_rQ) \geq 0$. This is sufficient to conclude that 
\begin{align}
    |b_j(z)| \leq \frac{|z|}{v} \text{, and } |\psi_j(z)| \leq \frac{|z|}{v}.
\end{align}
Now, we establish the bounds on $\|R_j(z)\|$. For arbitrary $\boldsymbol{\lambda}=(\lambda_1, \ldots, \lambda_K)^\top \in \mathbb{R}_+^K$, we have
\begin{align*}
    &\Im\bigg(z - \frac{1}{n}\sum_{i\neq j}\psi_i(z)\boldsymbol{\theta}_i^\top \boldsymbol{\lambda}\bigg)\\
    &= \Im\bigg(z - \frac{1}{n}\sum_{i\neq j}\frac{\boldsymbol{\theta}_i^\top \boldsymbol{\lambda}}{1 + c_n \boldsymbol{\theta}_j^\top \mathbb{E}\textbf{h}_n}\bigg)\\ 
    &=v - \Im\bigg(\frac{1}{n}\sum_{i\neq j}\frac{(1+c_n\boldsymbol{\theta}_i^\top \overline{\mathbb{E}\textbf{h}_n})\boldsymbol{\theta}_i^\top \boldsymbol{\lambda}}{|1+c_n\boldsymbol{\theta}_i^\top \mathbb{E}\textbf{h}_n|^2}\bigg) \geq v \text{, since } \Im(\mathbb{E}\textbf{h}_n) \geq 0.
\end{align*}
Therefore, $\|R_j(z)\|$ and similarly $\|R(z)\|$ are bounded by $v^{-1}$.
\end{proof}

\begin{lemma}\label{lemma:derivative_quadratic_forms}
    Recall the definitions of $\epsilon_j(z), \gamma_j(z)$ for $r \in [K], j \in [n]$ from (\ref{essential_quantities}). For any $z \not \in \operatorname{supp}(F^{S_n})$, the following holds:
    $$\frac{d}{dz}\bigg(\epsilon_j(z)\bigg) = \gamma_j(z).$$
\end{lemma}

\begin{proof}
    Using a spectral decomposition of $S_{nj}$, we rewrite $Q_j = P(\Lambda_j - zI_p)P^*$. This allows to simplify the tracial term as follows:
$$\operatorname{tr}A_rQ_j(z) = \operatorname{tr}(\underbrace{P^*A_rP}_{=D_r}(\Lambda_j - zI_p))^{-1} = \sum_{k=1}^p \frac{d_{kk}^{(r)}}{\lambda_k^{(j)} - z},$$ 
where $\{\lambda_k^{(j)}\}_{k=1}^p$ are the eigenvalues of $S_{nj}$ and $\{d_{kk}^{(r)}\}_{k=1}^p$ are the diagonal elements of $D_r$. Since $S_{nj}$ is independent of $z$, we have
\begin{align*}
    Q_j^{-1} = (S_{nj} - zI_p) \implies \frac{d}{dz}Q_j^{-1} = -I_p.
\end{align*}
As a result,
\begin{align*}
    Q_jQ_j^{-1} = I_p \implies Q_j\bigg(\frac{d}{dz}Q_j^{-1}\bigg) + \bigg(\frac{d}{dz}Q_j\bigg)Q_j^{-1}=0 \implies \frac{d}{dz}Q_j = Q_j^{2}.
\end{align*}
Since $\{\lambda_k^{(j)}\}_{k=1}^p \subset \mathbb{R}_+$, the map $z \mapsto 1/(\lambda_k^{(j)}-z)$ is holomorphic for each $k \in [p]$. Now we differentiate $\epsilon_j(z)$ to get 
\begin{align*}
    \frac{d}{dz}\big(\epsilon_j(z)\big)= &\frac{d}{dz} \bigg(x_j^*Q_j(z)x_j - n^{-1}\sum_{r=1}^K\operatorname{tr}(A_rQ_j)\bigg)\\
    =&  x_j^*\bigg(\frac{d}{dz}Q_j\bigg)x_j -\frac{1}{n}\sum_{r=1}^K\sum_{k=1}^p\frac{d}{dz}\bigg(\frac{d_{kk}^{(r)}}{\lambda_k^{(j)} - z}\bigg)\\
    =& x_j^*Q_j^{-2}x_j-\frac{1}{n}\sum_{r=1}^K\sum_{k=1}^p\frac{d_{kk}^{(r)}}{(\lambda_k^{(j)} - z)^2}\\
    =& x_j^*Q_j^2(z)x_j - n^{-1}\sum_{r=1}^K\operatorname{tr}(A_rQ_j^2) = \gamma_j(z).
\end{align*}
\end{proof}

\begin{lemma}\label{lemma:transfer_of_derivative}
    Denote $f_n(z_1, z_2) = \sum_{j=1}^n \mathbb{E}_{j-1}[\mathbb{E}_j(\widetilde{\beta}_j(z_1)\epsilon_j(z_1)) \mathbb{E}_j(\widetilde{\beta}_j(z_2)\epsilon_j(z_2))]$. If there exists $\{\widehat{\sigma}_n\}$ such that
    $$\frac{f_n(z_1, z_2)}{\widehat{\sigma}_n^2(z_1, z_2)} \xrightarrow{p} 1,$$
    then we also have
    $$\frac{\Phi_n(z_1, z_2))}{\sigma_n^2(z_1, z_2)} \xrightarrow{p} 1, \quad \text{where } \sigma_n^2(z_1, z_2) = \frac{\partial^2}{\partial z_2 \partial z_1} \widehat{\sigma}_n^2(z_1, z_2).$$
\end{lemma}
\begin{proof}
By Lemma \ref{lemma:Bounds_Important_Quantities} and Lemma \ref{lemma:bounds_tightness_section}, for any $z_1, z_2 \in \mathcal{C}_u$, we observe that 
\begin{align*}
    |f_n(z_1, z_2)| \leq \frac{|z_1||z_2|}{v_0^2}\bigg(\mathbb{E}_{j-1}|\epsilon_j(z)|^2\mathbb{E}_{j-1}|\epsilon_j(z)|^2\bigg)^\frac{1}{2} \leq \frac{C\tau^2|z_1|||z_2}{nv_0^2}.
\end{align*}
So, $|f_n(z_1, z_2)|$ is upper bounded by a quantity depending on $z_1, z_2, v_0$. The rest of the proof uses Lemma \ref{lemma:Vitali_convergence} as done in page 571 of \cite{BaiSilv04}. The details are extremely similar and  thus omitted.
\end{proof}

\subsection{Proof of Lemma \ref{lemma:bounds_tightness_section}}
\begin{proof}
    \textbf{Proof of 1:}
    For $z \in \mathcal{C}_u(v_0)$, we have the stronger (and trivial) result
    $$\|Q_*(z)\| \leq v_0^{-1}.$$
    For $z \in \mathcal{C}_l \bigcup \mathcal{C}_r$ and any positive integer $k$, we can choose an integer $l$ ($>k+\alpha$, see (\ref{bad_event_probability})) large enough such that 
\begin{align}\label{norm_resolvent_moment_bound}
    \mathbb{E}\|Q(z)\|^k &= \mathbb{E}\|Q(z)\|^k \mathbbm{1}_{\mathcal{G}}
    + \mathbb{E}\|Q(z)\|^k\mathbbm{1}_{\mathcal{G}^c}\\\notag
    &\leq M_{good}^k + \frac{1}{v^k}\mathbb{P}(\mathcal{G}^c) \\\notag
    &\leq C_1 + C_2 n^k\epsilon_n^{-k}n^{-l} \leq C_k.
\end{align}
The proofs for $Q_j$ and $Q_{ij}$ are similar.

\textbf{Proof of 2:} 
    We have $$\psi_j(z) = \frac{1}{1 + c_n\boldsymbol{\theta}_j^\top \mathbb{E}\textbf{h}_n(z)}.$$
It suffices to show the denominator is lower bounded. Noting that $\|Q(z)\| \leq n/\epsilon_n \leq n^{1+\alpha}$ on the \textit{bad} event, we see that for any $r \in [K]$,
\begin{align}\label{sup_exp_hrn_bound}
    \underset{n;z \in \mathcal{C}_n}{\sup} |\mathbb{E}h_{rn}(z)| &= \underset{n;z \in \mathcal{C}_n}{\sup}\bigg|\frac{1}{p}\mathbb{E}\operatorname{tr}(A_rQ(z))\mathbbm{1}_{\mathcal{G}}\bigg| + \underset{n;z \in \mathcal{C}_n}{\sup}\bigg|\frac{1}{p}\mathbb{E}\operatorname{tr}(A_rQ(z))\mathbbm{1}_{\mathcal{G}^c}\bigg|\\ \notag
    &\leq \tau M_{good} + \tau n^{1+\alpha}\mathbb{P}(\mathcal{G}^c) \\ \notag
    &< \frac{\tau}{4Kc\tau^2}  + \frac{1}{8Kc\tau} = \frac{3}{8Kc\tau},
\end{align}
where we used the technical requirement made in (\ref{contour_gap_technical_condition}). By triangle and Cauchy Schwarz inequalities, we get
    \begin{align*}
        |1 + c_n\boldsymbol{\theta}_j^\top \mathbb{E}\textbf{h}_n(z)| \geq 1 - 2c\|\boldsymbol{\theta}_{j}\|\, \|\mathbb{E}\textbf{h}_{n}(z)\| \geq 1 - 2Kc\tau\frac{3}{8Kc\tau} = \frac{1}{4}.
    \end{align*}
For the uniform bound on $|b_j(z)|$, by the same argument, we obtain 
\begin{align}
    \sup_{n;j\in [n],z\in \mathcal{C}_n} \bigg|\frac{1}{p}\mathbb{E}\operatorname{tr}(A_rQ_j(z))\bigg| < \frac{3}{8Kc\tau}.
\end{align}

\textbf{Proof of 3:}
By Woodbury's formula, we have
\begin{align*}
        x_j^*Q(z)x_j = \beta_j(z)x_j^*Q_j(z)x_j = 1 - \beta_j(z).
\end{align*}
Therefore, we have the following generic bound:
\begin{align*}
    |\beta_j(z)| \leq 1 + |x_j^*Q_jx_j| \leq 1 + \|x_j\|^2\,\|Q_j\|.
\end{align*}
On $\mathcal{G} \bigcap \mathcal{G}_0$, we have 
\begin{align*}
    \|x_j\| = \|S_n -S_{nj}\| \leq \lambda_{\max}(S_n) + \lambda_{\min}(S_{nj}) \leq 2(|\eta_r| + |\eta_l|) := C_{good}.
\end{align*}
On $\mathcal{G}^c \bigcup \mathcal{G}_0^c$, we have
\begin{align*}
    \|x_j\|^2 \leq \sum_{r=1}^K K\|x_{r,j}\|^2 \leq \frac{K\tau}{n}\sum_{i=1}^p |z_{ij}^{(r)}|^2 \leq \frac{K\tau pn\delta_n^2}{n} \leq Cn.
\end{align*}
Finally, we have
\begin{align*}
    \mathbb{E}|\beta_j(z)|^k &\leq 3^{k-1}\bigg(1 + C_{good}^{2k}M_{good}^k + Cn^{k(2+\alpha)}\bigl(\mathbb{P}(\mathcal{G}^c) + \mathbb{P}(\mathcal{G}_j^c)\bigr)\bigg) \leq C_p.
\end{align*}
This takes care of the moments of $|\beta_j(z)|$. Recall that 
$$\widetilde{\beta}_j(z) = \frac{1}{1 + n^{-1}\operatorname{tr}\textbf{A}_jQ_j(z)}.$$
On $\mathcal{G}_0$, we see that 
\begin{align}\label{sup_hrnj_bound}
    \bigg|\frac{1}{n}\operatorname{tr}(\textbf{A}_jQ_j(z))\bigg| \leq c_n\sum_{r=1}^K \theta_{rj}\bigg|\frac{1}{p}\operatorname{tr}(A_rQ_j(z))\bigg| \leq 2c\tau^2 M_{good} < \frac{1}{2K}.
\end{align}
As a result,
\begin{align*}
    |\widetilde{\beta}_j(z)| \leq \frac{1}{1 - 1/2K} = \frac{2K}{2K - 1}.
\end{align*}

On $\mathcal{G}_0^c$, which is a \textit{small} event, using the fact that
\begin{align*}
    |\widetilde{\beta}_j(z)| \leq \frac{|z|}{v} \leq |z|n^{1+\alpha},
\end{align*}
and $|z|$ is bounded on $\mathcal{C}_n$, we therefore have
\begin{align*}
    \mathbb{E}|\widetilde{\beta}_j(z)|^k \leq C_k.
\end{align*}
\textbf{Proof of 4:} Conditioning on $Q_j$, we first use Lemma \ref{lemma:Burkholder} and then part (1) to get
    \begin{align*}
        \mathbb{E}|\epsilon_j(z)|^k &= \mathbb{E}\bigg(\mathbb{E}(|\epsilon_j(z)|^k\bigg|Q_j)\bigg)\\
        &\leq  \mathbb{E} \bigg(\frac{M_4n^\frac{k}{2}\|Q_j\|^k + M_4n(\delta_n\sqrt{n})^{2k-4}\|Q_j\|^k}{n^k}\bigg) = O(\delta_n^{2k-4}n^{-1}),
    \end{align*}
    where we used the fact that for $k > 2$, the second term starts to dominate the first since
$$n\delta_n^4 \longrightarrow \infty \iff n^{-\frac{k}{2}} \ll n^{-1}\delta_n^{2k-4}.$$
The proof for $\gamma_j$ is similar.

If the innovations are (truncated) Guassian, then we can replace the term $M_4(\delta_n\sqrt{n})^{2k-4}$ by a constant $M_{2k}>0$ and observe that the first term begins to dominate the second as  a result of which, we get
    $$\mathbb{E}|\epsilon_j(z)| = O(n^{-\frac{k}{2}}).$$
\end{proof}

\subsection{Proof of Lemma \ref{lemma:bound_centered_trace_resolvent_squared}}
\begin{proof}
We have 
\begin{align*}
    \operatorname{tr}(MQ) - \mathbb{E}\operatorname{tr}(MQ)
    &= \sum_{j=1}^n(\mathbb{E}_j - \mathbb{E}_{j-1})\operatorname{tr}M(Q - Q_j)\\
    &= \sum_{j=1}^n\Delta_j\beta_j\underbrace{x_j^*Q_jMQ_jx_j}_{:=q_j}\\
    &= \sum_{j=1}^n\Delta_j(\beta_jq_j - \widetilde{\beta_j} \underbrace{\bigl(n^{-1}\operatorname{tr}(MQ_j^2)\bigr)}_{:=\mu_j}\\
    &= \sum_{j=1}^n \Delta_j\beta_j(q_j-\mu_j) + \sum_{j=1}^n \Delta_j(\beta_j - \widetilde{\beta_j})\mu_j.
\end{align*}
Note that $q_j-\mu_j$ is an object similar to $\gamma_j$ of Lemma \ref{lemma:bounds_tightness_section} and thus 
$$\mathbb{E}|q_j - \mu_j|^4 \leq \frac{C}{n^2}.$$
Now observe that  
$\beta_j - \widetilde{\beta_j} = -\beta_j\widetilde{\beta_j}\epsilon_j,$
which leads to
\begin{align*}
    \mathbb{E}|\operatorname{tr}(MQ_j) - \mathbb{E}\operatorname{tr}(MQ_j)|^2
    \leq C_1\sum_{j=1}^n\mathbb{E}|\beta_j(q_j-\mu_j)|^2 + C_2\sum_{j=1}^n\mathbb{E}|\beta_j\widetilde{\beta_j}\mu_j\epsilon_j|^2.
\end{align*}
Note that $|\mu_j|$ has bounded expectation by Lemma \ref{lemma:bounds_tightness_section}. By similar arguments, we also obtain the following bound:
\begin{align*}
    \mathbb{E}|\beta_j\widetilde{\beta}_j\mu_j\epsilon_j|^2 = O(n^{-1}),
\end{align*}
and this concludes the proof.
\end{proof}

\subsection{Proof of Lemma \ref{lemma:hn_Ehn_hn0_uniform_results}}
\begin{proof}
By Lemma \ref{lemma:bound_centered_trace_resolvent_squared}, we have
\begin{align*}
    \sup_{z\in\mathcal{C}_n}\mathbb{E}|h_{rn}(z) - \mathbb{E}h_{rn}(z)|^2 = \sup_{z\in\mathcal{C}_n}\frac{1}{p^2}\mathbb{E}|\operatorname{tr}(A_rQ(z)) - \mathbb{E}\operatorname{tr}(A_rQ(z))|^2 = O(n^{-2}).
\end{align*}
For the second result, denote $F_n^0 := F^{c_n, H_n, G_n}$. By Theorem 5.2 of \cite{JHDP26}, we have 
$$(s_n^0, \textbf{h}_n^0, \textbf{g}_n^0) \longrightarrow (s_F, \textbf{h}, \textbf{g}),$$
for each $z \in \mathbb{C}_+$. Observing that for large $n$, 
$$\operatorname{supp}(F_n^0), \operatorname{supp}(F) \subset [e_l, e_r],$$
and the (compact) contour $\mathcal{C}$ does not intersect $[e_l, e_r]$, $\sup_{z\in \mathcal{C}} \|h_n^0(z)\|$ is therefore uniformly bounded. By Vitali's theorem, we therefore have
$$\sup_{z\in \mathcal{C}}\|\textbf{h}_n^0(z) - \textbf{h}(z)\| \to 0.$$
\end{proof}

\subsection{Proof of Lemma \ref{lemma:exp_betaj_psij_bj_difference}}
\begin{proof}
    We first simplify $\beta_j - b_j$ using (\ref{beta_b_relationship}) and (\ref{defining_rhoj}) as follows:
\begin{align*}
    \mathbb{E}(\beta_j(z) - b_j(z)) &= b_j^3(z)\mathbb{E}\rho_j^2(z) - b_j^3(z)\mathbb{E}\beta_j(z)\rho_j^3(z)\\
    &=b_j^3(z)\mathbb{E}\epsilon_j^2(z) - b_j^3(z)\mathbb{E}\beta_j(z)\rho_j^3(z) + \frac{1}{n^2}b_j^3(z)\mathbb{E}\bigl(\operatorname{tr}(Q_j\textbf{A}_j)-\mathbb{E}\operatorname{tr}(Q_j\textbf{A}_j)\bigr)^2\\
    &:= \mathcal{H}_1 + \mathcal{H}_2 + \mathcal{H}_3.
\end{align*}
By Lemma \ref{lemma:bounds_tightness_section}, we obtain
$$|\mathcal{H}_2| \leq C\bigg(\mathbb{E}|\beta_j(z)|^2\bigg)^\frac{1}{2}\bigg(\mathbb{E}|\rho_j(z)|^6\bigg)^\frac{1}{2} = Cn^{-\frac{3}{2}} = o(n^{-1}).$$

By Lemma \ref{lemma:bound_centered_trace_resolvent_squared}, we have 
$$|\mathcal{H}_3| \leq Cn^{-2} = o(n^{-1}).$$
Combining these two bounds, we have
\begin{align}\label{betaj_minus_bj}
    \mathbb{E}(\beta_j(z) - b_j(z)) = b_j^3(z)\mathbb{E}\epsilon_j^2(z) + o(n^{-1}).
\end{align}

Again using (\ref{beta_b_relationship}), we find that
\begin{align*}
    b_j(z) - \psi_j(z) &= \frac{1}{n}b_j\psi_j\mathbb{E}\operatorname{tr}(Q- Q_j)\textbf{A}_j\\
    &= -\frac{1}{n}b_j\psi_j\mathbb{E}\beta_jx_j^*Q_j\textbf{A}_jQ_jx_j\\
    &= -\frac{1}{n}b_j^2\psi_j\mathbb{E}x_j^*Q_j\textbf{A}_jQ_jx_j + \frac{1}{n}b_j^2\psi_j\mathbb{E}\beta_j\rho_jx_j^*Q_j\textbf{A}_jQ_jx_j.
\end{align*}
Observe that 
$$\frac{1}{n}\bigl|b_j^2\psi_j\mathbb{E}\beta_j\rho_jx_j^*Q_j\textbf{A}_jQ_jx_j\bigr| \leq \frac{1}{n}\mathbb{E}\bigl(|\beta_j|^2\bigr)^\frac{1}{2}\bigl(\mathbb{E}|\rho_jx_j^*Q_j\textbf{A}_jQ_jx_j|^2\bigr)^\frac{1}{2}\leq Cn^{-3/2}.$$
This implies that
\begin{align}\label{bj_minus_psij}
    b_j(z) - \psi_j(z) &= -\frac{1}{n}b_j^2\psi_j\mathbb{E}x_j^*Q_j\textbf{A}_jQ_jx_j + o(n^{-1})\\\notag
    &= -\frac{1}{n^2}b_j^2\psi_j\mathbb{E}\operatorname{tr}Q_j\textbf{A}_jQ_j\textbf{A}_j + o(n^{-1}).
\end{align}
Using the fact that $|\operatorname{tr}A| \leq n\|A\|$, we get the first result as follows:
\begin{align}\label{bj_minus_psij_1}
    |b_j(z) - \psi_j(z)| \leq \frac{C}{n} \mathbb{E}\|Q_j\textbf{A}_jQ_j\textbf{A}_j\| + o(n^{-1}) \leq \frac{C}{n}.
\end{align}

Combining (\ref{betaj_minus_bj}) and (\ref{bj_minus_psij}), we get
\begin{align}\label{exp_beta_j_minus_psi_j}
    \mathbb{E}(\beta_j(z) - \psi_j(z)) = b_j^3(z)\mathbb{E}\epsilon_j^2(z) -\frac{1}{n^2}b_j^2\psi_j\mathbb{E}\operatorname{tr}Q_j\textbf{A}_jQ_j\textbf{A}_j + o(n^{-1}).
\end{align}

By Lemma \ref{lemma:bounds_tightness_section}, $|b_j(z)|, |\psi_j(z)|$ are uniformly bounded across $j \in [n]$, $z \in \mathcal{C}_n$ and also
$$\mathbb{E}|\epsilon_j(z)|^2 = O(n^{-1}).$$
This takes care of the first term in (\ref{exp_beta_j_minus_psi_j}). The proof follows from using (\ref{tr_ABCD}), Lemma \ref{lemma:bounds_bias_section} to obtain the following bound for the second term in (\ref{exp_beta_j_minus_psi_j}):
$$\mathbb{E}\operatorname{tr}(Q_j\textbf{A}_jQ_j\textbf{A}_j) \leq \mathbb{E}\|Q_j\|^2\operatorname{tr}(\textbf{A}_j^2) \leq Cn\mathbb{E}\|Q_j\|^2 \leq Cn.$$
\end{proof}

\subsection{Proof of Lemma \ref{lemma:bounds_bias_section}}
\begin{proof}
\textbf{Proof of (1):} 
From (\ref{primal_dual_relationship}), we have 
$$-\frac{1}{nz}\sum_{j=1}^n \beta_j(z) = \ubar{s}_n(z).$$
Now observe that
\begin{align*}
     &c_n\sum_{r=1}^K  |\mathbb{E}h_{rn}(z) - h_{rn}^0(z)| \bigg|\frac{1}{n}\sum_{j=1}^n \theta_{rj}\psi_j(z)b_j^0(z)\bigg|\\
     =& \, \bigg|\frac{1}{n}\sum_{j=1}^n \bigl(\psi_j(z) - b_j^0(z)\bigr)\bigg| \\
     \leq&\, \bigg|\frac{1}{n}\sum_{j=1}^n (\psi_j(z) - \mathbb{E}\beta_j(z))\bigg| + \bigg|\frac{1}{n}\sum_{j=1}^n\mathbb{E}\beta_j(z) + z\ubar{s}(z)\bigg| + \bigg|-z\ubar{s}(z) - \frac{1}{n}\sum_{j=1}^nb_j^0(z)\bigg|\\
     \leq&\, \frac{1}{n}\sum_{j=1}^n |\mathbb{E}\beta_j(z) - \psi_j(z)| + |z|\,|\mathbb{E}\ubar{s}_n(z) - \ubar{s}(z)| + |z|\,|\ubar{s}_n^0(z) - \ubar{s}(z)|
\end{align*}
By Lemma \ref{lemma:exp_betaj_psij_bj_difference}, the first term goes to $0$. Similar to the proof of (4.1) and (4.2) in \cite{BaiSilv04}, we also have
$$\sup_{n; z \in \mathcal{C}_n} |\mathbb{E}\ubar{s}_n(z) - \ubar{s}(z)| \longrightarrow 0, \quad \sup_{n; z \in \mathcal{C}} |\ubar{s}_n^0(z) - \ubar{s}(z)| \longrightarrow 0.$$
This takes care of the second and third terms. Therefore, it suffices to show that
$$\bigg|\frac{1}{n}\sum_{j=1}^n\theta_{rj}\psi_j(z)b_j^0(z)\bigg|$$
is lower bounded for $r \in [K]$. Suppose not, then there exists a sequence $\{z_l\} \subset \mathcal{C}_n$ such that $z_l \rightarrow z_0$ and 
$$\int \frac{\theta_r dG(\boldsymbol{\theta})}{(1 + c\boldsymbol{\theta}^\top \textbf{h}(z_l))^2} \longrightarrow 0.$$

We have,
\begin{align*}
    & -z_l\ubar{s}(z_l) = \int \frac{dG(\boldsymbol\theta)}{1 + c\boldsymbol{\theta}^\top \textbf{h}(z_l)}\\
    \implies& (z_l\ubar{s}(z_l))^\prime = \sum_{r=1}^K h_r^\prime(z_l)\int \frac{\theta_r dG(\boldsymbol{\theta})}{(1 + c\boldsymbol{\theta}^\top \textbf{h}(z_l))^2}.
\end{align*}
By the continuity of $h_r(z)$, we have 
$$(z_l\ubar{s}(z_l))^\prime \longrightarrow (z_0\ubar{s}(z_0))^\prime = 0.$$
Then from $s(z_0) + z_0s^\prime(z_0) = 0$, we obtain another solution of $s(z)$ that is independent of $c$, $H$ and $G$ which contradicts the uniqueness of the solution of the Matrix Dyson system. Therefore for any $r \in [K]$, we must have
$$\bigg|\int \frac{\theta_r dG(\boldsymbol{\theta})}{(1 + c\boldsymbol{\theta}^\top \textbf{h}(z))^2}\bigg| > 0.$$
By the continuity and convergence of $\mathbb{E}\textbf{h}_n(z)$ to \textbf{h}(z) and Lemma \ref{lemma:hn_Ehn_hn0_uniform_results}, we use DCT to infer that
$$\bigg|\frac{1}{n}\sum_{j=1}^n \frac{\theta_{rj}}{(1 + c_n\boldsymbol{\theta}_j^\top \mathbb{E}\textbf{h}_n(z))(1 + c\boldsymbol{\theta}_j^\top \textbf{h}_n^0(z))}\bigg| = \bigg|\frac{1}{n}\sum_{j=1}^n \theta_{rj}\psi_j(z)b_j^0(z)\bigg| > 0.$$
This implies that for any $r \in [K]$, we have
$$\sup_{n;z\in \mathcal{C}_n}|\mathbb{E}h_{rn}(z) - h_{rn}^0(z)| \rightarrow 0.$$

\textbf{Proof of (2):} 
We have
$$b_j^0(z) = \frac{1}{1 + c_n\boldsymbol{\theta}_j^\top \textbf{h}_n^0(z)}.$$ 
It suffices to show that the denominator is lower bounded. From (\ref{sup_exp_hrn_bound}), we have
$$\sup_{n; z\in \mathcal{C}_n} |\mathbb{E}h_{rn}(z)| < \frac{3}{8K\tau}.$$
Therefore from part (1), for any $r \in [K]$, we also have 
$$\sup_{n; z\in \mathcal{C}_n} |h_{rn}^0(z)| < \frac{3}{8K\tau}.$$
The proof for the bound on $|b_j^0(z)|$ is exactly similar to that of part (2) of Lemma \ref{lemma:bounds_tightness_section}. The second part of (2) is immediate upon observing the following:
\begin{align*}
    |\psi_j(z) - b_j^0(z)| = |\psi_j(z)| \, |b_j^0(z)| \, \boldsymbol{\theta}_j^\top (\mathbb{E}\textbf{h}_n(z) - \textbf{h}_n^0(z))|.
\end{align*}
\end{proof}

\subsection{\texorpdfstring{Proof of (\ref{vanishing_contribution_vertical_arms}).}{}}\label{proof_vanishing_contribution_vertical_arms}
\begin{proof}
Let $\mathcal{C}_v := \mathcal{C}_l \cup \mathcal{C}_r$ denote the vertical arms of the contour. By Cauchy Schwarz and Fubini, we have
\begin{align*}
    \mathbb{E}\bigg|\int_{\mathcal{C}_v} f(z)M_{n1}(z)dz\bigg|^2 \leq C\mathbb{E}\int_{\mathcal{C}_v} |M_{n1}(z)|^2dz = C\int_{\mathcal{C}_v} \mathbb{E}|M_{n1}(z)|^2dz.
\end{align*}
Therefore, it suffices to show that
$$\sup_{n; z \in \mathcal{C}_v}\mathbb{E}|M_{n1}(z)|^2 = O(1).$$
Fix $z_0 \in \mathcal{C}_u$. From (\ref{goal_tightness}), for any $z \in \mathcal{C}_v$, we have
\begin{align*}
    \mathbb{E}|M_{n1}(z)|^2 &\leq \mathbb{E}|M_{n1}(z_0)|^2 + \bigg(\sup_{z_1,z_2 \in \mathcal{C}_n} \frac{\mathbb{E}|M_{n1}(z_1) - M_{n1}(z_2)|^2}{|z_1 - z_2|^2} \bigg)|z - z_0|^2\\
    &\leq \mathbb{E}|\operatorname{tr}Q(z_0) - \mathbb{E}\operatorname{tr}Q(z_0)|^2 + C\operatorname{diam}^2(\mathcal{C}_n) < \infty,
\end{align*}
where the last inequality follows from Lemma \ref{lemma:bound_centered_trace_resolvent_squared}.
\end{proof}

\subsection{Proof of Claim \ref{claim:standardized_vector}}
\begin{proof}
    Define $y := \textbf{U}^+x + (I - \textbf{U}^+\textbf{U})w$, where $\textbf{U}^+$ is the Moore-Penrose Inverse of $\textbf{U}=\textbf{A}^\frac{1}{2}$ and $w$ is any random vector independent of $x$ with $\mathbb{E}(w) = \textbf{0}$ and $\operatorname{Var}(w) = I$. Then, $\textbf{U} y = x$ and $\mathbb{E}(y) =\textbf{0}$ holds trivially. Since $\textbf{UU}^+ = (\textbf{U}^+\textbf{U})^* = \textbf{U}^+\textbf{U}$ by the properties of Moore-Penrose inverses, we therefore have
\begin{align*}
    \operatorname{var}(y) &= \textbf{U}^+\textbf{AU}^+ + (I - \textbf{U}^+\textbf{U})^2\\
    &= (\textbf{U}^+\textbf{U})(\textbf{UU}^+) + I - 2\textbf{U}^+\textbf{U} + \textbf{U}^+\textbf{UU}^+\textbf{U} \\
    &= \textbf{U}^+\textbf{U}\textbf{U}^+\textbf{U} + I - 2\textbf{U}^+\textbf{U} + \textbf{U}^+\textbf{U} = I.
\end{align*}
\begin{remark}
     Note that by choosing $w_i$ to be $O(\delta_n\sqrt{n})$, we can also ensure $y_i=O(\delta_n\sqrt{n})$. 
\end{remark}
\end{proof}

\subsection{Proof of Lemma \ref{lemma:tracial_bounds}}\label{sec:proofLemma4.2}
\begin{proof}
\textbf{Bound on $\operatorname{tr}(C_j(z)M)$:}
Let $M$ be a non-random matrix. Define $H_{ij} = \textbf{A}_iQ_{ij}MR_j$. Using Lemma \ref{lemma:Bounds_Important_Quantities}, we have
\begin{align*}
    \mathbb{E}|\operatorname{tr}C_j(z)M| &= \mathbb{E}\absmod{\sum_{i\neq j} \psi_i \bigg(x_{i}^*Q_{ij}MR_jx_{i}-n^{-1}\operatorname{tr}(\textbf{A}_iQ_{ij}MR_j)\bigg)}\\
    &\leq C\sum_{i\neq j} \bigg(\mathbb{E}|x_{i}^*Q_{ij}MR_jx_{i}-n^{-1}\operatorname{tr}H_{ij}^{(r)}|^2\bigg)^\frac{1}{2}= O(n^\frac{1}{2}).
\end{align*}

\textbf{Bound on $\operatorname{tr}(D_j(z)M)$:}
Let $M$ be any (possibly random) matrix with non-random bound on its spectral norm. Then,
\begin{align*}
    \mathbb{E}|\operatorname{tr}D_j(z)M| 
    &\leq \sum_{i \neq j} \mathbb{E} |(\beta_{ij}(z) - \psi_i(z))x_i^*Q_{ij}MR_jx_i|\\
    &\leq \sum_{i\neq j} \bigg(\mathbb{E}|\beta_{ij}(z) - \psi_i(z)|^2\bigg)^\frac{1}{2} \bigg(\mathbb{E}|x_i^*Q_{ij}MR_jx_i|^2\bigg)^\frac{1}{2}.
\end{align*}
By Lemma \ref{lemma:bounds_tightness_section}, we have 
$$\mathbb{E}|x_i^*Q_{ij}MR_jx_i|^2 \leq 2\mathbb{E}|x_i^*Q_{ij}MR_jx_i - n^{-1}\operatorname{tr}H_{ij}|^2 + 2\mathbb{E}|n^{-1}\operatorname{tr}H_{ij}|^2 \leq C.$$
Together with Lemma \ref{lemma:betaij_minus_psii}, we have $$\mathbb{E}|\operatorname{tr}D_j(z)M| = O(n^\frac{1}{2}).$$

\textbf{Bound on $\operatorname{tr}(E_j(z)M)$:} By Lemma \ref{lemma:Bounds_Important_Quantities}, we have
\begin{align*}
    &E_j(z) := \frac{1}{n}\sum_{i\neq j}\psi_i(z)R_j\textbf{A}_i(Q_{ij} - Q_j)\\
    \implies&|\mathbb{E}\operatorname{tr}E_j(z)M| \leq \frac{1}{n} \sum_{i\neq j} \frac{|z| \times \|M\|}{v_0^3} = O(1). 
\end{align*}
\end{proof}

\begin{lemma}\label{lemma:betaij_minus_psii}
    Recall $Q_{ij}, \beta_{ij}, b_{ij}$ from (\ref{essential_quantities1}). Further define
    \begin{align}
        \widetilde{\beta}_{ij} = (1 + n^{-1}\operatorname{tr}\textbf{A}_iQ_{ij})^{-1}, \quad \epsilon_{ij} = x_{i}^*Q_{ij}x_{i} - n^{-1}\operatorname{tr}\textbf{A}_iQ_{ij}.
    \end{align}
Then, we have
    $$\sup_{z\in\mathcal{C}_n}\mathbb{E}|\beta_{ij}(z) - \psi_i(z)|^2 = O(n^{-1}).$$
\end{lemma}

\begin{proof}
We have 
\begin{align*}
    |\beta_{ij}(z) - \psi_i(z)| \leq |\beta_{ij}(z) - \widetilde{\beta}_{ij}(z)| + |\widetilde{\beta}_{ij}(z) - b_{ij}(z)| + |b_{ij}(z) - b_i(z)| + |b_{i}(z) - \psi_i(z)|.
\end{align*}
Similar to Lemma \ref{lemma:bounds_tightness_section}, we can establish that $|\beta_{ij}|,|\widetilde{\beta}_{ij}|$ have uniformly bounded moments and $|b_{ij}(z)|$ is uniformly bounded over $\mathcal{C}_n$ and $\epsilon_{ij}(z)$ behaves similarly to $\epsilon_j(z)$ uniformly over $\mathcal{C}_n$. As a result,
    $$\sup_{z\in\mathcal{C}_n}\mathbb{E}|\beta_{ij}(z) - \widetilde{\beta}_{ij}(z)|^2 = O(n^{-1}).$$
    By Lemma \ref{lemma:bound_centered_trace_resolvent_squared}, we have
    \begin{align*}
        \sup_{z\in\mathcal{C}_n}\mathbb{E}|\widetilde{\beta}_{ij}(z) - b_{ij}(z)| \leq \frac{C}{n}\mathbb{E}(|\widetilde{\beta}_{ij}(z)|^2)^\frac{1}{2}\bigg(\mathbb{E}|\operatorname{tr}\textbf{A}_iQ_{ij}(z) - \mathbb{E}\operatorname{tr}\textbf{A}_iQ_{ij}(z)|^2\bigg)^\frac{1}{2} = O(n^{-1}).
    \end{align*}
Now 
\begin{align*}
\sup_{z\in\mathcal{C}_n}|b_{ij}(z) - b_i(z)| &\leq \frac{C}{n}|\mathbb{E}\operatorname{tr}\textbf{A}_i(Q_i(z) - Q_{ij}(z))|= \frac{C}{n}\mathbb{E}|\beta_{ij}x_j^*Q_{ij}\textbf{A}_iQ_{ij}x_j| = O(n^{-1}).
\end{align*}
From (\ref{bj_minus_psij_1}), we have $|b_i(z) - \psi_i(z)| = O(n^{-1}).$ Combining everything, we have the result.
\end{proof}

\begin{lemma}\label{lemma:Uj_U0_U_element_wise_difference}
    $|U_j[r,s] - U^0[r,s]| = o(1), \quad |U^0[r,s] - U_{2,2}[r,s]| = o(1).$
\end{lemma}
\begin{proof}
\begin{align*}
    U_j[r,s] - U^0[r,s] &= \int \frac{\lambda_r\lambda_sdH_n(\boldsymbol{\lambda})}{\prod_{l=1}^2(-z + \frac{1}{n}\sum_{i\neq j}\psi_i(z_l)\boldsymbol{\lambda}^\top \boldsymbol{\theta}_j)} - \int \frac{\lambda_r\lambda_sdH_n(\boldsymbol{\lambda})}{\prod_{l=1}^2(-z + \frac{1}{n}\sum_{i\neq j}b_i^0(z_l)\boldsymbol{\lambda}^\top \boldsymbol{\theta}_j)} .
\end{align*}

By Lemma \ref{lemma:bounds_bias_section} we have
\begin{align*}
    \sup_{z\in\mathcal{C}_n}\absmod{\frac{1}{n}\sum_{i\neq j}\psi_i(z)\boldsymbol{\lambda}^\top \boldsymbol{\theta}_j) - \frac{1}{n}\sum_{i\neq j}b_i^0(z)\boldsymbol{\lambda}^\top \boldsymbol{\theta}_j)} \longrightarrow 0.
\end{align*}
Since the denominators are uniformly bounded away from 0, we use DCT to conclude. The proof of the second part follows similarly by applying Lemma \ref{lemma:hn_Ehn_hn0_uniform_results}.
\end{proof}

\section{Results related to the expansion of \texorpdfstring{$I - cU_{2,2}W_2$}{}}\label{sec:invertibility_statbility_operator}
\textbf{Note on notation:} Throughout this section,  $\rho(\cdot)$ denotes the spectral radius, $\|\cdot\|$ denotes the spectral norm, $\circ$ denotes the Hadamard product, and $\sqrt{M}$ denotes the entry-wise square root of a matrix $M$ with non-negative entries. For vectors/ matrices \textbf{x} and \textbf{y}, we will use the notation $\textbf{x} \prec \textbf{y}$ to indicate that $x_i < |y_i|$. 
\begin{lemma}\label{lemma:spectral_radius_le_1}
     Let $\rho(M)$ denote the spectral radius of the matrix $M$. Recalling the definitions of $U_{2,2}(\cdot, \cdot)$ and $W_2(\cdot , \cdot)$ from (\ref{sec:Glossary}), define $M(z_1, z_2) = c\,U_{2,2}(z_1, z_2)\,W_2(z_1, z_2).$  For any $z_1, z_2 \in \mathbb{C}\backslash \mathbb{R}$, the following holds: $$\rho(M(z_1, z_2)) <1.$$\end{lemma}
\begin{proof}
\textbf{Note}: Throughout this proof, we will denote $U_{2,2}$ as $U$ and $W_2$ as $W$ for simplicity.

We will use a contraction argument. For every $z \in \mathbb{C}_+$, let $(\textbf{h}(z), \textbf{g}(z))$ represent the unique solution of the matrix Dyson system. Recall that $\textbf{h}(\bar{z}) = \overline{\textbf{h}(z)}$, $\textbf{g}(\bar{z}) = \overline{\textbf{g}(z)}$. Denoting $v = \Im(z)$, we have
\begin{align*}
    \Im(\textbf{h}(z)) &= v\int \frac{\boldsymbol{\lambda}dH(\boldsymbol{\lambda})}{|z(1 + \boldsymbol{\lambda}^\top\textbf{g}(z))|^2} + \int \frac{\boldsymbol{\lambda\lambda^\top} dH(\boldsymbol{\lambda})}{|z(1 + \boldsymbol{\lambda}^\top\textbf{g}(z))|^2}\Im(z\textbf{g}(z))\\
    &= v\int \frac{\boldsymbol{\lambda}dH(\boldsymbol{\lambda})}{|z(1 + \boldsymbol{\lambda}^\top\textbf{g}(z))|^2} + U(z,\bar{z})\Im(z\textbf{g}(z)).
\end{align*}
Similarly, 
\begin{align*}
    &\Im(z\textbf{g}(z)) =  c\int \frac{\boldsymbol{\theta\theta^\top} dG(\boldsymbol{\theta})}{|1 + c\boldsymbol{\theta}^\top\textbf{h}(z)|^2}\Im(\textbf{h}(z)) = cW(z,\bar{z})\Im(\textbf{h}(z)).
\end{align*}
Combining the two, we have
\begin{align*}    
\Im(\textbf{h}(z)) &= v\int \frac{\boldsymbol{\lambda}dH(\boldsymbol{\lambda})}{|z(1 + \boldsymbol{\lambda}^\top\textbf{g}(z))|^2} + cU(z,\bar{z})W(z, \bar{z})\Im(\textbf{h(z)})\\
&= v\int \frac{\boldsymbol{\lambda}dH(\boldsymbol{\lambda})}{|z(1 + \boldsymbol{\lambda}^\top\textbf{g}(z))|^2} + M(z, \bar{z})\Im(\textbf{h(z)}).
\end{align*}
We will show that the first term on the RHS is bounded away from $0$ element wise by a quantity that depends on $K, z, \tau$ and  $H$. Observe that for any $r \in [K]$, $\lambda_r \in [0,\tau]$ for any $\boldsymbol{\lambda} \in \operatorname{supp}(H)$ and $|g_r(z)| \leq \tau/v$. As a result,
\begin{align*}
   |1 + \boldsymbol{\lambda}^\top \textbf{g}(z)| \leq 1 + |\boldsymbol{\lambda}^\top \textbf{g}(z)| \leq 1 + \frac{K\tau^2}{v}.
\end{align*}
Since $H$ is not degenerate at $\textbf{0}$ across all marginals, there exists a vector $C_H \in \mathbb{R}_+^K$ such that $$\int \boldsymbol{\lambda} dH(\boldsymbol{\lambda}) \succeq C_H \succ \textbf{0}.$$
From these inequalities, we get
\begin{align}\label{Imh_contraction_bounded_away}
&\Im(\textbf{h}(z)) \succeq v\int \frac{\boldsymbol{\lambda} dH(\boldsymbol{\lambda})}{|z|^2(1 + K\tau^2/v)^2} + M(z, \bar{z})\Im(\textbf{h(z)})\\ \notag
\implies& (I - M(z,\bar{z}))\Im(\textbf{h}(z)) \succeq \frac{vC_H}{|z|^2(1 + K\tau^2/v)^2} \succ \textbf{0}.
\end{align}
Note that $M(z, \bar{z})$ is entry-wise non-negative since $U(z, \bar{z})$ and $W(z, \bar{z})$ are entry-wise non-negative (and p.s.d.) matrices. Therefore, by Theorem \ref{theorem:CollatzWielandt}, for any $z \in \mathbb{C}_+$, we have
\begin{align}\label{Collatz_Wielandt}
    \rho(M(z, \bar{z})) < 1.
\end{align}
Since $M(z, \bar{z}) = M(\bar{z}, z)$, the result extends to all $z \in \mathbb{C}\backslash\mathbb{R}.$ Now let $z_1, z_2 \in \mathbb{C}\backslash\mathbb{R}$ be arbitrary. Let $(\textbf{h}_i, \textbf{g}_i) = (\textbf{h}(z_i), \textbf{g}(z_i))$ represent the unique fixed point solutions at $z_i, i = 1,2$. By Cauchy Schwarz, for any $r,s, \in [K]$, we get
\begin{align}\label{element_wise_CS_U}
  |U_{rs}(z_1, z_2)|
  &=\absmod{\int \frac{\lambda_r \lambda_sdH(\boldsymbol{\lambda})}{z_1z_2(1+\boldsymbol{\lambda^\top \textbf{g}_1})(1+\boldsymbol{\lambda^\top \textbf{g}_2})}}\\ \notag
  &\leq \bigg(\int \frac{\lambda_r\lambda_sdH(\boldsymbol{\lambda})}{|z_1(1+\boldsymbol{\lambda^\top \textbf{g}_1})|^2}\bigg)^\frac{1}{2}\bigg(\int \frac{\lambda_r\lambda_sdH(\boldsymbol{\lambda})}{|z_2(1+\boldsymbol{\lambda^\top \textbf{g}_2})|^2}\bigg)^\frac{1}{2}\\\notag
  &= \sqrt{U_{rs}(z_1, \bar{z_1})}\sqrt{U_{rs}(z_2, \bar{z_2})}\\ \notag
  \implies |U(z_1, z_2)| &\preceq \sqrt{U(z_1, \bar{z}_1)} \circ \sqrt{U(z_2, \bar{z}_2)}, 
\end{align}
where $A \circ B$ denotes the Hadamard product of equi-dimensional matrices $A$ and $B$ and $\sqrt{A}$ denotes the entry-wise square root for the positive matrix $A$. Similarly, we also have 
\begin{align}\label{element_wise_CS_W}
    |W(z_1, z_2)| \prec \sqrt{W(z_1, \bar{z}_1)} \circ \sqrt{W(z_2, \bar{z}_2)}.
\end{align}  
Noting that $$M_{rs}(z_1, z_2) = c(U(z_1, z_2)W(z_1, z_2))_{rs} = \sum_{k=1}^KcU_{rk}(z_1, z_2)W_{ks}(z_1, z_2),$$ we use (\ref{element_wise_CS_U}), (\ref{element_wise_CS_W}) and invoke Cauchy Schwarz to get
\begin{align*}
    |M_{rs}(z_1, z_2)| &\leq \sum_{k=1}^Kc\sqrt{U_{rk}(z_1, \bar{z}_1)U_{rk}(z_2, \bar{z}_2) W_{ks}(z_1, \bar{z}_1)W_{ks}(z_2, \bar{z}_2)}\\
    &\leq \bigg(\sum_{k=1}^K cU_{rk}(z_1, \bar{z}_1)W_{ks}(z_1, \bar{z}_1)\bigg)^\frac{1}{2}\bigg(\sum_{k=1}^K cU_{rk}(z_2, \bar{z}_2)W_{ks}(z_2, \bar{z}_2)\bigg)^\frac{1}{2}\\
    &= \sqrt{M_{rs}(z_1, \bar{z}_1)} \circ \sqrt{M_{rs}(z_2, \bar{z}_2)}.
\end{align*}
Therefore, we have 
\begin{align}\label{element_wise_CS_rho_M}
    M(z_1, z_2) \preceq |M(z_1, z_2)| \preceq \sqrt{M(z_1, \bar{z}_1)}\circ \sqrt{M(z_2, \bar{z}_2)}.
\end{align}
The log-convexity of the Perron root (see Lemma \ref{lemma:log_convexity}), the monotonicity of $\rho(\cdot)$, (\ref{element_wise_CS_rho_M}) and (\ref{Collatz_Wielandt}) finally leads to
\begin{align*}
    \rho(M(z_1, z_2)) \leq \rho(|M(z_1, z_2)|) \leq \sqrt{\rho(M(z_1, \bar{z_1}))\rho(M(z_2, \bar{z}_2))} < 1.
\end{align*}
\end{proof}

\begin{remark}\label{radius_uniformly_away_from_1_on_contour}
    \textbf{Stability of the operator:} On any compact subset of $\mathbb{C}\backslash\mathbb{R}$ and in particular on the horizontal part of the contour, i.e. $\mathcal{C}_u$, we have
    $$\rho(M(z_1, z_2)) \leq \delta < 1,$$
    where $\delta$ depends on the compact subset. This follows from (\ref{Imh_contraction_bounded_away}) and the Collatz-Wielandt characterization of spectral radius. Let $\textbf{x} = \Im(\textbf{h}(z))$. We have
    \begin{align*}
        &\rho(M) \leq \max_{r \in [K]} \frac{(M\textbf{x})_r}{x_r}\\
        \implies& 1 - \rho(M) \geq \min_{r \in [K]} \frac{v_0C_H(r)}{|z|^2(1 + K\tau^2/v_0)^2x_r} \geq \frac{v_0^2C_H(r)}{\tau |z|^2(1 + K\tau^2/v_0)^2} > 0,
    \end{align*}
    where we used the fact that $|x_r| = |\Im(h_r(z))|\leq \tau/v_0.$
    
    We now extend the argument to the vertical arms of the contour. Observe that it suffices to establish (\ref{Collatz_Wielandt}) when $z\in \mathcal{C}, \Im(z) = 0$ as well.
    Since $M(z, \overline{z})$ is continuous for $z \in \mathcal{C}$, and $\rho(\cdot)$ is a continuous function, by compactness,
    there exists $\delta$ such that $$\sup_{z\in\mathcal{C}}\rho(M(z,\overline{z})) \leq \delta < 1.$$
    
    For the rest of this proof, we assume that $z \in {x_l, x_r}$. Choose $x_l, x_r$ as described previously in but additionally, we ensure that $$x_0 := \min\{|x_l|, |x_r|\}$$
    is very large which ensures that $|z|$ is large on the vertical arms. 
    Note that we established
    $$|1 + c_n\boldsymbol{\theta}^\top \mathbb{E}\textbf{h}_n(z)| \geq 1/4,$$
    during the proof of Lemma \ref{lemma:bounds_tightness_section}. 
    By continuity of $\mathbb{E}\textbf{h}_n(z)$ to $\textbf{h}(z)$ over $z\in\mathcal{C}_n$ and the analytic extension of $\textbf{h}(z)$ from $\mathbb{C}\backslash\mathbb{R}$ to $\mathbb{C}\backslash\operatorname{supp}(F)$ (and similar arguments for $\textbf{g}(z)$), we therefore have
    $$|1 + c\boldsymbol{\theta}^\top \textbf{h}(z)| \geq \frac{1}{4}, \quad |1  + \boldsymbol{\lambda}^\top \textbf{g}(z)| \geq \frac{1}{4}.$$
    Then for any $r,s \in [K]$, we have
    \begin{align*}
         |W_{r,s}(z, \bar{z})| \leq \frac{16\tau^2}{9} , \quad |U_{r,s}(z, \bar{z})| \leq \frac{16\tau^2}{9x_0^2}.
    \end{align*}
Since $x_0$ can be chosen to be arbitrarily large, in particular we can ensure that quantity
\begin{align*}
    \max_{r \in [K]} \sum_{s,t=1}^K c|U_{r,t}(z, \bar{z})|\,|W_{t,s}(z, \bar{z})| 
\end{align*}
is sufficiently small and in particular, smaller than $1$. By the Perron Frobenius theorem, we finally have the result.
\end{remark}

\begin{lemma}\label{lemma:log_convexity}
    For non-negative matrices $A, B$, we have
$$ \rho\bigl(\sqrt{A} \circ \sqrt{B}\bigr) \leq \rho(AB).$$
\end{lemma}

\begin{proof}
Let $N := \sqrt{A} \circ \sqrt{B}$, so that $N_{rs}=\sqrt{A_{rs}B_{rs}}\ge 0$. For arbitrary $k \in\mathbb{N}$, we  consider $(r,s)$ entry of the $k$-th power of $N$. Let $\pi=(t_0=r,t_1,\ldots,t_k=s)$ where $\{t_i\}_{i=1}^{k-1} \in [K]$ denote a \textit{path} from $r$ to $s$. Then,
$$(N^k)_{rs} = \sum_{\pi}\prod_{i=1}^{k}
\sqrt{A_{t_{i-1}t_i}B_{t_{i-1}t_i}} = \sum_{\pi}
\left(
\prod_{i=1}^{k}A_{t_{i-1}t_i}
\right)^{1/2}
\left(
\prod_{i=1}^{k}B_{t_{i-1}t_i}
\right)^{1/2}.$$
Since all terms are nonnegative, we apply Cauchy Schwarz to get
$$(N^k)_{rs}
\le
\left(
\sum_{\pi}
\prod_{i=1}^{k}
A_{t_{i-1}t_i}
\right)^{1/2}
\left(
\sum_{\pi}
\prod_{i=1}^{k}
B_{t_{i-1}t_i}
\right)^{1/2}
=
\left(
(A^k)_{rs}(B^k)_{rs}
\right)^{1/2}.$$
By basic properties of the spectral norm, we have
$$\max_{r,s}(A^k)_{rs}\le \|A^k\|, \qquad \max_{r,s}(B^k)_{rs}\le \|B^k\|.$$
Therefore, 
$$\|N^k\| \leq K\max_{r,s}(N^k)_{rs} \leq K\,\|A^k\|^{1/2}\|B^k\|^{1/2}.$$
Taking \(k\)-th roots and letting \(k\to\infty\) and using Gelfand's formula (Lemma \ref{lemma:Gelfand_formula}) for spectral radius
we have
$$ \rho(N) \leq \underset{k \rightarrow \infty}{\lim} K^{1/k} \bigl(\|A^k\|^{1/k}\bigr)^{1/2} \bigl(\|B^k\|^{1/k}\bigr)^{1/2}= \sqrt{\rho(A)\rho(B)}.$$
\end{proof}

\begin{lemma}\label{lemma:Gelfand_formula}
    Gelfand's formula: For any matrix $M$ and any norm $\|\cdot\|$, the spectral radius of $M$ is given by $$\rho(M)=\lim_{k\to\infty}\|M^k\|^{1/k}.$$
\end{lemma}

\begin{theorem}\label{theorem:CollatzWielandt}
    \textbf{Collatz–Wielandt:} Let $M$ be a nonnegative matrix. Then $$\rho(A) = \inf_{y \succ \textbf{0}} \max_i \frac{(My)_i}{y_i}.$$
Consequently, if there exists $x \succ 0$ such that $Ax \prec x$, then $\rho(A)<1$.
\end{theorem}

\begin{proof}
    If such an \textbf{x} exists, it immediately follows that
    \begin{align*}
      \rho(A) = \underset{y \succ \textbf{0}}{\inf} \underset{i}{\max}\frac{(My)_i}{y_i} \leq \frac{(Mx)_i}{x_i} < 1.
    \end{align*}
\end{proof}

\begin{lemma}\label{lemma:tracial_integration}
    Suppose $\rho(X) < 1$. Then,
    $$\int_0^1 t\operatorname{tr}\bigg((I-tX)^{-1}X\bigg)^2 dt = \operatorname{tr}\bigg((I-X)^{-1}X  + \log(I - X)\bigg).$$
\end{lemma}

\begin{proof}
    For $|t|<1$, we can write $(I - tX)^{-2} = \sum_{k=0}^\infty (k+1)t^kX^k$. Then the LHS can be expressed as
    \begin{align*}
        \int_0^1 \sum_{k=0}^\infty (k+1)t^{k+1}\operatorname{tr}(X^{k+2})dt
        =& \sum_{k=0}^\infty \operatorname{tr}(X^{k+2}) \int_0^1 (k+1)t^{k+1}dt\\
        =& \sum_{k=1}^\infty \frac{k}{k+1}\operatorname{tr}(X^{k+1})\\
        =& \sum_{k=1}^\infty \operatorname{tr}(X^k) - \sum_{k=1}^\infty \frac{1}{k}\operatorname{tr}(X^k)\\
        =&\operatorname{tr}\bigg((I - X)^{-1}X\bigg) + \operatorname{tr}\bigg(\log(I-X)\bigg).
    \end{align*}
\end{proof}

\section{Results related to the bias term}
\subsection{Proof of Lemma \ref{lemma:norm_Rz_uniform_bound}}\label{proof_Rz_uniform_bound}
\begin{proof}
Recall that 
$$R(z) = \bigg(\frac{1}{n}\sum_{j=1}^n \psi_j(z)\textbf{A}_j - zI_p\bigg)^{-1},$$
where $\textbf{A}_j = \sum_{r=1}^K \theta_{rj}A_r$ for $j \in [n]$. Note that for any $j \in [n]$, since $\Im(\mathbb{E}\textbf{h}_n(z)) \in \mathbb{R}_+^K$ by basic properties of Stieltjes transforms. As a consequence, we have
\begin{align*}
    \Im(\psi_j(z)) = \frac{c_n\boldsymbol{\theta}_j^\top \overline{\mathbb{E}\textbf{h}_n(z)}}{|1 + c_n\boldsymbol{\theta}_j^\top \mathbb{E}\textbf{h}_n(z)|^2} \leq 0.
\end{align*}
Therefore, for any $\boldsymbol{\lambda} \in \mathbb{R}_+^K$, we have
\begin{align*}
    \Im\bigg( \frac{1}{n}\sum_{j=1}^n\psi_j(z)\boldsymbol{\theta}_j^\top\boldsymbol{\lambda} - z\bigg) = -v - \frac{1}{n}\sum_{j=1}^n\Im(\psi_j(z))\boldsymbol{\theta}_j^\top\boldsymbol{\lambda} \leq -v.
\end{align*}
As a result, $\sup_{z \in \mathcal{C}_u}\|R(z)\| \leq v_0^{-1}.$ For $z 
\in \mathcal{C}_v$, we establish another lower bound. Since $\Im(zg_r(z))>0$, there exists a constant $\delta_1>0$ such that for any $\boldsymbol{\lambda} \in \operatorname{supp}(H)$, we have
$$\inf_{z \in \mathcal{C}_n} |z\boldsymbol{\lambda}^\top \textbf{g}(z) + z| \geq \delta_1.$$
This $\delta_1$, in particular, depends on the compact set $\mathcal{C}_n$ and therefore on $v_0$. By the continuity of $z\textbf{g}(z)$ on $\mathcal{C}_v$, there exists $C_0>0$ such that
$$\sup_{z \in \mathcal{C}_n} \|z\textbf{g}(z)\| < C_0.$$
Since $H_n \overset{d}{\longrightarrow} H$, for large $n$ and any $\boldsymbol{\lambda} \in \operatorname{supp}(H)$, $S_n$ has an eigenvalue $\boldsymbol{\lambda}_n^0$ such that 
$$\|\boldsymbol{\lambda}_n^0 - \boldsymbol{\lambda}\| < \frac{\delta_1}{4C_0}.$$
The last ingredient for this proof is the following upper bound. Observe that for any $r \in [K]$, we have
\begin{align*}
    &\frac{1}{n}\sum_{j=1}^n \psi_j(z)\theta_{rj} + zg_r(z)\\
    =& \int \frac{\theta_r dG_n(\boldsymbol{\theta})}{1 + c_n \boldsymbol{\theta}^\top \mathbb{E}\textbf{h}_n(z)} - \int \frac{\theta_r dG(\boldsymbol{\theta})}{1 + c\boldsymbol{\theta}^\top \textbf{h}(z)}\\
    =& \int \frac{c\theta_r\boldsymbol{\theta}^\top \textbf{h}(z) - c_n\theta_r\boldsymbol{\theta}^\top \mathbb{E}\textbf{h}_n(z)}{(1 + c_n \boldsymbol{\theta}^\top \mathbb{E}\textbf{h}_n(z))(1 + c\boldsymbol{\theta}^\top \textbf{h}(z))}dG_n(\boldsymbol{\theta})+ \int \frac{\theta_r}{1 + c\boldsymbol{\theta}^\top h(z)} d\{G_n(\boldsymbol{\theta}) - G(\boldsymbol{\theta})\}\\
    =&-c\sum_{s=1}^K (\mathbb{E}h_{sn}(z) - h_{sn}(z)) \int \frac{\theta_r\theta_s dG_n(\boldsymbol{\theta})}{(1 + c_n \boldsymbol{\theta}^\top \mathbb{E}\textbf{h}_n(z))(1 + c\boldsymbol{\theta}^\top \textbf{h}(z))}\\ 
    &-(c_n - c)\sum_{s=1}^K \mathbb{E}h_{sn}(z) \int \frac{\theta_r\theta_s dG_n(\boldsymbol{\theta})}{(1 + c_n \boldsymbol{\theta}^\top \mathbb{E}\textbf{h}_n(z))(1 + c\boldsymbol{\theta}^\top \textbf{h}(z))}\\
    &+\int \frac{\theta_r}{1 + c\boldsymbol{\theta}^\top h(z)} d\{G_n(\boldsymbol{\theta}) - G(\boldsymbol{\theta})\}
\end{align*}

First of all, the function classe defined by 
$$z \longmapsto \theta_r/(1 + c\boldsymbol{\theta}^\top \textbf{h}(z)); \boldsymbol{\theta} \in \{\boldsymbol{\theta}_j; j \in [n]\},$$
is bounded and equi-continuous for $z \in \mathcal{C}_n$. Since $G_n \xrightarrow{d} G$, the last term goes to $0$. By Lemma \ref{lemma:bounds_bias_section}, we have 
$$\sup_{n; z \in \mathcal{C}_n} |\mathbb{E}h_{rn}(z) - h_r(z)| = o(1).$$
Using DCT, this implies that the first term goes to 0. Similarly, the second term goes to $0$ as $c_n \rightarrow c$ and $\mathbb{E}|h_{rn}(z)|$ is bounded for $z \in \mathcal{C}_n$ from (\ref{sup_exp_hrn_bound}). Therefore, we have
$$\sup_{n;z \in \mathcal{C}_n} \bigg|\frac{1}{n}\sum_{j=1}^n \psi_j(z)\theta_{rj} + zg_r(z)\bigg| \longrightarrow 0,$$
and in particular bounded by $\delta_1/4K\tau$ for large $n$. Therefore, by triangle inequality, we finally have
\begin{align*}
    &\inf_{z \in \mathcal{C}_v} \bigg|\frac{1}{n}\sum_{j=1}^n \psi_j(z)\boldsymbol{\theta}_j^\top \boldsymbol{\lambda}_n^0 - z\bigg|\\ 
    \geq& \,\inf_{z \in \mathcal{C}_v} |z\boldsymbol{\lambda}^\top \textbf{g}(z) +z| - \sup_{z\in \mathcal{C}_v} \|z\textbf{g}(z)\|\, \|\boldsymbol{\lambda}_n^0 - \boldsymbol{\lambda}\| - \sup_{z \in \mathcal{C}_v} \|\boldsymbol{\lambda}_n^0\|\, \bigg\|\frac{1}{n}\sum_{j=1}^n \psi_j(z)\boldsymbol{\theta}_j + z\textbf{g}(z)\bigg\|\\ 
    \geq& \,\delta_1 - C_0\frac{\delta_1}{4C_0} -  K\tau \frac{\delta_1}{4K\tau} = \frac{\delta_1}{2}.
\end{align*}
This concludes the proof.
\end{proof}

\subsection{Proof of Lemma \ref{lemma:tracial_expansion_bias}}
\begin{proof}
    Similar to the decomposition of $Q_j(z)$ involving R(z) in (\ref{second_order_expansion_Qj}), we perform the following decomposition:
    \begin{align}\label{decomposing_G1_G2_G3}
        R(z) - Q(z) &= \sum_{j=1}^n \psi_jR(x_jx_j^*-n^{-1}\textbf{A}_j)Q_j\\ \notag
        &+ \sum_{j=1}^n (\beta_j - \psi_j)R(x_jx_j^* - n^{-1}\textbf{A}_j)Q_j\\ \notag
        &+ \frac{1}{n}\sum_{j=1}^n \psi_jR\textbf{A}_j(Q_j-Q) := G_1(z) + G_2(z) + G_3(z).
    \end{align}

By arguments similar to those in Lemma \ref{lemma:tracial_bounds} and using the fact that $\|R(z)\|$ is uniformly bounded across $\mathcal{C}_n$, we have 
$$\mathbb{E}\operatorname{tr}G_2(z)A_rQM = O(n^\frac{1}{2}), \quad \mathbb{E}\operatorname{tr}G_3(z)A_rQM = O(1).$$
Observing that $\mathbb{E}\operatorname{tr}RA_jG_1 = 0$, we therefore have
\begin{align*}
    \frac{1}{n}|\mathbb{E}\operatorname{tr}(R\textbf{A}_jQ - R\textbf{A}_jQ)| = o(1).
\end{align*}
It also follows that:
\begin{align*}
    \mathbb{E}\operatorname{tr}(QA_rQM) = \operatorname{tr}RA_rRM - \operatorname{tr}G_1A_rQM + a(z),
\end{align*}
where $\mathbb{E}|a(z)| = O(n^\frac{1}{2}).$ Following the steps in (\ref{object_interest}), the second term simplifies to
\begin{align*}
    \operatorname{tr}G_1A_rQM &= \sum_{j=1}^n \psi_j \operatorname{tr}R(x_jx_j^* - n^{-1}\textbf{A}_j)Q_jA_rQM\\
    &= \sum_{j=1}^n \psi_j \operatorname{tr}R(x_jx_j^* - n^{-1}\operatorname{tr}\textbf{A}_j)Q_jA_r(Q-Q_j + Q_j)M\\
    &= -\sum_{j=1}^n \psi_j^2 x_j^*Q_jA_rQ_jx_j\, x_j^*Q_jMRx_j \\
    &+ \sum_{j=1}^n \psi_j\bigg(x_j^*Q_jA_rQ_jMRx_j  - \frac{1}{n}\operatorname{tr}R\textbf{A}_jQ_jA_rQ_jM\bigg)\\
    &= -\frac{1}{n^2}\sum_{j=1}^n\psi_j^2\bigl(\operatorname{tr}Q_jA_rQ_j\textbf{A}_j\bigr) \,\bigl(\operatorname{tr} Q_jMR\textbf{A}_j\bigr) + a(z).
\end{align*}
Therefore, 
\begin{align*}
    \operatorname{tr}QA_rQM &= \operatorname{tr}RA_rRM + \frac{1}{n^2}\sum_{j=1}^n\psi_j^2 \bigl(\mathbb{E}\operatorname{tr}Q\textbf{A}_jQA_r\bigr)\bigl(\mathbb{E}\operatorname{tr}RMR\textbf{A}_j\bigr) + a(z)\\
    &=\operatorname{tr}RA_rRM + \sum_{s,t=1}^K\bigg(\frac{1}{n^2}\sum_{j=1}^n\psi_j^2\theta_{sj}\theta_{tj}\bigg) \bigl(\mathbb{E}\operatorname{tr}QA_sQA_r\bigr)\bigl(\mathbb{E}\operatorname{tr}RMRA_t\bigr) + a(z).
\end{align*}
\end{proof}

\begin{lemma}\label{lemma:norm_control}
    Let $A_n \in \mathbb{C}^{K \times K}$ be a sequence of random matrices such that $$\mathbb{E}|a_{ij}|^2 = O(n^{-\frac{1}{2}}).$$ Then, $\mathbb{E}\|A_n\| = O(n^{-\frac{1}{2}})$.
\end{lemma}

\begin{proof}
    The proof is clear from the following inequality:
    \begin{align*}
        \mathbb{E}\|A_n\|^2 \leq \mathbb{E}\|A_n\|_F^2 \leq CK^2/n.
    \end{align*}
\end{proof}

\begin{lemma}\label{lemma:Billingsley_prob8_page17}
    Suppose $\{f_\gamma: \gamma \in \Gamma\}$ is a bounded equicontinuous class of functions with a compact domain $D$. If $F_n \xrightarrow{d} F$, then
    $$\sup_{\gamma \in \Gamma}\bigg|\int f_\gamma(x)d\{F_n(x) - F(x)\}\bigg| \rightarrow 0.$$
\end{lemma}

\begin{proof}
    Let $\epsilon>0$ be arbitrary. By equicontinuity, there exists $\delta>0$ such that $d(x, y) < \delta$ implies that 
    $$\sup_{\gamma \in \Gamma}|f_\gamma(x) - f_\gamma(y)| < \epsilon.$$
By compactness, there exists a finite $2\delta-$net of $D$ to be denoted by $\{B_i\}_{i=1}^n$. Let $\{p_i\}_{i=1}^n$ be arbitrary points in $B_i$. For arbitrary $x \in D$, there exists $i(x) \in [n]$ such that $d(x, p_{i(x)}) < 2\delta$. Then for any $\gamma \in \Gamma$, we observe that
\begin{align*}
    &\bigg|\int f_\gamma(x)dF_n(x) - \int f_\gamma(x)dF(x)\bigg|\\
    \leq& \int |f_\gamma(x)-f_\gamma(p_{i(x)})|dF_n(x) +  \int |f_\gamma(p_{i(x)})|d\{F_n(x)-F(x)\} +  \int |f_\gamma(x)-f_\gamma(p_{i(x)})|dF(x)\\
    \leq&\, 2\epsilon + \max_{i \in [n]} \int |f_\gamma(p_i(x))| d\{F_n(x) - F(x)\}.
\end{align*}
The second term goes to $0$ because we are looking at the maximum over a finite number of quantities, each of which converges to $0$. As a result, the RHS is free of $\gamma \in \Gamma$.
\end{proof}

\end{document}